\newif\ifpdf
\ifx\pdfoutput\undefined
\pdffalse 
\else
\pdfoutput=1 
\pdftrue \fi

\newif\iffinal
\finalfalse	
\finaltrue 

\documentclass[reqno,twoside,11pt]{amsart}
\iffinal\else\usepackage[notref,notcite]{showkeys}\fi
\iffinal\else\IfFileExists{pdfsync.sty}{\usepackage{pdfsync}}{}\fi
\usepackage{cite}
\usepackage{amsmath}
\usepackage{amsfonts}
\usepackage{amssymb}
\usepackage{epsfig}
\usepackage{verbatim}

\usepackage{mathrsfs}
\usepackage{mathpazo}
\usepackage{mathabx}

\usepackage{amsthm,mathtools}
\usepackage{xcolor}
\usepackage{ esint }

\DeclareFontFamily{OT1}{eusb}{} \DeclareFontShape{OT1}{eusb}{m}{n} {<5> <6> <7> <8> <9> <10> <11> <12> <14.4> eusb10}{}
\DeclareMathAlphabet{\eusb}{OT1}{eusb}{m}{n}

\DeclareFontFamily{OT1}{eusm}{} \DeclareFontShape{OT1}{eusm}{m}{n} {<5> <6> <7> <8> <9> <10> <11> <12> <14.4> eusm10}{}
\DeclareMathAlphabet{\eusm}{OT1}{eusm}{m}{n}

\DeclareFontFamily{OT1}{eufm}{} \DeclareFontShape{OT1}{eufm}{m}{n} {<5> <6> <7> <8> <9> <10> <11> <12> <14.4> eufm10}{}
\DeclareMathAlphabet{\mathfrak}{OT1}{eufm}{m}{n}

\DeclareFontFamily{OT1}{fraktura}{}
\DeclareFontShape{OT1}{fraktura}{m}{n} {<5> <6> <7> <8> <9> <10> <11> <12> <13> <14.4> [1.1] eufm10}{}
\DeclareMathAlphabet{\fraktura}{OT1}{fraktura}{m}{n}

\DeclareFontFamily{OT1}{cmfi}{} \DeclareFontShape{OT1}{cmfi}{m}{n} {<5> <6> <7> <8> <9> <10> <11> <12> <13> <14.4> [0.9] cmfi10}{}
\DeclareMathAlphabet{\cmfi}{OT1}{cmfi}{b}{n}

\DeclareFontFamily{OT1}{cmss}{} \DeclareFontShape{OT1}{cmss}{m}{n} {<5> <6> <7> <8> <9> <10> <11> <12> <13> <14.4> cmss10}{}
\DeclareMathAlphabet{\cmss}{OT1}{cmss}{m}{n}

\newtheoremstyle{thm}{1.5ex}{1.5ex}{\itshape\rmfamily}{} {\bfseries\rmfamily}{}{2ex}{}

\newtheoremstyle{def}{1.5ex}{1.5ex}{\rmfamily\sl}{} {\bfseries\rmfamily}{}{2ex}{}

\newtheoremstyle{rem}{1.3ex}{1.3ex}{\rmfamily}{} {\itshape}
{} {1.5ex}{}

\newenvironment{proofsect}[1] {\vskip0.1cm\noindent{\rmfamily\itshape#1.}}{\qed\vspace{0.15cm}}

\theoremstyle{thm}
\newtheorem{theorem}{Theorem}[section]
\newtheorem{lemma}[theorem]{Lemma}
\newtheorem{proposition}[theorem]{Proposition}

\newtheorem*{Main Theorem}{Main Theorem.}
\newtheorem{corollary}[theorem]{Corollary}

\newtheorem{assumption}[theorem]{Assumption}

\theoremstyle{def}

\theoremstyle{rem}
\newtheorem{remark}[theorem]{{Remark}}

\numberwithin{equation}{section}

\renewcommand{\section}{\secdef\sct\sect}
\newcommand{\sct}[2][default]{\refstepcounter{section}
\addcontentsline{toc}{section}
{{\tocsection {}{\thesection}{\!\!\!\!#1\dotfill}}{}}
\vspace{0.7cm}
\centerline{ 
\scshape\arabic{section}.\ #1} \nopagebreak \vspace{0.2cm}}
\newcommand{\sect}[1]{
\vspace{0.4cm} \centerline{\large\scshape\rmfamily #1}
\vspace{0.2cm}}

\renewcommand{\subsection}{\secdef\subsct\sbsect}
\newcommand{\subsct}[2][default]{\refstepcounter{subsection}
\addcontentsline{toc}{subsection}
{{\tocsection{\!\!}{\hspace{1.2em}\thesubsection}{\!\!\!\!#1\dotfill}}{}}
\vspace{0.45\baselineskip} {\flushleft\bf
\arabic{section}.\arabic{subsection}~\bf #1.~}
\\*[3mm]\noindent
\nopagebreak}
\newcommand{\sbsect}[1]{\vspace{0.1cm}\noindent
\textbf{#1.~}\vspace{0.1cm}}

\renewcommand{\subsubsection}{%
\secdef \subsubsect\sbsbsect}
\newcommand{\subsubsect}[2][default]{%
\refstepcounter{subsubsection} 
\addcontentsline{toc}{subsubsection}{{\tocsection{\!\!}
{\hspace{3.05em}\thesubsubsection}{\!\!\!\!#1\dotfill}}{}}
\nopagebreak
\vspace{0.15\baselineskip} \nopagebreak {\flushleft\rmfamily
\itshape\arabic{section}.\arabic{subsection}.\arabic{subsubsection}
\ \rmfamily #1\/.}\ }
\newcommand{\sbsbsect}[1]{\vspace{0.1cm}\noindent
\rmfamily \itshape
\arabic{section}.\arabic{subsection}.\arabic{subsubsection} \
\sffamily #1\/.\ }

\iffinal

\else

\fi

\renewcommand{\caption}[1]{%
\vglue0.5cm
\refstepcounter{figure}
\begin{minipage}{0.9\textwidth}\small {\sc Figure~\thefigure. }#1\end{minipage}}

\newcommand{\diam}{\operatorname{diam}}

\newcommand{\textd}{\text{\rm d}\mkern0.5mu}
\newcommand{\texti}{\text{\rm  i}\mkern0.7mu}
\newcommand{\texte}{\text{\rm  e}\mkern0.7mu}

\newcommand{\BB}{\mathcal B}

\newcommand{\FF}{\mathcal F}

\newcommand{\LL}{\mathcal L}

\newcommand{\NN}{\mathcal N}

\newcommand{\N}{\mathbb N}

\newcommand{\R}{\mathbb R}

\newcommand{\T}{\mathbb T}

\newcommand{\Z}{\mathbb Z}

\newcommand{\twoeqref}[2]{(\ref{#1}--\ref{#2})}

\def\myffrac#1#2 in #3{\raise 2.6pt\hbox{$#3 #1$}\mkern-1.5mu\raise 0.8pt\hbox{$#3/$}\mkern-1.1mu\lower 1.5pt\hbox{$#3 #2$}}
\newcommand{\ffrac}[2]{\mathchoice%
	{\myffrac{#1}{#2} in \scriptstyle}
	{\myffrac{#1}{#2} in \scriptstyle}
	{\myffrac{#1}{#2} in \scriptscriptstyle}
	{\myffrac{#1}{#2} in \scriptscriptstyle}
}

\newcommand{\cc}{\text{\rm c}}

\newcommand{\wt}{\widetilde}
\newcommand{\wh}{\widehat}

\newcommand{\frakp}{\fraktura p}
\newcommand{\frakc}{\fraktura c}

\newcommand{\eMB}{\normalcolor}

\usepackage{ulem}
\newcommand{\stkout}[1]{\ifmmode\text{\sout{\ensuremath{#1}}}\else\sout{#1}\fi}

\newcommand{\bm}{\boldsymbol}

\newcommand{\eHY}{\normalcolor}

\newcommand{\frakd}{\fraktura d}
\newcommand{\frake}{\fraktura e}
\newcommand{\fraks}{\fraktura s}
\newcommand{\fraku}{\fraktura u}
\newcommand{\frakv}{\fraktura v}
\newcommand{\frakw}{\fraktura w}

\begin{document}

\vglue-2mm

\title[Phase diagram of DG~model \hfill]{
Phase transition and critical behavior in hierarchical\\integer-valued Gaussian and Coulomb gas models}
\author[\hfill M.~Biskup, H.~Huang]
{Marek~Biskup \,{\tiny and}\, Haiyu Huang}
\thanks{\hglue-4.5mm\fontsize{9.6}{9.6}\selectfont\copyright\,\textrm{2025}\ \ \textrm{M.~Biskup, H.~Huang. Reproduction for non-commercial use is permitted.}}
\maketitle

\vglue-5mm
\centerline{\it Department of Mathematics, UCLA, Los  Angeles, California, USA}



\begin{abstract}
 Given a square box~$\Lambda_n\subseteq\Z^2$ of side length~$L^n$ with~$L,n>1$, we study hierarchical random fields $\{\phi_x\colon x\in\Lambda_n\}$ with law proportional to $\texte^{\frac12\beta(\phi,\Delta_n\phi)}\prod_{x\in\Lambda_n}\nu(\textd\phi_x)$, where~$\beta>0$ is the inverse temperature,~$\Delta_n$ is a hierarchical Laplacian on~$\Lambda_n$, and $\nu$ is a non-degenerate $1$-periodic measure on~$\R$. Our setting includes the integer-valued Gaussian field (a.k.a.\  DG~model  or Villain Coulomb gas) and the sine-Gordon model. Relying on renormalization group analysis we derive sharp asymptotic formulas, in the limit as $n\to\infty$, for the covariance $\langle\phi_x\phi_y\rangle$ and the fractional charge $\langle\texte^{2\pi\texti\alpha(\phi_x-\phi_y)}\rangle$ in the subcritical $\beta<\beta_\cc:=\pi^2/\log L$, critical $\beta=\beta_\cc$ and slightly supercritical $\beta>\beta_\cc$ regimes. The field exhibits logarithmic correlations throughout albeit with a distinct $\beta$-dependence of both the covariance scale and the fractional-charge exponents in the sub/supercritical regimes. Explicit logarithmic corrections appear at the critical point.
\end{abstract}


\section{Introduction and results}
\label{sec-1}\noindent
\vglue-3mm
\subsection{The model and assumptions}
The aim of this paper is to study a class of random fields on~$\Z^2$ with periodically modulated values. The general setting of these models  is  as follows: Fix an integer~$L\ge2$ and, for each integer $n\ge1$, let $\Lambda_n:=\{0,\dots,L^{n}-1\}^2$ be a box of  side-length  $L^{n}$ in~$\Z^2$. Then consider a family $\{\phi_x\colon x\in\Lambda_n\}$ of real-valued random variables with joint law
\begin{equation}
\label{E:1.1}
P_{n,\beta}(\textd\phi):=\frac1{Z_n(\beta)}\texte^{\frac12\beta(\phi,\Delta_n\phi)}\prod_{x\in\Lambda_n}\nu(\textd\phi_x),
\end{equation}
where $\beta>0$ is the inverse temperature, $Z_n(\beta)$ is a normalization constant,  $(\cdot,\cdot)$ denotes the canonical inner product in $\ell^2(\Lambda_n)$ and $\Delta_n$ is a Laplacian or, in probabilistic terms, the generator of a Markov chain on~$\Lambda_n$. The periodic modulation enters through $\nu$ which is assumed to be a $1$-periodic locally-finite  positive Borel measure on~$\R$.

Throughout we focus on hierarchical models, for which the Markov chain defined by~$\Delta_n$ jumps from~$x$ to~$y$ at a rate that depends only on the coefficients in base-$L$ expansion of the coordinates of $x$ and~$y$. To state this precisely, set $b:=L^2$ and identify~$\Lambda_n$ with the set of sequences $(x_1,\dots,x_n)\in\{0,\dots,b-1\}^n$. For $x=(x_1,\dots,x_n)\in\Lambda_n$, let~$\BB_k(x)$ denote the set of~$y=(y_1,\dots,y_n)$ such that~$y_i=x_i$ for all $i=1,\dots,n-k$. Then take~$\Delta_n$ to be the operator that acts on~$f\colon\Lambda_n\to\R$ as
\begin{equation}
\label{E:1.2}
\Delta_n f(x) :=- \frakc_{n+1}\, f(x)+\sum_{k=1}^n \,
\sum_{y\in \BB_k(x)}\frakc_{k}\,\bigl[\,f(y) - f(x)\bigr],
\end{equation}
where~$\{\frakc_k\}_{k=1}^{n+1}$ is a sequence of positive numbers subject to specific decay conditions; see Assumption~\ref{ass-1} below. 
The positivity ensures that the hierarchical Laplacian $\Delta_n$ is strictly negative definite on~$\ell^2(\Lambda_n)$  which  is needed for~$Z_n(\beta)$ to be finite.

A well-studied example of above field is the hierarchical \textit{Gaussian Free Field} (GFF) for which~$\nu$ is the Lebesgue measure on~$\R$. A primary example of our interest in this work is the hierarchical \textit{integer-valued Gaussian} model, a.k.a.\ \textit{DG~model}, for which~$\nu$ is the counting measure on~$\Z$.
The two models are interpolated by a continuous family of \textit{sine-Gordon} models defined by
\begin{equation}
\label{E:1.4}
\nu(\textd\phi):=\texte^{-\kappa[1-\cos( 2\pi \phi)]}\textd\phi,
\end{equation}
where~$\kappa>0$ is a parameter. Indeed, the GFF corresponds to $\kappa=0$ while the DG~model arises via the weak limit as  $\kappa\to\infty$ under the scaling of~$\nu$ by $(2\pi\kappa)^{ 1/2}$. 

The models \eqref{E:1.1} turn out to be dual to Coulomb gas systems whenever the Fourier coefficients of~$\nu$ are non-negative. A remarkable fact is that two-dimensional Cou\-lomb gas models, and thus also our fields, undergo a \textit{BKT phase transition} at some~$\beta_\cc$ (named after Berezinskii~\cite{Berezinskii}, Kosterlitz and Thouless~\cite{KT}) as soon as~$\nu$ is distinct from the Lebesgue measure; see Section~\ref{sec-2.1} for more discussion. Various aspects of this transition have previously been addressed in hierarchical setting (e.g., by Benfatto, Gallavotti and Nicol\`o~\cite{BGN}, Marchetti and Perez~\cite{MP}, Benfatto and Renn~\cite{BR}, Guidi and Marchet\-ti~\cite{GM}) albeit subject to limitations that generally exclude the DG~model. Our aim here is to provide a robust treatment of the transition, and establish heretofore uncontrolled aspects of the critical and near-critical behavior.

Similarly to references \cite{BGN,MP,BR,GM}, our analysis  relies  on the  renormalization group  technique whose implementation requires  some  regularity of the coefficients $\{\frakc_k\}_{k=1}^{n+1}$. We collect these requirements in:

\begin{assumption}
\label{ass-1}
There exists a positive sequence $\{\frakd_k\}_{k\ge0}$ satisfying
$\sum_{k\ge0}\frakd_k<\infty$ such that, for each~$n\ge1$, the sequence $\{\frakc_k\}_{k=1}^{n+1}$ takes the form
\begin{equation}
\label{E:1.5u}
\frakc_{n+1}=\biggl(\,\sum_{j=0}^n b^j\sigma_j^2\biggr)^{-1}
\end{equation}
and
\begin{equation}
\label{E:1.6u}
\frakc_k = \frac1{b^k}\Biggl[\biggl(\,\sum_{j=0}^{k-1} b^j\sigma_j^2\biggr)^{-1}-\biggl(\,\sum_{j=0}^k b^j\sigma_j^2\biggr)^{-1}\Biggr],\quad k=1,\dots,n,
\end{equation}
for a strictly positive sequence $\{\sigma_k^2\}_{k=0}^n$ satisfying
\begin{equation}
\label{E:1.4u}
\bigl|\sigma_k^2-1\bigr|\le \frakd_{\min\{k,n-k\}},\quad k=0,\dots,n.
\end{equation}
Moreover, we have $\inf_{n\ge k\ge0}\sigma_k^2>0$.

\end{assumption}

In addition, we also need a bit of regularity of the measure~$\nu$:

\begin{assumption}
\label{ass-2}
~$\nu$ is a  $1$-periodic Radon measure on~$\R$ whose Fourier coefficients defined by $a(q):=\int_{[0,1)}\texte^{ -  2\pi\texti q z}\nu(\textd  z)$ are real-valued,  strictly positive and satisfy
\begin{equation}
\label{E:1.6w}
a(-q)=a(q),\quad q\in\Z,
\end{equation}
along with
\begin{equation}
\label{E:1.7a}
\sup_{q\ge0}\frac{a(q+1)}{a(q)}<\infty.
\end{equation}
In particular, $\nu$ is reflection symmetric and $1$-periodic but not~$1/p$-periodic for any~$p\ge2$.
\end{assumption}

Note that, while the analysis by way of the renormalization technique is easiest when $\{\sigma_k^2\}_{k=0}^n$ are all equal to a positive constant (which we take to be~$1$), permitting more general coefficients allows us to remain flexible as to what specific operator we take to be the  hierarchical Laplacian~$\Delta_n$; see Remark~\ref{rem-hierarchy}. The formulation using $\{\frakd_k\}_{k\ge0}$ is done to ensure uniformity. Observe also that \eqref{E:1.4u} along with $\frakd_k\to0$ imply that~$\frakc_k$, for both~$k$ and~$n-k$ large, decays proportionally to~$k\mapsto b^{-2k}$. The term $\frakc_{n+1}$ scales only as~$b^{-n}$ due to its role of a ``mass''; see again Remark~\ref{rem-hierarchy}.
(The same asymptotic arises if we think of~$\frakc_{n+1}$ as an aggregate killing rate $\sum_{k>n}\frakc_k b^k$ for a Markov chain on~$\Z^2$ with conductances~$\{\frakc_k\}_{k>n}$ given as in~\eqref{E:1.6u}.)

As to the conditions on measure~$\nu$, here the DG~model corresponds to $a(q)=1$ for each~$q\in\Z$ while for the sine-Gordon model \eqref{E:1.4}, we get 
\begin{equation}
\label{E:1.7c}
a(q)=\sum_{\ell=0}^\infty\frac{(\kappa/2)^{2\ell+|q|}}{(\ell+|q|)!\ell!},\quad q\in\Z.
\end{equation}
In particular, these models satisfy the conditions \twoeqref{E:1.6w}{E:1.7a}. The GFF is excluded but so is unfortunately the \textit{hard-core Coulomb gas} that corresponds to
\begin{equation}
\label{E:1.9}
\nu(\textd\phi) := [1+2\kappa\cos(2\pi\phi)]\textd\phi,
\end{equation}
where~$\kappa\in[0,1/2]$. While our conclusions (to be stated next) definitely fail for the GFF, we still expect them to apply to the model \eqref{E:1.9}.

\subsection{Covariance structure}
Our first result concerns the asymptotic covariance structure of the field. Recall that~$b$ denotes the ``branching number'' of the hierarchical model which in the description based on a box in~$\Z^2$ relates to the base scale~$L$ of~$\Lambda_n$ as $b=L^2$. The representation of elements of~$\Lambda_n$ as sequences leads to a hierarchical  metric on~$\Lambda_n$,  defined   for any two  distinct  vertices $x=(x_1,\dots,x_n)$ and~$y=(y_1,\dots,y_n)$ by
\begin{equation}
\label{E:1.9i} 
 d(x,y):= b^{\frac12(\min\{\,j\in \{0,\dots, n\}: \,x_i = y_i \,\forall i = 1,\dots, n- j\,\})},
\end{equation}
with the proviso $d(x,y):=0$ when~$x=y$. Under the embedding of~$\Lambda_n$ into~$\Z^2$ we have $d(x,y)\ge\Vert x-y\Vert_\infty$ for all~$x,y\in\Lambda_n$ with both sides comparable for generic~$x$ and~$y$.

As a consequence of the hierarchical structure of~$\Delta_n$, all of our models undergo a phase transition at the same value of the inverse temperature; namely, at
\begin{equation}
\beta_\cc:=\frac{2\pi^2}{\log b}.
\end{equation}
We will henceforth write $\langle-\rangle_{n,\beta}$ to denote expectation with respect to~$P_{n,\beta}$. Our result on  the covariance structure of~$P_{n,\beta}$ is then as follows:

\begin{theorem}[Covariance structure]
\label{thm-1}
For each~$b\ge2$ there exists a non-increasing continuous function $\sigma^2\colon\R_+\to\R_+$ satisfying
\begin{equation}
\label{E:1.10}
\sigma^2(\beta)\,\begin{cases}
=1/\beta,\qquad&\text{if }\beta\le\beta_\cc,
\\
<1/\beta,\qquad&\text{if }\beta>\beta_\cc,
\end{cases}
\end{equation}
such that the following holds for all models satisfying Assumptions~\ref{ass-1}--\ref{ass-2} with~$\{\frakd_k\}_{k\ge0}$ decaying exponentially when~$\beta>\beta_\cc$  and obeying $\sum_{j\ge1}\frakd_j\log(j)<\infty$ when~$\beta=\beta_\cc$: 
There exists $\epsilon>0$ and, for all~$\beta>0$ with $1/\beta>1/\beta_\cc-\epsilon$,
all~$n\ge1$ and all~$x,y\in\Lambda_n$,
\begin{equation}
\label{E:1.12}
\langle\phi_x\phi_y\rangle_{n,\beta}=\begin{cases}
\displaystyle \sigma^2(\beta)\log_{b^{1/2}}\Bigl(\frac{\diam(\Lambda_n)}{1+d(x,y)}\Bigr)+O(1),&\text{if }\beta\ne\beta_\cc,
\\*[4mm]
\displaystyle \frac1{\beta_\cc}\log_{b^{1/2}}\Bigl(\frac{\diam(\Lambda_n)}{1+d(x,y)}\Bigr)-\bar c\log\Bigl(\frac{\log\diam(\Lambda_n)}{\log[2+d(x,y)]}\Bigr)+ O(1),&\text{if }\beta=\beta_\cc,
\end{cases}
\end{equation}
where
\begin{equation}
\label{E:1.15ui}
\bar c:=\frac{8\pi^2}{\beta_\cc^2} \frac{b(b^3-1)}{(b-1)^3(b+1)^2}
\end{equation}
and~$O(1)$ are quantities bounded uniformly in~$n\ge1$ and~$x,y\in\Lambda_n$.
\end{theorem}

The above shows that models \eqref{E:1.1} subject to Assumptions~\ref{ass-1}--\ref{ass-2} exhibit logarithmic decay of correlations at all $\beta>0$ (with~$1/\beta> 1/\beta_\cc-\epsilon$). This makes them qualitatively similar to the GFF, for which the covariances behave exactly as in the~$\beta<\beta_\cc$ regime above. The connection to GFF at $\beta<\beta_\cc$ is very tight; indeed, in our earlier work~\cite{BHu} we showed that one can couple~$P_{n,\beta}$ to the law of GFF so closely that the two fields are within order unity at most points. 

For $\beta>\beta_\cc$, \eqref{E:1.10} shows that the overall scale of the fluctuations is strictly smaller than what GFF would give and, indeed, $P_{n,\beta}$ is far from the law of GFF both in terms of global scaling properties as well as other correlations (see Theorem~\ref{thm-2} below). A reader looking for exponential decay for~$\beta>\beta_\cc$ should note that the Laplacian \eqref{E:1.2} is long range with the matrix coefficient for the pair~$x$ and~$y$ decaying proportionally to~$d(x,y)^{-4}$ (see Section~\ref{sec-2.3}) so exponential decay is not to be expected. 

\begin{figure}[t]
\refstepcounter{obrazek}
\label{fig1}
\vglue5mm
\centerline{\includegraphics[width=5.2in]{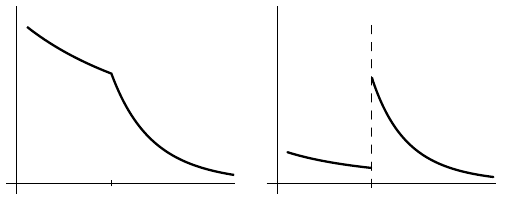}}
\vspace{.2in}
\begin{quote}
\fontsize{9}{5}\selectfont
{\sc Fig~\theobrazek.} Left: A plot of $\beta\mapsto\sigma^2(\beta)$ for~$b=2$ and $\beta\in[20,40]$, which (roughly) corresponds to~$\theta\in[0.37, 0.61]$. (Note that $\beta_\cc \approx 28.477$.) Right: A plot of the negative derivative of~$\sigma^2$. The numerical value of the jump is $0.00716$ which is in good agreement with the exact value $\approx 0.0073$ obtained from \eqref{E:1.16ui}.
\normalsize
\end{quote}
\end{figure}

The difference in the overall variance scale arises from the fact that, at $\beta>\beta_\cc$, the field ``senses'' the $1$-periodicity of~$\nu$ at all spatial scales. Technically, this is seen in  renormalization group iterations that draw the system towards a ``non-trivial'' fixed point --- meaning one that does not correspond to GFF --- rather than the ``trivial'' one as happens for~$\beta\le\beta_\cc$. The quantity~$\sigma^2(\beta)$  admits a formula, see \eqref{E:4.39},  that  makes the inequality in \eqref{E:1.10} quite apparent. We even get the asymptotic expansion
\begin{equation}
\label{E:1.16ui}
\sigma^2(\beta) = \frac1\beta-\frac{ 32\pi^4}{\beta_\cc^4}\frac{b(b^3-1)}{(b-1)^3(b+1)^2}(\beta-\beta_\cc)+O\bigl((\beta-\beta_\cc)^{ 3/2}\bigr),\quad \beta\downarrow\beta_\cc,
\end{equation}
see Remark~\ref{rem-sigma-asymp}. In particular, $\beta\mapsto\sigma^2(\beta)$ is not differentiable at~$\beta_\cc$, see Fig.~\ref{fig1}. The apparent numerical closeness of \eqref{E:1.15ui}  to the coefficient of $\beta-\beta_\cc$ in  \eqref{E:1.16ui} is not a coincidence; see Remark~\ref{rem-crit-near-crit1}  for an explanation. 

\newcounter{obrazek}

The behavior at~$\beta_\cc$ is yet different as an iterated-log correction arises in the covariance structure. This can be attributed to the fact that, while the renormalization group iterations still draw the model towards the ``trivial'' fixed point, the convergence is polynomially slow and a residue of $1$-periodicity of~$\nu$ survives to the macroscopic scale. We expect that the iterated-log correction is reflected in the extremal behavior of the field. For instance, the maximum $\max_{x\in\Lambda_n}\phi_x$, which for~$\beta<\beta_\cc$  scales exactly as that of the GFF (see \cite[Corollary~2.2]{BHu}), should have a different second-order (i.e., order $\log n$) term at~$\beta=\beta_\cc$.  Controlling the maximum at (and beyond) the critical~$\beta$ is an interesting open problem; see Remark~\ref{rem-max} for further discussion.

\subsection{Fractional charge asymptotic}
The connection of our model with Coulomb gas naturally leads us to the  so-called   fractional charge correlation $\langle \texte^{2\pi\texti \alpha(\phi_x-\phi_y)}\rangle_{n,\beta}$, where~$\alpha$ is a parameter that, due to the underlying $1$-periodicity and also interpretation as an electric charge, is taken generally real-valued; see Section~\ref{sec-2.4}. Here we get:

\begin{theorem}[Fractional charge asymptotic]
\label{thm-2}
For each~$b\ge2$ there exists a continuous function $\kappa\colon(0,\ffrac12)\times\R_+\to\R_+$ satisfying
\begin{equation}
\label{E:1.16}
\kappa(\alpha,\beta)\,\begin{cases}
=
\frac{4\beta_\cc}\beta \alpha^2,\qquad&\text{if }\beta\le\beta_\cc,
\\*[2mm]
<
\frac{4\beta_\cc}\beta \alpha^2,\qquad&\text{if }\beta>\beta_\cc,
\end{cases}
\end{equation}
such that the following is true for all models satisfying
Assumptions~\ref{ass-1}--\ref{ass-2} with~$\{\frakd_k\}_{k\ge0}$ decaying exponentially when~$\beta>\beta_\cc$  and obeying $\sum_{j\ge1}\frakd_j\log(j)<\infty$ when~$\beta=\beta_\cc$:  For all $\alpha_0\in(0,\ffrac12)$ there exists $\epsilon>0$ and, for all $\alpha\in(0,\alpha_0]$, all $\beta>0$ with~  $1/\beta > 1/\beta_\cc-\epsilon$   and all~$n\ge1$, there exists $C_n\in(0,\infty)$ satisfying~$0<\inf_{n\ge1}C_n\le\sup_{n\ge1}C_n<\infty$ such that
\begin{equation}
\label{E:1.17}
\bigl\langle \texte^{2\pi\texti \alpha(\phi_x-\phi_y)}\bigr\rangle_{n,\beta}
=\bigl[C_n+o(1)\bigr]\begin{cases}
d(x,y)^{-\kappa(\alpha,\beta)},&\text{if }\beta\ne\beta_\cc,
\\*[3mm]
d(x,y)^{-\kappa(\alpha,\beta)} [\log d(x,y)]^{\tau(\alpha)},&\text{if }\beta=\beta_\cc,
\end{cases}
\end{equation}
holds for all~$x,y\in\Lambda_n$ with~$x\ne y$, where
\begin{equation}
\label{E:tau}
\tau(\alpha):=2\frac{b^3-1}{(b-1)(b+1)^3}\biggl[\frac{b-1}{b^{1+2\alpha}-1}+\frac{b-1}{b^{1-2\alpha}-1}- 2\biggr]
\end{equation}
and~$o(1)\to0$ in the limit as~$\min\{d(x,y),\diam(\Lambda_n)/d(x,y)\}\to\infty$.
\end{theorem}

\begin{figure}[b]
\refstepcounter{obrazek}
\label{fig2}
\vglue5mm
\centerline{\includegraphics[width=4in]{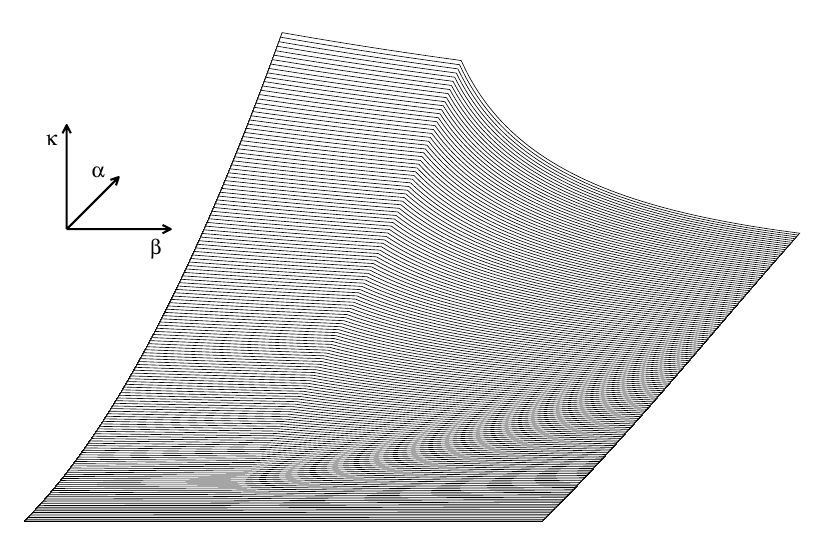}}
\vspace{.2in}
\begin{quote}
\fontsize{9}{5}\selectfont
{\sc Fig~\theobrazek.} A plot of the graph of $(\alpha,\beta)\mapsto\kappa(\alpha,\beta)$ for $b = 2$. Here $\beta$ increases to the right along the horizontal axis while~$\alpha$ increases in the direction front to back. The plot shows the portion of the graph for~$\beta\in[25,35]$ and~$\alpha\in[0,0.4]$. (Recall that $\beta_\cc \approx 28.477$ when $b = 2$.) The small-$\alpha$ asymptotic is  governed by the variance; see \eqref{E:1.21u}. 
\normalsize
\end{quote}
\end{figure}

The $n$-dependence of~$C_n$ stems from potential variability in~$n$ of the sequence~$\{\sigma_k^2\}_{k=0}^n$. The quantity~$\kappa(\alpha,\beta)$ is determined, albeit somewhat implicitly, by \eqref{Eq:5.56} which we are only able to solve numerically, see Fig.~\ref{fig2}.  As discussed in Remark~\ref{rem-kappa}, we have the asymptotic~form 
\begin{equation}
\label{E:1.20i}
\kappa(\alpha,\beta)=
\frac{4\beta_\cc}\beta \alpha^2
-
\frac4{\beta_\cc} \tau(\alpha)(\beta-\beta_\cc)+O\bigl((\beta-\beta_\cc)^2\bigr),\quad\beta\downarrow\beta_\cc,
\end{equation}
where~$\tau(\alpha)$ is as in \eqref{E:tau}. (As~$\kappa(\alpha,\beta)\ge0$, the fact that~$\tau(\alpha)$ diverges as~$|\alpha|$ increases to~$\ffrac12$ only attests the lack of uniformity.) Again, the connection between the critical and near-critical asymptotic is not a coincidence; see Remark~\ref{rem-kappa}.

Since \eqref{E:1.17} is, in effect, a calculation of the characteristic function of~$\phi_x-\phi_y$, a comparison with \eqref{E:1.12} yields the asymptotic form
\begin{equation}
\label{E:1.21u}
\kappa(\alpha,\beta) = 4\beta_\cc\,\sigma^2(\beta)\alpha^2 +o(\alpha^2),\quad \alpha\to0.
\end{equation}
The advantage of \eqref{E:1.17} is that we learn a bit more about the distribution than just variance. Indeed, \eqref{E:1.17} implies that $\phi_x-\phi_y$ normalized by its standard deviation is Gaussian in the limit $n\to\infty$ and~$d(x,y)\to\infty$. We expect that the martingale calculations underpinning the proof of Theorem~\ref{thm-1} would allow us to quantify the corrections beyond the leading order which we expect to be non-Gaussian once~$\beta\ge\beta_\cc$. 

In the language of Coulomb gas models, $\log\langle \texte^{2\pi\texti \alpha(\phi_x-\phi_y)}\rangle_{n,\beta}$ represents the energetic cost of inserting a charge~$\alpha$ at~$x$ and a charge~$-\alpha$ at~$y$ into a system of integer-valued charges at thermal equilibrium. Inserting just a single charge~$\alpha$ at~$x$ has energetic cost $\log\langle \texte^{2\pi\texti \alpha\phi_x}\rangle_{n,\beta}$ for which our proof similarly shows
\begin{equation}
\label{E:1.21}
\bigl\langle \texte^{2\pi\texti \alpha\phi_x}\bigr\rangle_{n,\beta}
=\bigl[C_n'+o(1)\bigr]\begin{cases}
N^{-\kappa(\alpha,\beta)/2},&\text{if }\beta\ne\beta_\cc,
\\*[3mm]
N^{-\kappa(\alpha,\beta)/2}[\log N]^{\tau(\alpha)/2},&\text{if }\beta=\beta_\cc,
\end{cases}
\end{equation}
where we abbreviated $N:=\diam(\Lambda_n)$ and where~$o(1)\to0$ as~$N\to\infty$; see Remark~\ref{rem-single-charge}. The drop in the value of $\kappa(\alpha,\beta)$ marked by the inequality in \eqref{E:1.16} is indicative of a \textit{charge screening} taking place above~$\beta_\cc$ which (unlike for the lattice model) is only partial due to the long-range structure of the hierarchical Laplacian. See again Section~\ref{sec-2.4}.

We note that some aspects of the above result are already known. For instance, the subcritical part of \eqref{E:1.17} appears as an upper bound in Marchetti and Perez~\cite[Theorem~4.3]{MP}, albeit assuming that~$\nu$ is suitably close to the Lebesgue measure when~$\beta$ is close to~$\beta_\cc$. For~$\beta\gtrsim\beta_\cc$, \cite[Theorem~5.1]{MP} shows existence and stability of a non-trivial  renormalization group fixed point and, for the model  with~$\nu$ corresponding to the  fixed point, gives the leading order expansion of the fractional charge  exponent  as~$\beta\downarrow\beta_\cc$,  albeit less explicitly than  \eqref{E:1.20i}. (The paper~\cite{MP} works in the language of Coulomb gasses, so translations described in Section~\ref{sec-2.4} are needed to identify their conclusions with ours.)

Another relevant paper is  that  by Benfatto and Renn~\cite{BR} who (while working in our framework) established existence of a non-trivial renormalization fixed point for $\beta\gtrsim\beta_\cc$. Then, for the model with~$\nu$ corresponding to this fixed point, they studied the \textit{integer-charge correlations}; namely, truncated correlations of $1$-periodic functions~$f$ of the field, for which they proved the asymptotic
\begin{equation}
\bigl\langle \,f(\phi_x)f(\phi_y)\bigr\rangle_{n,\beta}-\bigl\langle\, f(\phi_x)\bigr\rangle_{n,\beta}\bigl\langle \,f(\phi_y)\bigr\rangle_{n,\beta}  \asymp  d(x,y)^{-2}
\end{equation}
as~$d(x,y)\to\infty$ regardless of~$\beta\gtrsim\beta_\cc$. 
This coincides with the behavior of the massive hierarchical~GFF. It will be of interest to find an argument that proves the same for more general initial~$\nu$.

\subsection{Summary and main ideas}
Theorems~\ref{thm-1} and~\ref{thm-2} capture the character of the phase transition in $\Z$-modulated hierarchical fields by way of asymptotic form of two important correlation functions. The main novelty is uniformity in the underlying model, and thus \textit{universality}, which we achieve (in Theorems~\ref{thm-3}--\ref{thm-5}) by following directly the exponential of the renormalized potentials, rather than the potentials themselves. This avoids arguments based on linearization, whose accuracy deteriorates close to the critical point, and thus also the restrictions on the model assumed in earlier work. Our method works uniformly all the way up to and even slightly beyond the critical point revealing heretofore unattended aspects of the critical behavior.

Our conclusions for the subcritical and critical regimes apply solely under Assumptions~\ref{ass-1} and~\ref{ass-2}. In the supercritical regime we restrict to $\beta-\beta_\cc$ small, but we think of this as a mere technicality whose purpose is to keep (already very long) proof of Theorem~\ref{thm-5} to a manageable length. A more serious restriction (for all~$\beta$) comes in the assumption that the Fourier coefficients of~$\nu$ obey Assumption~\ref{ass-2}. This is natural for the connection with Coulomb gas but not necessarily so for the field itself. We take this as a price to pay for the precision of our conclusions.

In tracking the renormalization group ``flow'' of the effective potentials we rely on Fourier representation and thus effectively work in the Coulomb gas picture.  A key observation, stated in Lemma~\ref{lemma-3.5} which itself draws on~\cite[Lemma~4.2]{BHu}, is that the iterations preserve the structure in Assumption~\ref{ass-2} and, in fact, improve the estimate on the ratios in~\eqref{E:1.7a}. For~$\beta\le\beta_\cc$ this leads to a full asymptotic analysis while, for $\beta>\beta_\cc$, we at least eventually dominate the ratios by a quantity of order~$\sqrt{\beta-\beta_\cc}$. Assuming that to be small, a suitable fixed-point argument extracts the desired asymptotic behavior.

The proofs of Theorems~\ref{thm-1} and~\ref{thm-2} themselves hinge on the observation that the Gibbs measure \eqref{E:1.1} with hierarchical~$\Delta_n$ can be viewed as the law of a tree-indexed Markov chain after~$n$ steps. The transition probabilities of this (time-inhomogeneous) chain are simple functions of the effective potentials, see~\eqref{E:1.23}, and so one can extract a good amount of information about the chain just from the renormalization group ``flow.'' The details unfortunately still require lengthy calculations.

\subsection{Outline}
The remainder of this paper is organized as follows. First, in Section~\ref{sec-2}, we discuss the broader context of the above models while providing additional (or missing) details for various remarks made in the text above. In Section~\ref{sec-3} we then introduce the renormalization group approach and state the key convergence theorems; see Theorems~\ref{thm-3}--\ref{thm-5} in Section~\ref{sec-3.2}. Sections~\ref{sec-4} and~\ref{sec-5} are devoted to the proofs of our main results (namely, Theorems~\ref{thm-1} and~\ref{thm-2}) based on these convergence theorems. The final section (Section~\ref{sec-6}) supplies the proof of Theorem~\ref{thm-5} on supercritical renormali\-za\-tion-group flow which, unlike Theorems~\ref{thm-3}--\ref{thm-4}, could not be efficiently reduced to estimates proved in our previous work~\cite{BHu}.

\section{Connections and references}
\label{sec-2}\noindent
We proceed to discuss the broader context of our work; specifically, connections to lattice interface models, the BKT transition, Coulomb gas systems, and hierarchical models. This will also give us the opportunity to cite additional relevant literature.

\subsection{Lattice interface models}
\label{sec-2.1}\noindent
The Gibbsian distributions of the kind \eqref{E:1.1} arise as models of fluctuating interfaces in statistical mechanics, albeit with the ``harmonic'' energy term~$\frac12(\phi,-\Delta_n\phi)$ often generalized to the ``anharmonic'' expression of the form
\begin{equation}
\label{E:2.1}
\frac12\sum_{x,y\in\Lambda_n}V_{x,y}(\phi_x-\phi_y)
\end{equation}
for some collection of potentials $\{V_{xy}\colon x,y\in\Lambda_n\}$ --- typically, convex or close to convex, translation invariant and decaying sufficiently fast with~$|x-y|$; see e.g., Velenik~\cite{Velenik}, Funaki~\cite{Funaki} or Sheffield~\cite{Sheffield}.  In this language our  setting corresponds to
\begin{equation}
\label{E:2.2a}
V_{x,y}(\eta):=\frakc(x,y)\eta^2
\end{equation}
 for a collection ~$\{\frakc(x,y)=\frakc(y,x)\colon x,y\in\Lambda_n\}$ of non-negative quantities called conductances, due to a natural connection of this problem to resistor-network theory (see, e.g., Biskup~\cite{B-RCM} for a review). 
 
There are two canonical choices for the ``single-spin'' measure~$\nu$: the  Lebesgue measure on~$\R$ and the counting measure on~$\Z$. In the  former  case, the field corresponding to \eqref{E:2.2a} is the \textit{Gaussian Free Field} (GFF) associated with the generator
\begin{equation}
\label{E:2.2b}
\LL f(x):=\sum_{y\in\Lambda_n}\frakc(x,y)[f(y)-f(x)]
\end{equation}
of a Markov chain defined by the conductances $\{\frakc(x,y)\colon x,y\in\Lambda_n\}$. Here, often but not always, $\frakc(x,y)=1$ when~$x$ and~$y$ neighbors and zero otherwise.
 
The GFF is special among models \eqref{E:2.1} for the fact that many relevant quantities are explicitly computable. A continuum version of the GFF also arises as the limit process at large spatial scales for many of the above models. This was first shown for models with uniformly strictly convex potentials by Naddaf and Spencer~\cite{NS} and Giacomin, Olla and Spohn~\cite{GOS} and later extended to various cases beyond; e.g., Biskup and Spohn~\cite{BS11}, Brydges and Spencer~\cite{BS12}, Cotar, Deuschel and M\"uller~\cite{CDM}, Ye~\cite{Ye}, Adams, Buchholtz, Koteck\'y and M\"uller~\cite{AKM,ABKM}, Dario~\cite{Dario,Dario2} and Armstrong and Wu~\cite{AW}.

The integer-valued models (i.e., for~$\nu$ being the counting measure on~$\Z$) exhibit richer behavior and are thus less well understood. One clear distinction is that any perturbation of a ground state costs a uniformly positive amount of energy. A Peierls-type argument then shows that, for~$\beta$ very large, a sample from the corresponding Gibbs measure deviates from a ground state configuration only by localized perturbations whose density decreases exponentially with their size. In particular, two-point correlations decay exponentially and we have
\begin{equation}
\label{E:2.2}
\sup_{n\ge1}\sup_{x,y\in\Lambda_n}\langle (\phi_x-\phi_y)^2\rangle_{n,\beta}<\infty,
\end{equation}
for all~$\beta$ large. 

As it turns out, for $\Z$-valued fields over~$\Z^d$ with~$d\ge3$, the salient part of the previous conclusion is not limited to large~$\beta$. Indeed, the interface is  expected to be  \textit{localized} in the sense \eqref{E:2.2} for all $\beta>0$; see, e.g., Bricmont, Fontaine, and Lebowitz~\cite{BFL} for a proof for the SOS model (where~$V_{x,y}(\eta):=|\eta|$ for nearest neighbors and zero otherwise). On the other hand, in spatial dimension $d=1$ the interface is always \textit{delocalized} in the sense that the $\lim_{n\to\infty}\langle(\phi_x-\phi_y)^2\rangle_{n,\beta}$ grows linearly with~$|x-y|$.

\subsection{Roughening transition for 2D interfaces} 
The behavior of integer-valued models and even just \textit{$\Z$-modulated} ones, for which~$\nu$ is a $1$-periodic measure, in spatial dimension~$d=2$ is special and has been the source of much interest and effort of mathematical physicists and probabilists alike. Indeed, here one expects both types of behavior to arise depending on the value of~$\beta$. Specifically, localization in the sense \eqref{E:2.2} should occur for~$\beta>\beta_\cc$ and delocalization for~$\beta<\beta_\cc$, where~$\beta_\cc$ is a positive and finite critical value. The phase transition at~$\beta_\cc$ is referred to as \textit{roughening}; see e.g.~\cite{BFL} for a discussion of this phenomenon.

The roughening transition bears a close connection to another remarkable feature of two-dimensional spin models; namely, the \textit{Berezinskii-Kosterlitz-Thouless} (BKT) phase transition predicted independently by Berezinskii~\cite{Berezinskii} and Kos\-terlitz and Thouless~\cite{KT} for, e.g., the XY-model on~$\Z^2$. A common point is a power-law decay of correlations on one side of~$\beta_\cc$ in contrast to an exponential decay on the other side. In the XY-model the power-law decay occurs at~$\beta$ large as a ``residue'' of long-range order precluded by the Mermin-Wagner phenomenon. In $\Z$-modulated interface models the power-law decay takes place at~$\beta$ small where it reflects on the discrete nature of the fields being washed out at large spatial scales.

The first mathematical treatment of a BKT phase transition was achieved by Fr\"ohlich and Spencer~\cite{FS} who proved that, in the DG~model as well as sine-Gordon and other models of this type, the fractional charge correlations,
\begin{equation}
\label{E:1.5w}
x,y\mapsto\langle\texte^{2\pi\texti\alpha(\phi_x-\phi_y)}\rangle_{n,\beta}
\end{equation}
with~$\alpha$ small, exhibit power-law decay in $|x-y|$ when~$n\gg|x-y|\gg1$ at high temperatures;  i.e., for $\beta$ small. (There is no decay when~$\beta$ is large.) The argument of~\cite{FS} was later extended throughout the ``asymptotic subcritical regime'' by Marchetti and Klein~\cite{MK} although this is not the same as controlling the model up to  the conjectural  critical value~$\beta_\cc$. Alternative presentations appeared in Ph.D. thesis  of  Braga~\cite{Braga} and a recent paper by Kharash and Peled~\cite{KP}.

Fr\"ohlich and Spencer's result (see \cite[Theorem~1.1]{KP}) implies that~$\phi$ is loga\-rithmically-correlated at small~$\beta$ while the Peierls argument shows that it  exhibits exponential decay of correlations at large~$\beta$.
A different point of view has been pursued by Aizenman, Harel, Peled and Shapiro~\cite{AHPS} (drawing on Lammers~\cite{Lammers}) who focus on the asymptotic properties of the variance function
\begin{equation}
\label{E:2.4}
x\mapsto\bigl\langle\phi_x^2\bigr\rangle_{n,\beta}
\end{equation}
in the limit as~$n\to\infty$. Here monotonicity arguments based on Ginibre-type inequalities (Fr\"ohlich and Park~\cite{FP}) or the technique developed in~\cite{AHPS} prove existence of a sharp threshold $\tilde\beta_\cc\in(0,\infty)$ such that, for~$x$ deep inside~$\Lambda_n$,
\begin{equation}
\lim_{n\to\infty}\,\langle\phi_x^2\rangle_{n,\beta}\begin{cases}
<\infty,\qquad&\text{if }\beta>\tilde\beta_\cc,
\\
=\infty,\qquad&\text{if }\beta<\tilde\beta_\cc.
\end{cases}
\end{equation}
For the fields, the BKT transition thus manifests itself as a localization-delocalization  transition (although, technically, the transition at~$\tilde\beta_\cc$ has yet to be linked to the conjectural threshold~$\beta_\cc$ for the polynomial decay rate of the fractional charge).

Another context where the above has been studied is one-dimensional ``discrete-Gaus\-si\-an chain'' with $1/r^2$-decaying interaction. Here, based on a duality discovered by Kjaer and Hillhorst~\cite{KH}, Fr\"ohlich and Zegarlinski~\cite{FZ} with later extensions due to Mar\-chetti~\cite{Marchetti} proved existence of both localized and delocalized phases. This class of models has recently been studied by Garban~\cite{Garban} who identified the scaling limit of the field for a range of exponents of interaction decay. Quantitative delocalization results have in turn been proved by Coquille, Dario and Ny~\cite{CDN}. Garban's work~\cite{Garban} also highlights the phenomenon of ``invisibility to integers'' in these models which refers to the effective~$\beta$ for the process at large scale being the same as for the base model; i.e.,~$\beta$ does not ``renormalize.'' The same is true for the hierarchical models but, as shown by Garban and Sep\'ulveda~\cite{GS}, not for the models on~$\Z^2$.

\newcommand{\betaeff}{{\beta_{\text{\rm eff}}}}

The extreme ends of the two phases have in the meantime been studied by perturbative methods. As mentioned earlier, the very low-temperature regime ($\beta\gg1$) can be analyzed by contour expansions. (Notably, an interesting remnant of the GFF-connection persists in the behavior of the maximum; see Lubetzky, Martinelli and Sly~\cite{LMS}.) Important inroads have also been made into the high-temperature regime ($\beta\ll1$) using the renormalization group technique. Here the field is expected to scale to a continuum Gaussian Free Field, albeit at some effective inverse temperature.
For the sine-Gordon model \eqref{E:1.4} with small~$\kappa$ this was shown by Dimock and Hurd~\cite{Dimmock-Hurd} and for the DG~model by Bauerschmidt, Park and Rodriguez~\cite{BPR1,BPR2}. Park~\cite{Park} extended the latter approach to include information about the covariance structure of the limit field.

The behavior at~$\beta_\cc$ is yet different. Indeed,  the convergence to continuum Gaussian Free Field is expected to persist but only with additional logarithmic corrections popping up in correlation functions. The only context in which this seems to have been controlled mathematically is the remarkable work of Falco~\cite{Falco1,Falco2} who determined the asymptotic form of the  fractional charge for the lattice sine-Gordon model~\eqref{E:1.4} at~$\beta=\beta_\cc(\kappa)$ for $\kappa>0$ small.  

\subsection{Hierarchical models}
\label{sec-2.3}\noindent
The present work focuses on $\Z$-modulated interface models with interactions having a hierarchical structure. Hierarchical models were originally introduced by Dyson~\cite{Dyson} as systems that are friendly to coarse-graining arguments. They soon became a testing ground for the study of critical behavior (e.g., Bleher and Sinai~\cite{BS1,BS2}). For similar reasons, they also served well in the analysis of interacting fields using the real-space renormalization group method; see e.g., Brydges~\cite{Brydges-notes}. 

Mathematicians often resort to hierarchical models when the actual model of interest is just too hard but one still wishes to make serious predictions about its behavior. This was the case in, e.g., the classical studies of ``triviality'' of the four-dimensional $\varphi^4$ and Ising models (Gaw\c edzki and Kupiainen~\cite{GK}, Hara, Hattori and Watanabe~\cite{HHW}) whose lattice counterparts have now been established as well, albeit along rather different lines. The trend to test a hierarchical setting first continues; see e.g., Hutchcroft's recent work~\cite{Hutchcroft1,Hutchcroft2} on hierarchical critical percolation.  The present  paper is a similar attempt for two-dimensional $\Z$-modulated interface models.

Our hierarchical models fall under the umbrella of GFF-like interface systems discussed after \eqref{E:2.1} but with the conductances of the associated Markovian generator \eqref{E:2.2b} taking constant values on annuli $\BB_k(x)\smallsetminus\BB_{k-1}(x)$; namely,
\begin{equation}
\frakc(x,y) := \frakc_k\quad\text{for}\quad k:=\log_{b^{1/2}}d(x,y)
\end{equation}
whenever~$x\ne y$. For $\{\frakc_k\}_{k\ge1}$ as in \eqref{E:1.6u} of Assumption~\ref{ass-1},  calculations show $\frakc_k \asymp  b^{-2k}$ and so  we have 
\begin{equation}
\frakc(x,y)\asymp d(x,y)^{-4}
\end{equation}
at large separations of~$x$ and~$y$. (Recall that $d(x,y)$ is comparable with $\Vert x-y\Vert_\infty$ at typical vertices of~$\Lambda_n$ under the embedding in~$\Z^2$.) It is worth noting that long-range conductance/percolation models over~$\Z^2$ with this kind of decay are known to exhibit interesting scaling phenomena; e.g., in the scaling of the graph-theoretical distance (B\"aumler~\cite{Baumler}) and, conjecturally, in superdiffusive behavior of random walks on such percolation graphs. The polynomial decay built into the interaction naturally amplifies the critical properties of two-dimensional hierarchical interface models.
Still, as noted to us by C.~Garban, the long-range nature does not account for everything due to the one-dimensional DG~model with $1/r^2$-decaying interaction being localized for~$\beta$ large.

\subsection{Duality with Coulomb gas}
\label{sec-2.4}\noindent
As noted earlier, the $\Z$-modulated interface models are dual to Coulomb gas models, which describe systems of charged particles interacting via Coulomb forces. A configuration of such a system is an assignment $\{q_x\colon x\in\Lambda_n\}$ of $\Z$-valued electrostatic charges to vertices of~$\Lambda_n$. The Coulomb electrostatic energy is then given by $\frac12(q,(-\Delta_n)^{-1}q)$ and the equilibrium  distribution  of the charge configuration at inverse temperature~$\beta$ is thus given by the Gibbs law
\begin{equation}
\label{E:1.5a}
\wt P_{n,\beta}(\textd q):=\frac1{\wt Z_n(\beta)}\,\texte^{\frac{\beta}{2}(q,\Delta_n^{-1}q)}\prod_{x\in\Lambda_n}\frakw(\textd q_x),
\end{equation}
where $\frakw$ is an \textit{a priori} measure on charge configurations at each vertex. Physical reasons dictate that $\frakw$ is concentrated on~$\Z$ and obeys $\frakw(-\textd q)=\frakw(\textd q)$.

The link between \eqref{E:1.1} and~\eqref{E:1.5a} is facilitated by the  so-called   \textit{sine-Gordon transformation} (also known as Siegert transformation after~\cite{Siegert}) which amounts to the following: For any~$f\colon\Lambda_n\to\R$ a calculation  shows
\begin{equation}
\label{E:1.7}
\bigl\langle\texte^{2\pi\texti (\phi,f)}\bigr\rangle_{n,\beta} = \int_{\Z^{\Lambda_n}}
\texte^{\frac{\beta'}2(q+f,\,\Delta_n^{-1}(q+f))}\prod_{x\in\Lambda_n}\frakw(\textd q_x),
\end{equation}
where
\begin{equation}
\label{E:2.12}
\beta':=4\pi^2/\beta
\end{equation}
 while
\begin{equation}
\frakw(\textd q) := a(q)\#(\textd q)
\end{equation}
for~$\{a(q)\}_{q\in\Z}$ the Fourier coefficients of~$\nu$ and $\#$ the counting measure on~$\Z$. Writing~$\langle-\rangle^\sim_{n,\beta'}$ for expectation with respect to~$\wt P_{n,\beta'}$, this becomes
\begin{equation}
\bigl\langle\texte^{2\pi\texti (\phi,f)}\bigr\rangle_{n,\beta} = \texte^{\frac{\beta'}2(f,\,\Delta_n^{-1}f)}\bigl\langle \texte^{\beta'(q,\Delta_n^{-1}f)}\bigr\rangle^\sim_{n,\beta'},
\end{equation}
where we noted that taking~$f:=0$ in \eqref{E:1.7} gives $\wt Z_n(\beta')=1$. In particular, the measures~$P_{n,\beta}$ and~$\wt P_{n,\beta'}$ determine each other. 

Through the above connection, the DG~model is dual to the  so-called   \textit{Villain gas}, which corresponds to~$a(q)=1$ for all~$q\in\Z$ and both~$\nu$ and~$\frakw$ being the counting measure on~$\Z$. For the sine-Gordon models \eqref{E:1.4} we get \eqref{E:1.7c} while for the hard-core Coulomb gas \eqref{E:1.9} we get~$a(0):=1$, $a(\pm1):= \kappa\in[0,1/2]$ and $a(q)=0$ for~$q\ne -1,0,+1$. (This is interpreted as a rule that  at most one particle can  appear at  each  vertex, giving the model its name.) Note that, by \eqref{E:2.12}, the high-temperature regime of the fields corresponds to the low-temperature regime of the Cou\-lomb gas, and \textit{vice versa}.

The connection of our models to the Coulomb gas is a central motivation for the consideration (and reason for the name) of the fractional charge correlation~\eqref{E:1.5w}. Indeed, setting $f:=\alpha\delta_x-\alpha\delta_y$ in \eqref{E:1.7} gives
\begin{equation}
\bigl\langle\texte^{2\pi\texti \alpha(\phi_x-\phi_y)}\bigr\rangle_{n,\beta} = \int_{\Z^{\Lambda_n}}
\texte^{\frac{\beta'}2(q+\alpha\delta_x-\alpha\delta_y,\,\Delta_n^{-1}(q+\alpha\delta_x-\alpha\delta_y))}\prod_{x\in\Lambda_n}\frakw(\textd q_x).
\end{equation}
The quantity
\begin{equation}
\frac12\bigl(q+\alpha\delta_x-\alpha\delta_y,\,(-\Delta_n)^{-1}(q+\alpha\delta_x-\alpha\delta_y)\bigr),
\end{equation}
has the interpretation of the Coulomb energy of the charge configuration $q+\alpha\delta_x-\alpha\delta_y$; namely, the fluctuating ``background'' distribution~$q$ with a ``static'' charge~$\alpha$ inserted at~$x$ and a ``static'' charge~$-\alpha$ inserted at~$y$. 

A power-law decay of the fractional charge correlation is indicative of a logarithmic growth  of this energy as the separation of~$x$ and~$y$ increases to infinity, while an exponential decay to a non-zero constant (which is what is expected in lattice models) makes the energy gain bounded. The change in the behavior for~$\beta'$  small  is explained by the  so-called   \textit{Debye screening} which is a mechanism through which the ambient charges shield the monopole at~$x$ from the monopole at~$y$ to make their existence at large separation less costly than if these monopoles were placed in a vacuum.

As is well known (see Brydges~\cite[Section~3.1]{Brydges-notes}), the Debye screening is far less pronounced in the hierarchical models than what is expected in lattice models. Indeed, as shown in Theorem~\ref{thm-2}, for~$\beta>\beta_\cc$ the energy still increases logarithmically but now with a smaller overall scale than for the $\beta<\beta_\cc$, where it behaves as in a vacuum. An additional iterated-log correction to the energy appears at~$\beta=\beta_\cc$ similarly as shown in the lattice sine-Gordon model with small~$\kappa$ by Falco~\cite{Falco1,Falco2}.

\section{Renormalization group flow}
\label{sec-3}\noindent
We are now ready to commence the proofs of Theorems~\ref{thm-1}--\ref{thm-2}. As noted earlier, we rely on the renormalization group method that works particularly well in the hierarchical setting. Here we review the steps that turn the model \eqref{E:1.1} into the form  amenable to analysis by this method and state the relevant conclusions. The results for $\beta\le\beta_\cc$ largely follow from our earlier work~\cite{BHu} so we give the needed proofs here. The proofs for~$\beta>\beta_\cc$ are deferred to Section~\ref{sec-6}.

\subsection{Representation as a tree-indexed Markov chain}
\label{sec-3.1}\noindent
The ($x$-space) renormalization group analysis of a Gibbs measure of the form \eqref{E:1.1} typically consists of repeated applications of two steps: a coarse-graining step and a renormalization step. In the coarse-graining step we partition the system into disjoint blocks and integrate the configuration on each block conditioned on a suitable ``representative'' value. The renormalization step then casts the integrated Gibbs weight (which is a function of the ``representative'' values) as the Gibbs weight for a new energy function with suitably adjusted, or ``renormalized,'' potentials  or  coefficients. This is done with the expectation that the resulting ``flow'' of the energy functions captures the large-scale correlations of the original Gibbs measure.
 
For the coarse-graining step in the hierarchical models \eqref{E:1.1} we use blocks that are just balls $\BB_k(x)$ in the ultrametric distance \eqref{E:1.9i}. Note that two such balls either coincide  or are disjoint and so~$\Lambda_n$ partitions, for each $k=1,\dots,n$, into~$b^{n-k}$ of such disjoint balls which we will refer to as $k$-blocks. The choice of the ``representative value''  relies on a ``finite-range'' decomposition of the inverse Laplacian~$\Delta_n^{-1}$ stated in: 

\begin{lemma}
Given $n\ge1$, suppose that~$\{\frakc_k\}_{k=1}^{n+1}$ is related to a positive sequence
$\{\sigma_k^2\}_{k=0}^n$ as in \twoeqref{E:1.5u}{E:1.6u}. Writing $Q_kf(x):=b^{-k}\sum_{y\in\BB_k(x)}f(y)$ for the orthogonal projection of~$f\colon\Lambda_n\to\R$ on its averages over $k$-blocks, we then have
\begin{equation}
\label{E:3.2ui}
(-\Delta_n)^{-1} = \sigma_0^2 Q_0 +\sum_{k=1}^n \sigma_k^2 \,b^k Q_k.
\end{equation}
\end{lemma}

\begin{proofsect}{Proof}
We start by recalling facts from the proof of~\cite[Lemma~3.1]{BHu}: The family of operators $\{Q_k-Q_{k+1}\colon k=0,\dots,n\}$, subject to the convention~$Q_{n+1}:=0$, are orthogonal projections on orthogonal subspaces of~$\ell^2(\Lambda_n)$ and $Q_0=\sum_{k=0}^n (Q_k-Q_{k+1})$ is a partition of the identity. As a consequence, any operator~$\LL_n$ on~$\ell^2(\Lambda_n)$ of the form
\begin{equation}
\label{E:3.4i}
\LL_n:=-u_n^{-1}Q_0+\sum_{k=1}^n(u_{k-1}^{-1}-u_k^{-1})(Q_k-Q_0)
\end{equation}
for $\{u_k\}_{k=0}^n$ such that~$u_0>0$ and~$u_k>u_{k-1}$ for~$k=1,\dots,n$ inverts to
\begin{equation}
\label{E:3.5ui}
\LL_n^{-1}:=-u_0Q_0-\sum_{k=1}^n(u_k-u_{k-1})Q_k,
\end{equation}
see~\cite[Eqs.~(3.11--3.13)]{BHu}.

Now observe that the operator $\Delta_n$ from \eqref{E:1.2} takes the form \eqref{E:3.4i} if
\begin{equation}
\label{E:3.6i}
u_n^{-1}=\frakc_{n+1}
\end{equation}
and
\begin{equation}
u_{k-1}^{-1}-u_k^{-1}= b^k\frakc_k, \quad k=1,\dots,n,
\end{equation}
while \eqref{E:3.5ui} matches \eqref{E:3.2ui} if
\begin{equation}
\label{E:3.8i}
u_k = \sum_{j=0}^k b^j\sigma_j^2,\quad k=0,\dots,n.
\end{equation}
As a calculation shows, under \twoeqref{E:1.5u}{E:1.6u} we get \twoeqref{E:3.6i}{E:3.8i} as desired. 
\end{proofsect}

\begin{remark}
\label{rem-hierarchy}
In the literature (see, e.g.,~\cite[Section~1.3]{BB}) the massive hierarchical Laplacian is sometimes presented in the form
\begin{equation}
\label{E:3.7ii}
-m^2Q_0-\sum_{k=0}^{n-1}L^{-2k}(Q_k-Q_{k+1}),
\end{equation}
where $m^2$ is the ``mass-squared'' and where we write the coefficient using the scale~$L$. Noting that \eqref{E:3.4i} rewrites as $-\sum_{k=0}^nu_k^{-1}(Q_k-Q_{k+1})$, the form \eqref{E:3.7ii} agrees with \eqref{E:3.4i} provided we set $u_k^{-1}:=m^2+L^{-2k}$ for~$k=0,\dots,n-1$ and~$u_n^{-1}:=m^2$. (Observe that \eqref{E:3.6i} then gives $\frakc_{n+1}=m^2$.) This in turn matches \eqref{E:3.8i} with~$b:=L^2$ provided that
\begin{equation}
\sigma_k^2 = \frac{L^{-4k}(L^2-1)}{(m^2+L^{-2k})(m^2+L^{2-2k})}, \quad k= 1,\dots,n-1,
\end{equation}
 with the ``boundary'' cases given as
\begin{equation}
\sigma_0^2=\frac1{m^2+1} \quad\text{ and }\quad \sigma_n^2 = \frac{L^{2-4n}}{m^2(m^2+L^{2-2n})}.
\end{equation}
 Assuming that $m^2/L^{-2n}$ is bounded between two positive constants uniformly in~$n\ge1$, a calculation shows that, for $k=1,\dots,n-1$,
\begin{equation}
\sigma_k^2-(1-L^{-2}) = O(L^{-2(n-k)})
\end{equation}
while  $\sigma_0^2=1+O(L^{-2n})$ and  $\sigma_n^2=O(1)$. Modulo scaling by $1-L^{-2}$, the operator \eqref{E:3.7ii}  thus  conforms to Assumption~\ref{ass-1} with~$\{\frakd_k\}_{k\ge0}$ decaying exponentially.  This example is actually the prime motivation for the setting in Assumption~\ref{ass-1}. 
\end{remark}

The representation \eqref{E:3.2ui} allows us to view $\texte^{\frac12\beta(\phi,\Delta_n\phi)}$ in \eqref{E:1.1} as a convolution of~$n+1$ probability densities of Gaussian fields on~$\Lambda_n$ with covariances $\sigma_0^2Q_0$, $\sigma_1^2bQ_1$,\dots, $\sigma_n^2b^nQ_n$, respectively. (A caveat is that these densities are singular because the field with covariance~$Q_k$ is constant on each~$k$-block.) Adding the integral over~$\phi$ with respect to the product measure $\prod_{x\in\Lambda_n}\nu(\textd\phi_x)$, we perform one integral after another starting from that for~$\phi$ itself and proceeding to the field with covariance~$\sigma_0^2Q_0$, then to the field with covariance~$\sigma_1^2bQ_1$, etc. The ``representative'' value of a block is, at each step, the sum of the Gaussian fields yet to be integrated. Since the field is constant on each  $k$-block, we simply take its value at any point within the block.
 
As a consequence (see \cite[Section~3.3]{BHu}), after~$k$ integrals have been performed, the resulting ``renormalized'' Gibbs measure admits a density with respect to the law of the Gaussian field $\phi^{(k)}$ on~$\Lambda_{n-k}$ with covariance $\sigma_k^2Q_{0}+\sigma_{k+1}^2 b Q_{1}+ \dots+\sigma_n^2 b^{n-k}Q_{n-k}$ that  is proportional to 
\begin{equation}
\exp\Bigl\{-\sum_{x\in\Lambda_{n-k}}bv_{k-1} (\phi^{(k)}_x) \Bigr\}.
\end{equation}
Here $\{v_k\}_{k=0}^n$ is a sequence of potentials defined, for~$k=0,\dots,n-1$, recursively  by
\begin{equation}
\label{E:3.5i}
\texte^{-v_{k+1}(z)}:=\int \texte^{-b v_k(z+\zeta)}\mu_{\sigma_{k+1}^2/\beta}(\textd\zeta)
\end{equation}
where~$\mu_{\sigma^2}$ denotes the law of~$\NN(0,\sigma^2)$, with the ``initial value'' set as
\begin{equation}
\texte^{-v_0(z)}:=\int \texte^{-\frac12\beta\sigma_0^{-2}(z-\zeta)^2}\nu(\textd\zeta).
\end{equation}
The $1$-periodicity of~$\nu$ implies that all~$v_k$'s are $1$-periodic. The renormalization group flow is thus encoded by the sequence $\{v_k\}_{k=0}^n$ of functions of one variable. (This sequence depends on~$n$, $\beta$ and $\{\sigma_k^2\}_{k=0}^n$ but we suppress that from the notation.)

The primary output of the above procedure is a representation of the normalization constant of the Gibbs measure \eqref{E:1.1} as
\begin{equation}
Z_n(\beta)=\frac{\beta^{\frac{|\Lambda_n|}2}\sqrt{\det(-\Delta_n)}}{\bigl(\frac{2\pi\sigma_0^2}\beta\bigr)^{\frac{|\Lambda_n|}2}}\,\texte^{-v_n(0)},
\end{equation}
 where the numerator, resp., the denominator in the prefactor are the quantities that normalize $\phi\mapsto\texte^{\frac\beta2(\phi,\Delta_n\phi)}$, resp., $\phi\mapsto \texte^{-\frac\beta2\sigma_0^{-2}\sum_{x\in\Lambda_n}\phi_x^2}$ into probability densities. (Note that all of the $\nu$-dependence is now hidden inside~$v_n$.) 
 In order to control expectations of relevant local observables, standard treatments of rigorous renormalization group proceed by incorporating the observable into suitably modified potentials  whose ``flow''  then needs  to be controlled alongside~$\{v_k\}_{k=0}^n$.

In our earlier work~\cite{BHu} we instead took a different approach that is based on representing the full Gibbs measure \eqref{E:1.1} as a \textit{tree-indexed Markov chain}. Consider a $b$-ary rooted tree~$\T_n$ of depth~$n$ with the root denoted by~$\varrho$ and note that, keeping the same root, $\T_n$ naturally embeds~$\T_k$ for each~$k=0,\dots,n$. Let~$m\colon\T_n\smallsetminus\{\varrho\}\to\T_n$ be the map assigning to~$x$ its parent in the unique path to the root. For~$k=0,\dots,n$, define the probability kernels
\begin{equation}
\label{E:1.23}
\frakp_{k}(\textd\varphi|\varphi'):=\begin{cases}
\texte^{v_{k}(\varphi')-bv_{k-1}(\varphi)}\,\mu_{\sigma_{k}^2/\beta}(-\varphi'+\textd\varphi),\qquad&\text{if }k\ge1,
\\
\texte^{v_0(\varphi')-\frac\beta2\sigma_0^{-2}(\varphi-\varphi')^2}\,\nu(\textd\varphi),\qquad&\text{if }k=0.
\end{cases}
\end{equation} 
 Now  generate a family of random variables
\begin{equation}
\label{E:3.6}
\{\varphi_x\colon x\in\T_n\}
\end{equation}
as follows: Sample $\varphi_\varrho$ from~$\frakp_n(\cdot|\,0\,)$. Then, for each~$k=1,\dots,n$, assuming that the values of~$\varphi$ on~$\T_{k-1}$ have already been sampled, draw~$\varphi_x$ for each~$x\in\T_{k}\smallsetminus\T_{k-1}$ from $\frakp_{n-k}(\cdot|\varphi_{m(x)})$, independently for different~$x$. The joint law of the random variables thus takes the form
\begin{equation}
\frakp_n(\textd\varphi_\varrho|\,0\,)\prod_{k=1}^{n}\,\prod_{x\in\T_k\smallsetminus\T_{k-1}}\frakp_{n-k}(\textd\varphi_x|\,\varphi_{m(x)})
\end{equation}
We then have:

\begin{lemma}
\label{lemma-3.2}
Under the canonical identification of~$\Lambda_n$ with the leaves of~$\T_n$, the restriction of the family \eqref{E:3.6} to~$\Lambda_n$ is distributed according to~$P_{n,\beta}$ from \eqref{E:1.1}.
\end{lemma}

\begin{proofsect}{Proof}
This is a restatement of \cite[Lemma~3.2]{BHu} modulo the fact that there the proof was performed only for~$\{\sigma_k^2\}_{k=0}^n$ equal to one. We  leave the modifications to the reader.
\end{proofsect}

Note that the definition implies that the values of \eqref{E:3.6} along any path from the root to a leaf-vertex are an ordinary (time-inhomogeneous) Markov chain with transition probabilities \eqref{E:1.23}. This will be very useful in our later calculations.

\subsection{Results for renormalization group iterations}
\label{sec-3.2}\noindent
Our ability to control the above tree-indexed Markov chain depends very strongly on our ability to control the differences $v_{k}(\phi')-b v_{k-1}(\phi)$ for large values of~$k$ (and~$n$). In high-temperature approaches to this problem (see, e.g., Bauerschmidt and Bodineau~\cite{BB}) this is done by linearization of \eqref{E:3.5i}. However, linearization becomes inefficient if we want to work uniformly up to~$\beta_\cc$, or even beyond, so in~\cite{BHu} we instead followed the flow of the Fourier coefficients of  $\texte^{-v_k}$,  defined for~$k=0,\dots,n$ by 
\begin{equation}
\label{E:3.15ii}
a_k(q):=\int_0^1 \texte^{-v_{k}(z) - 2\pi\texti q z}\,\textd z.
\end{equation}
As shown in~\cite[Lemma-4.1]{BHu}, in light of \eqref{E:3.5i} these coefficients iterate as
\begin{equation}
\label{E:3.8}
a_{k+1}(q) := \sum_{\begin{subarray}{c}
\ell_1,\dots,\ell_b\in\Z\\\ell_1+\dots+\ell_b=q
\end{subarray}}
\biggl[\,\prod_{i=1}^b a_k(\ell_i)\biggr]\,\theta_{k+1}^{q^2},\quad q\in\Z,
\end{equation}
where
\begin{equation}
\label{Eq:3.14w}
\theta_k:=\texte^{-\frac{2\pi^2}\beta\sigma_k^2}
\end{equation}
 and where the ``initial'' value is set as
\begin{equation}
\label{E:3.18q}
a_{0}(q):=\sqrt{\frac{2\pi\sigma_0^2}\beta}\,a(q)\,\theta_0^{q^2}\quad\text{for}\quad a(q):=\int_{[0,1)}\texte^{-2\pi\texti q z}\nu(\textd z). 
\end{equation}
 As also noted in~\cite[Lemma-4.1]{BHu} (whose proof only needs that~$0\le\theta_{k+1}<1$), the  conditions in Assumption~\ref{ass-2} ensure that $a_k(q)>0$ for all~$k\ge0$ and~$q\in\Z$ and that $\{a_k(q)\}_{q\in\Z}\in\ell^1(\Z)$ for all~$k=0,\dots,n$. (As for the~$v_k$'s,  the~$a_k$'s  also  depend on~$n$, $\beta$ and $\{\sigma_k^2\}_{k=0}^n$ but we do not mark that explicitly in the notation.)

We will now state our results concerning the iterations \eqref{E:3.5i} and \eqref{E:3.8}. Our first theorem concerns the subcritical~$\beta$:

\begin{theorem}[Subcritical flow]
\label{thm-3}
Suppose that Assumptions~\ref{ass-1}--\ref{ass-2} hold. For each $b\ge2$ and each~$\beta>0$ with~$\beta<\beta_\cc$ there exist $\eta>0$ and $C>0$ such that for all~$n\ge k\ge0$,
\begin{equation}
\label{E:3.12}
0< \frac{a_k(q)}{a_k(0)}\le C\texte^{-\eta (k|q|+q^2)},\quad q\in\Z.
\end{equation}
Moreover, each $v_k$ is a $C^\infty$ function with 
\begin{equation}
\label{E:3.13}
\sup_{z,z'\in\R}\bigl|v_{k+1}(z)-b v_k(z')\bigr|\le C\texte^{-\eta k}
\end{equation}
 for all $n>k\ge0$,  and
\begin{equation}
\label{E:3.14}
\sup_{z\in\R}\max\bigl\{ |v_k'(z)|,|v_k''(z)|\bigr\}\le C\texte^{-\eta k}
\end{equation}
for all~$n\ge k\ge 0$.
\end{theorem}

The next result provides a similar statement at~$\beta_\cc$, where control of the iterations is  more subtle than in the subcritical cases.

\begin{theorem}[Critical flow]
\label{thm-4}
Suppose that Assumptions~\ref{ass-1}--\ref{ass-2} hold and assume $\beta=\beta_\cc$. For each~$b\ge2$ there exist~$\eta>0$ and $C>0$ such that for all~$n\ge k\ge0$,
\begin{equation}
\label{E:3.15}
0< \frac{a_k(q)}{a_k(0)}\le \frac{C\texte^{-\eta q^2}}{(1+\sqrt k)^{|q|}}, \quad  q\in\Z.
\end{equation}
Moreover,  each $v_k$ is a $C^\infty$ function with
\begin{equation}
\label{E:3.16}
\sup_{z,z'\in\R}\bigl|v_{k+1}(z)-b v_k(z')\bigr|\le \frac C{\sqrt{1+k}}
\end{equation}
 for all $n>k\ge0$,  and
\begin{equation}
\label{E:3.17}
\sup_{z\in\R}\max\bigl\{ |v_k'(z)|,|v_k''(z)|\bigr\}\le \frac C{\sqrt{1+k}}
\end{equation}
for all~$n\ge k\ge0$. Furthermore, we have
\begin{equation}
\label{Eq:3.18u}
\frac{a_k(1)}{a_k(0)}=\frac{1}{\sqrt k}\biggl[\frac{b^3-1}{(b-1)^2(b+1)^3}\biggr]^{1/2}+O(k^{-1})+O\biggl(\,\,k^{-1/2}\!\!\!\!\!\!\!\sum_{j\ge\min\{\sqrt k,n-k\}}\frakd_j\biggr),
\end{equation}
 where~$\{\frakd_j\}_{j\ge0}$ is the sequence from Assumption~\ref{ass-1}.  Consequently,
\begin{equation}
\label{Eq:1.44}
v_k'(z) = \frac{4\pi}{\sqrt k}\biggl[\frac{b^3-1}{(b-1)^2(b+1)^3}\biggr]^{1/2}\sin(2\pi z) +O(k^{-1})+O\biggl(\,\,k^{-1/2}\!\!\!\!\!\!\!\sum_{j\ge\min\{\sqrt k,n-k\}}\frakd_j\biggr),\end{equation}
where~$o(1)\to0$ as~$\min\{k,n-k\}\to\infty$, uniformly in~$z\in\R$. 
\end{theorem}

The asymptotic \eqref{Eq:1.44} implies that $\{|v_k'|^2\}_{k=0}^n$ is not summable uniformly in~$n$ which, as we will see in Section~\ref{sec-4}, is the root cause of the doubly-logarithmic correction to the covariance structure at~$\beta=\beta_\cc$. The logarithmic correction to the fractional charge asymptotic at~$\beta=\beta_\cc$ can in turn be traced to \eqref{Eq:3.18u}.

Our final theorem in this section deals with supercritical~$\beta$. The statement only applies to~$\beta$ slightly over~$\beta_\cc$. Denote
\begin{equation}
\label{Eq:3.27w}
\theta:=\texte^{-\frac{2\pi^2}\beta}
\end{equation}
and observe that~$\beta>\beta_\cc$ is equivalent to~$b\theta>1$.
We then have:

\begin{theorem}[Supercritical flow]
\label{thm-5}
Suppose Assumption~\ref{ass-1} holds. For each~$b\ge2$ there exists~$\epsilon>0$ and, for all~$\beta>0$ with~$1<b\theta<1+\epsilon$, there exist $\eta>0$, $C>0$ and a sequence $\{\lambda_\star(q)\}_{q\in\Z}$ with $\lambda_\star(0)=1$, $\lambda_\star(q)\le\theta^{q^2}$ for~$q\in\Z$ and
\begin{equation}
\label{E:3.28u}
0<\lambda_\star(q)=\lambda_\star(-q)\le \bigl(2b^{1/2}\sqrt{b\theta-1}\,\bigr)^{|q|},\quad q\in\Z,
\end{equation}
for which the following is true: For all initial~$\nu$ subject to Assumption~\ref{ass-2} there exists~$k_0\ge0$ such that for all integer $n$ and~$k$ with $\min\{k,n-k\}\ge k_0$ and all~$q\in\Z$,
\begin{equation}
\label{Eq:3.30w}
a_k(q)\le a_k(0) \bigl(2b^{1/2}\sqrt{b\theta-1}\,\bigr)^{|q|}
\end{equation}
and
\begin{equation}
\label{Eq:3.31}
\biggl|\frac{a_k(q)}{a_k(0)}-\lambda_\star(q)\biggr|\le C^{|q|}\biggl[\texte^{-\eta k}+\sum_{j=0}^k\texte^{-\eta(k-j)}\,\frakd_{\min\{j,n-j\}}\biggr]
\end{equation}
 hold  with~$\{\frakd_j\}_{j\ge0}$ denoting the sequence from Assumption~\ref{ass-1}. Moreover, each $v_k$ is a $C^\infty$ function with   with~$\{v_k'\}_{k\ge0}$ and~$\{v_k''\}_{k\ge0}$ uniformly bounded and, defining $v_\star\colon\R\to\R$ by 
\begin{equation}
\label{E:3.34a}
\texte^{-v_\star(z)}:=\Biggl(\,\sum_{\begin{subarray}{c}
q_1,\dots,q_b\in\Z\\q_1+\dots+q_b=0
\end{subarray}}
\prod_{i=1}^b \lambda_\star(q_i)\Biggr)^{-\frac1{b-1}}
\sum_{q\in\Z}\lambda_\star(q)\texte^{2\pi\texti q z},
\end{equation}
we have  $v_k(z)-b v_{k-1}(z')\to v_\star(z)-bv_\star(z')$ and $v_k'(z)\to v'_\star(z)$ as~$\min\{k,n-k\}\to\infty$, uniformly in $z,z'\in\R$.  (The existence of~$v_\star'$ is ensured by \eqref{E:3.28u}.) In addition, assuming $\{\frakd_k\}_{k\ge0}$ decays exponentially fast, there exist  $\eta'>0$ and~$C'>0$ such that
\begin{equation}
\label{Eq:3.32a}
\sup_{z,z'\in\R}\,\Bigl| v_{k+1}(z)-b v_{k}(z')- \bigl[v_\star(z)-bv_\star(z')\bigr]\Bigr|\le C'\texte^{-\eta'\min\{k,n-k\}}
\end{equation}
 for all $n > k \ge 0$ and 
\begin{equation}
\label{Eq:3.32}
\sup_{z\in\R}\,\bigl| v_k'(z)-v'_\star(z)\bigr|\le C'\texte^{-\eta'\min\{k,n-k\}}
\end{equation}
 for all~$n\ge k\ge0$. 
\end{theorem}

\begin{figure}[t]
\refstepcounter{obrazek}
\label{fig3}
\vglue5mm
\newcommand\skala{1.85in}
\centerline{\includegraphics[width=\skala]{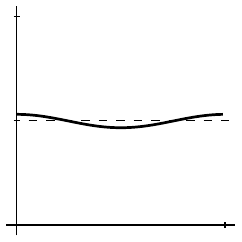}\hglue2mm \includegraphics[width=\skala]{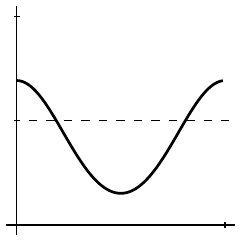}\hglue2mm\includegraphics[width=\skala]{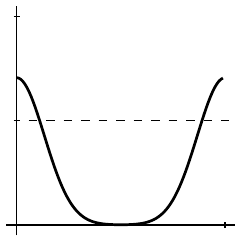}}
\vspace{.2in}
\begin{quote}
\fontsize{9}{5}\selectfont
{\sc Fig~\theobrazek.} Plots of $z\mapsto\texte^{-v_\star(z)}$ from \eqref{E:3.34a} (and solving \eqref{E:3.33}) for $b=2$ and $\theta=0.501$ (left), $\theta=0.6$ (middle) and $\theta=0.84$ (right). (Recall that~$\theta=0.5$ at~$\beta=\beta_\cc$.) The dashed line is at height~$1$, the mark on the vertical line as at height~$2$. 
\normalsize
\end{quote}
\end{figure}

We emphasize that~$\{\lambda_\star(q)\}_{q\in\Z}$ and thus also~$v_\star$ do not depend on~$\nu$; indeed, they represent a ``nontrivial'' fixed point of the renormalization group flow as in
\begin{equation}
\label{lambda-rec}
\lambda_\star(q)=\frac{\sum_{\bar\ell\in\Xi_{b}(q)}\prod_{i=1}^{b}\lambda_\star(\ell_i)}{\sum_{\bar\ell'\in\Xi_b(0)}\prod_{i=1}^{b}\lambda_\star(\ell_i')}\,\theta^{q^2},\quad q\in\Z,
\end{equation}
where we abbreviated
 \begin{equation}
\label{E:5.73i}
\Xi_b(q):=\bigl\{(\ell_1,\dots,\ell_{b})\in\Z^{b}\colon\ell_1+\dots+\ell_{b}=q\bigr\}.
\end{equation}
This means that~$v_\star$ is a non-zero solution to
\begin{equation}
\label{Eq:3.33}
\texte^{-v_\star(z)}=\int \texte^{-bv_\star(z+\zeta)}\mu_{1/\beta}(\textd\zeta),
\end{equation}
where, we recall, $\mu_{1/\beta}$ is the law of~$\NN(0,1/\beta)$.
As our proofs show (see Theorem~\ref{thm-6.1}), such a fixed point is \textit{unique} and attractive to all the $1$-periodic  measures  whose Fourier coefficients are positive and obey  \twoeqref{E:1.6w}{E:1.7a}.  (These assumptions matter as for $b\theta^{p^2}\le1$, functions that are $1/p$-periodic are still attracted to the ``trivial'' fixed point.) We do not have explicit expressions for~$v_\star$ or~$\{\lambda_\star(q)\}_{q\in\Z}$. The best we can offer is a characterization of  their  $b\to\infty$ asymptotic behavior; see Remark~\ref{rem-large-b}. Note that setting~$v_\star:=0$ for~$\beta\le\beta_\cc$ puts \eqref{E:3.13} and \eqref{E:3.16} in the same form as \eqref{Eq:3.32a}.

\subsection{Proof of Theorem~\ref{thm-3}}
\label{sec-3.3}\noindent
The above convergence statements for the subcritical and critical regimes require only relatively minor adaptations of the results already proved in~\cite{BHu}, and so we prove them right away. The main new obstacle is the fact that~\cite{BHu} assumed~$\sigma_k^2=1$ for all~$k\ge0$ while for us the equality holds only asymptotically.
We will suppose that Assumptions~\ref{ass-1}--\ref{ass-2} hold throughout and stop referencing these in the statements of lemmas.

We start with Theorem~\ref{thm-3}. First we recall an observation from~\cite{BHu} that drives many of the subsequent arguments:

\begin{lemma}
\label{lemma-3.5}
For all $\beta>0$, $n > k\ge0$ and~$q\ge0$,
\begin{equation}
\label{E:3.25}
\frac{a_{ k+1}(q+1)}{a_{ k+1}(q)}\le b\,\theta_{ k+1}^{(q+1)^2-q^2}\sup_{\ell\ge0}\frac{a_{ k}(\ell+1)}{a_{ k}(\ell)}.
\end{equation} 
\end{lemma}

\begin{proofsect}{Proof}
This is a restatement of \cite[Lemma~4.2]{BHu} with~$\theta$ allowed to depend on~$k$.
\end{proofsect}

As a consequence we obtain:

\begin{lemma}
\label{lemma-3.6}
For all $\beta>0$, $n\ge k\ge0$ and~$q\in\Z$,
\begin{equation}
\label{E:3.26}
\frac{a_k(q)}{a_k(0)}\le\theta_k^{q^2}\biggl(\,\hat c\,\prod_{j=0}^{ k-1}(b\theta_j)\biggr)^{|q|},
\end{equation}
where $\hat c:=\sup_{\ell\ge0}\frac{a(\ell+1)}{a(\ell)}$ for~$\ell\mapsto a(\ell)$ being the Fourier coefficients of~$\nu$.
\end{lemma}

\begin{proofsect}{Proof}
Denote $c_k:=\sup_{\ell\ge0}\frac{a_k(\ell+1)}{a_k(\ell)}$. Then \eqref{E:3.25} along with $(q+1)^2-q^2\ge1$ for~$q\ge0$ gives $c_{k+1}\le (b\theta_{k+1}) c_k$  with~$c_0\le\hat c \theta_0$.  Iterating, we get $c_k\le  b^{-1}\hat c\prod_{j=0}^{k}(b\theta_j)$ for all~$n\ge k\ge0$. Plugging this in \eqref{E:3.25} and iterating yields the claim.
\end{proofsect}

We are now ready to give:

\begin{proofsect}{Proof of Theorem~\ref{thm-3}}
Suppose~$0<\beta<\beta_\cc$. We start with \eqref{E:3.12}, which for~$\sigma_k^2$ equal to one was shown already in~\cite[Lemma~4.3]{BHu}.  Assumption~\ref{ass-1} shows $\theta_{\text{max}}:=\sup_{n\ge k\ge0}\theta_k<1$ and $\prod_{j=0}^{k-1}b\theta_j\le \theta^{-\tilde c}(b\theta)^k$ where~   $\tilde c:=2\sum_{j\ge0}\frakd_j$ . Using \eqref{E:3.26} we then get
\begin{equation}
\frac{a_k(q)}{a_k(0)}\le \theta_{\text{max}}^{q^2}\bigl[\hat c\theta^{-\tilde c}(b\theta)^k\bigr]^{|q|}.
\end{equation}
 Setting,  with some waste  for a later  convenience,  $\texte^{-\eta}:=\max\{\sqrt{b\theta},\sqrt{\theta_{\text{max}}}\}$,  we get \eqref{E:3.12} with~$C:=\sup_{q\ge0}\theta_{\text{max}}^{q^2/2}[\hat c\theta^{-\tilde c}]^q$.
 The positivity follows from iterations of \eqref{E:3.8} and the assumption that the Fourier coefficients of~$\nu$ are strictly positive.

Concerning \eqref{E:3.13}, we observe that, by $\{a_k(q)\}_{q\in\Z}\in\ell^1(\Z)$,
\begin{equation}
\label{E:3.27}
\texte^{-v_k(z)} = \sum_{q\in\Z}a_k(q)\texte^{2\pi\texti qz}
\end{equation}
 with the left-hand side continuous and, by \eqref{E:3.5i}, strictly positive for all~$z\in\R$. Hence, the~$v_k$'s are also continuous. 
The bound \eqref{E:3.12} gives 
\begin{equation}
\bigl|\texte^{-v_k(z)}-a_k(0)\bigr|\le 2a_k(0)\sum_{q\ge1}C\texte^{-\eta(kq+q^2)}\le\frac{2C\texte^{-k\eta}}{1-\texte^{-k\eta}} a_k(0).
\end{equation}
 For~$k$  so large that $\max\{1,C\}\texte^{-k\eta}\le1/8$  we get 
\begin{equation}
\label{E:3.29}
\sup_{z\in\R}\bigl|\,a_k(0)^{-1}\texte^{-v_k(z)}-1\bigr|\le 3C\texte^{-\eta k}\le\frac12.
\end{equation}
 Once $k$ is so large that  the quantity in the large parentheses in \eqref{E:3.26} is less than $\texte^{-\eta k}$,  this now implies \eqref{E:3.13}  along the same argument that proved \cite[Eq.~(3.40)]{BHu}. The remaining~$k$ and~$n$ are handled  directly by noting that, under Assumption~\ref{ass-1}, $\texte^{-v_k}$ is bounded above and below by positive constants that depend only on~$k$. Hence~$|v_k|$ is bounded uniformly in~$n$ with~$n\ge k$ and so \eqref{E:3.13} follows by relabeling~$C$. 

For the corresponding bound on the derivatives of~$v_k$, first note that \eqref{E:3.12} permits us to differentiate the series in \eqref{E:3.27} term-by-term to get
\begin{equation}
\label{E:3.30}
v_k'(z) =-\texte^{v_k(z)}\sum_{q\in\Z}(2\pi\texti q)a_k(q)\texte^{2\pi\texti qz}.
\end{equation}
 Using \eqref{E:3.29} along with the uniform boundedness of $a_k(0)\texte^{v_k}$ for each~$k$  we conclude that the~$v_k'$ are continuous and bounded on~$\R$, uniformly in~$n\ge k\ge0$. It thus suffices to prove \eqref{E:3.14} for~$k$ sufficiently large. Here we invoke the bound \eqref{E:3.12} to get
\begin{equation}
\bigl|v_k'(z)\bigr|\le  a_k(0)\texte^{v_k(z)}\sum_{q\ge1} 4\pi q C\texte^{-\eta kq}
\le  a_k(0)\texte^{v_k(z)}\frac{4\pi C\texte^{-\eta k}}{(1-\texte^{-\eta k})^2}.
\end{equation}
Since $a_k(0)\texte^{v_k(z)}\le 2$ whenever \eqref{E:3.29} is in force, the right-hand side is at most $32\pi C\texte^{-\eta k}$ as soon as $\texte^{-\eta k}\le1/8$. This proves \eqref{E:3.14} for the first derivative. For the bound on the second derivative we differentiate \eqref{E:3.30} one more time and apply a similar reasoning, along with the bound on the first derivative. We leave the details to the reader. 
\end{proofsect}

\subsection{Bounds on Fourier coefficients}
For the critical case, we first need to establish estimates on the Fourier coefficients $a_k(q)$. We start with a bound that drove the analysis of the critical case in~\cite{BHu}:

\begin{lemma}
\label{lemma-3.7}
For all $\beta>0$ and~$n > k\ge0$,
\begin{equation}
\label{E:3.28}
\frac{a_{ k+1}(1)}{a_{ k+1}(0)}\le\theta_{ k+1}\frac{\displaystyle\frac{a_{ k}(1)}{a_{ k}(0)}}{\displaystyle 1+\binom b2\Bigl(\frac{a_{ k}(1)}{a_{ k}(0)}\Bigr)^2}+(b-1)\theta_{ k+1}\sup_{\ell\ge0}\frac{a_{ k}(\ell+1)}{a_{ k}(\ell)}.
\end{equation}
\end{lemma}

\begin{proofsect}{Proof}
This is a restatement of \cite[Lemma~4.5]{BHu} with~$\theta$ allowed to depend on~$k$.
\end{proofsect}

Next we show that the supremum on the right of \eqref{E:3.28} exhibits polynomial decay:

\begin{lemma}
\label{lemma-3.8}
Assume~$\beta=\beta_\cc$. Then there exists a constant~$C>0$ such that for all~$n\ge k\ge0$,
\begin{equation}
\label{E:3.33}
\sup_{\ell\ge0}\frac{a_k(\ell+1)}{a_k(\ell)}
\le\frac{C}{\sqrt{1+k}}.
\end{equation}
\end{lemma}

\begin{proofsect}{Proof}
We will adapt the proofs of~\cite[Lemma~4.6]{BHu} and~\cite[Theorem~3.5]{BHu} to allow~$\sigma_k^2$ depend on~$n$ and~$k$. Abbreviate the supremum in \eqref{E:3.33} as~$c_k$ and let~$\alpha_k$ be the unique number in~$(0,1)$ such that
\begin{equation}
\label{E:3.46iu}
\alpha_k = \frac1{ 1+\binom b2 c_k^2\alpha_k^2}.
\end{equation}
Next  we will prove  that~$k\mapsto c_k$ is bounded. Indeed, at $\beta=\beta_\cc$ we have $\theta_k=b^{-\sigma_k^2}$ and so $b\theta_k= b^{1-\sigma_k^2}$. The argument from the proof of Lemma~\ref{lemma-3.6}  along with the inequality $c_0\le \theta_0 \hat c \le b\theta_0 \hat c$   implied by \eqref{E:3.18q}  then show 
\begin{equation}
\label{E:3.44w}
c_k\le\hat c\exp\biggl\{(\log b)\sum_{j=0}^k|1-\sigma_j^2|\biggr\}
\end{equation}
with the sum is bounded uniformly in~$n\ge k\ge0$ thanks to Assumption~\ref{ass-1}. This along with~$\alpha_k\le1$ implies $\inf_{n\ge k\ge0}\alpha_k\ge[1+\binom b2(\sup_{n\ge k\ge0}c_k)^2]^{-1}>0$.

We will now repeat the argument from the proof of~\cite[Lemma~4.6]{BHu} to get an iterative bound    on~$c_k$.  We start by the inequality 
\begin{equation}
\label{E:3.34}
\frac{a_{ k+1}(1)}{a_{ k+1}(0)}\le\frac{\theta_{ k+1}\,c_{ k}}{1+\binom b2 c_{ k}^2\alpha_{ k}^2}+(b-1)\theta_{ k+1} c_{ k}.
\end{equation}
For $a_k(1)/a_k(0)> \alpha_k c_k$, this is obtained by bounding the denominator in \eqref{E:3.28} from below by $1+\binom b2\alpha_k^2c_k^2$ and then applying $a_k(1)/a_k(0)\le c_k$ in the numerator. For $a_k(1)/a_k(0)\le \alpha_k c_k$, we instead drop the denominator in \eqref{E:3.28} altogether, invoke $a_k(1)/a_k(0)\le \alpha_k c_k$ in the numerator and then observe that, by \eqref{E:3.46iu}, the right-hand side of \eqref{E:3.34} equals $(b-1+\alpha_{ k})\theta_{k+1}c_{ k}$. With \eqref{E:3.34} in hand,  observe that Lemma~\ref{lemma-3.5} also gives
\begin{equation}
\label{E:3.36}
\frac{a_{ k+1}(q+1)}{a_{ k+1}(q)}\le b\,\theta_{ k+1}^{(q+1)^2-q^2}c_{ k}
\end{equation}
with $(q+1)^2-q^2\ge3$ for~$q\ge1$.
Noting  again  that the right-hand side of \eqref{E:3.34} can been written as $(b-1+\alpha_{ k})\theta_{k+1}c_{ k}$, we get
\begin{equation}
\label{E:3.37}
c_{ k+1}\le\frac{\theta_{ k+1}\,c_{ k}}{1+\binom b2 c_{ k}^2\alpha_{ k}^2}+(b-1)\theta_{ k+1} c_{ k}
\end{equation}
 as soon as $\min\{k,n-k\}$ is  so large that $b\theta_{ k+1}^3\le(b-1+\alpha_{ k+1})\theta_{ k}$.

 To deal with~$k$ dependence of~$\theta_k$ and~$\alpha_k$ in \eqref{E:3.37}, denote 
\begin{equation}
\label{E:3.48w}
\tilde c_k:=c_k\exp\biggl\{\log(1/\theta)\sum_{j=0}^k(\sigma_j^2-1)\biggr\},
\end{equation}
 and abbreviate 
\begin{equation}
\bar\alpha:=\binom b2\bigl(\inf_{n\ge k\ge0}\alpha_k^{ 2}\bigr)\exp\biggl\{-4\log (1/\theta)\sum_{j\ge0} \frakd_j\biggr\}.
\end{equation}
The bound \eqref{E:3.37} then gives
\begin{equation}
\label{E:3.50u}
\tilde c_{ k+1}\le\frac{\theta\,\tilde c_{ k}}{1+\bar\alpha\tilde c_{ k}^2}+(b-1)\theta\tilde c_{ k}
\end{equation}
once~$k$ and~$n-k$ exceed some~$k_0\ge1$. 

 Now observe that setting $\theta:=1/b$ reduces \eqref{E:3.50u}  to the conclusion of \cite[Lemma~4.6]{BHu}. The proof of \cite[Theorem~3.5]{BHu} then applies, resulting in  the bound  $\tilde c_k\le C'(1+\sqrt k)^{-1/2}$. The fact that $c_k/\tilde c_k$ is bounded uniformly in $n\ge k\ge0$  by Assumption~\ref{ass-1} then gives \eqref{E:3.33} once~$\min\{k,n-k\}\ge k_0$. Thanks to~$\sup_{n\ge k\ge0}c_k<\infty$, the extension to small~$k$ is achieved by choosing~$C$ sufficiently large. The extension to~$k$ close to~$n$ is in turn supplied by Lemma~\ref{lemma-3.5} along with Assumption~\ref{ass-1}.
\end{proofsect}

The next lemma shows that the supremum in \eqref{E:3.33} is actually order $k^{-1/2}$ and, in fact, so is even the ratio $a_k(1)/a_k(0)$:

\begin{lemma}
\label{lemma-3.9}
Assume~$\beta=\beta_\cc$. There exists a constant~$C'>0$ such that for all~$n\ge k\ge0$,
\begin{equation}
\label{E:3.42}
\frac{a_k(1)}{a_k(0)}
\ge\frac{C'}{\sqrt{1+k}}.
\end{equation}
\end{lemma}

\begin{proofsect}{Proof}
Let~$c_k$ continue to denote the supremum in \eqref{E:3.33}.
We first prove suitable bounds on $a_{k+1}(1)$ and~$a_{k+1}(0)$ using \eqref{E:3.8}. Indeed, neglecting all but the terms with just one $\ell_i$ non-zero in  \eqref{E:3.8}  yields  the  lower bound
\begin{equation}
\label{E:3.43}
a_{k+1}(1)\ge b\theta_{k+1} a_k(1) a_k(0)^{b-1}.
\end{equation}
For an upper bound on $a_{k+1}(0)$, we note that the sum in \eqref{E:3.8} contains a term with all zeros, then terms with exactly two indices equal to $\pm1$ (and all other equal to zero), and then the remaining terms in which either three indices are non-zero or two are non-zero but at least one is at least two in absolute value. Invoking the bound $a_k(\ell)\le c_k^{|\ell|} a_k(0)$ while assuming, thanks to Lemma~\ref{lemma-3.8}, that~$k$ is so large that $c_k\le1/2$, the contribution of the last two cases is estimated by
\begin{equation}
\bigl[b(b-1)+b(b-1)(b-2)\bigr]c_k^3\biggl(1+2\sum_{m\ge1}c_k^m\biggr)^b a_k(0)^b.
\end{equation}
 Here  the prefactors dominate the number of ways the  above  sets of  indices can appear in the given ordering of the $b$-tuple $\ell_1,\dots,\ell_b$ in \eqref{E:3.8} and the term in the large parentheses  dominates the sum over the remaining indices after   the restriction on $\ell_1+\dots+\ell_b$  has been dropped.  Hence we get
\begin{equation}
\label{E:3.45}
a_{k+1}(0)\le a_k(0)^b+b(b-1) a_k(1)^2 a_k(0)^{b-2}+ \alpha' c_k^3a_k(0)^b,
\end{equation}
where~$\alpha':=b(b-1)^2 3^b$. Abbreviating
\begin{equation}
\label{E:3.46}
\lambda_k:=\frac{a_k(1)}{a_k(0)},
\end{equation}
from \eqref{E:3.43} and \eqref{E:3.45} we then get
\begin{equation}
\label{Eq:3.55}
\lambda_{k+1}\ge \frac{b\theta_{k+1}\lambda_k}{1+b(b-1)\lambda_k^2+ \alpha' c_k^3}
\ge\frac{b\theta_{k+1}}{1+\alpha' c_k^3}\,\frac{\lambda_k}{1+b(b-1)\lambda_k^2}.
\end{equation}
At~$\beta=\beta_\cc$ we have~$b\theta_{k+1} = b^{1-\sigma_{k+1}^2}$.  Denote  
\begin{equation} 
\tilde\lambda_k:= \lambda_k\prod_{j=1}^k \frac{1+\alpha' c_{j-1}^3}{b \theta_j}  = \lambda_k\prod_{j=1}^k \frac{1+\alpha' c_{j-1}^3}{b^{1-\sigma_j^2}}, 
\end{equation} 
 and observe  that $\lambda_k /\tilde \lambda_k$ is bounded from above and below uniformly in $n\ge k \ge 1$ thanks to Assumption~\ref{ass-1} and $\prod_{j\ge 1} \left(1 + j^{-3/2}\right) < \infty$. Abbreviating $r:=b(b-1)\sup_{n\ge k\ge0}\lambda_k/\tilde\lambda_k$, the inequality \eqref{Eq:3.55}  then gives 
\begin{equation}
\tilde\lambda_{k+1}\ge \frac{\tilde\lambda_k}{1+r\tilde\lambda_k^2}.
\end{equation}
This now readily yields $\tilde\lambda_{k+1}^{-2}\le\tilde\lambda_k^{-2}+r'$ for  $k \ge 1$ and ~$r':=2r+r^2 \sup_{n\ge k\ge1}\tilde\lambda_k^2  < \infty$.  Hence we get that $\tilde\lambda_k^{-2}/k$ is bounded uniformly in~$n\ge k\ge k_0$, for some $k_0\ge1$. This, along with  $\lambda_k/\tilde\lambda_k$ being bounded from below  gives \eqref{E:3.42} for~$k$ sufficiently large.  The extension to $k < k_0$ is routine from the positivity of $a_k(1)/a_k(0)$ which by Assumption~\ref{ass-1} and \eqref{Eq:3.55} holds uniformly in~$n$ satisfying~$n \ge k$. 
\end{proofsect}

\subsection{Proof of Theorem~\ref{thm-4}}
Moving to the proof of Theorem~\ref{thm-4}, we now cast the iterations of~$\{\lambda_k\}_{k\ge0}$ from \eqref{E:3.46} in a form that is amenable to asymptotic analysis. This will require tracking also the second-order  Fourier coefficients  in the form
\begin{equation}
\rho_k:=\frac{a_k(2)}{a_k(0)}
\end{equation}
which, as suggested by \eqref{E:3.33}  and \eqref{E:3.42},  decays proportionally to~$\lambda_k^2$.
Here we get:

\begin{lemma}
\label{lemma-3.12}
Assume~$\beta=\beta_\cc$. For each~$n\ge1$ there exist positive sequences $\{r_k\}_{k=0}^{n-1}$, $\{s_k\}_{k=0}^{n-1}$ and $\{t_k\}_{k=0}^{n-1}$ that are bounded uniformly in~$n$ such that
\begin{equation}
\label{Eq:1.47b}
\begin{aligned}
\lambda_{k+1}&=b\theta_{k+1}\frac{\lambda_k+ (b-1)\lambda_k\rho_k+\frac12(b-1)(b-2)\lambda_k^3+r_k \lambda_k^4}{1+b(b-1)\lambda_k^2+s_k\lambda_k^3}
\\*[3mm]
\rho_{k+1}&=b^4\theta_{k+1}^4\frac{b^{-3}\rho_k+\frac12 b^{-3}(b-1)\lambda_k^2+ t_k \lambda_k^3}{1+b(b-1)\lambda_k^2+s_k\lambda_k^3}
\end{aligned}
\end{equation}
hold true for all~$n>k\ge 0$,   where $\lambda_k$ is as in \eqref{E:3.46}. 
\end{lemma}

\begin{proofsect}{Proof}
Note that, by Lemmas~\ref{lemma-3.8}--\ref{lemma-3.9}, when $\beta=\beta_\cc$ the quantity~$\lambda_k$ is proportional to~$c_k$ while $a_k(q)/a_k(0)$ is bounded by~$c_k^{|q|}$. Using this we now write \eqref{E:3.45} as equality
\begin{equation}
\label{E:3.52}
\frac{a_{k+1}(0)}{a_k(0)^b}=1+b(b-1)\lambda_k^2+s_k\lambda_k^3,
\end{equation}
where~$\{s_k\}_{k=0}^{n-1}$ is a positive sequence that is bounded uniformly in~$n\ge1$. Similarly, since the only  integer-valued $b$-tuples $(\ell_1,\dots,\ell_b)$ with $\ell_1+\dots+\ell_b=1$ and $|\ell_1|+\dots+|\ell_b|<4$ are permutations of $(1,0,\dots,0)$, $(-1,2,0,\dots,0)$ and $(1,1,-1,0,\dots,0)$, we get
\begin{equation}
\label{E:3.53}
\frac{a_{k+1}(1)}{a_k(0)^b} = \theta_{k+1}\Bigl(b\lambda_k+b(b-1)\lambda_k\rho_k+\tfrac12 b(b-1)(b-2)\lambda_k^3+  b r_k\lambda_k^4  \Bigr)
\end{equation}
for a positive sequence $\{r_k\}_{k=0}^{n-1}$ that is bounded uniformly in~$n\ge1$. Dividing \eqref{E:3.53} by \eqref{E:3.52} then yields the first equality in \eqref{Eq:1.47b}.

For the second equality in \eqref{Eq:1.47b} we first observe that the only integer-valued $b$-tuples $(\ell_1,\dots,\ell_b)$ with $\ell_1+\dots+\ell_b=2$ and  $|\ell_1|+\dots+|\ell_b|<4$  are permutations of $(2,0,\dots,0)$ and $(1,1,0,\dots,0)$. This implies
\begin{equation}
\label{E:3.66ii}
\frac{a_{k+1}(2)}{a_k(0)^b} = \theta_{k+1}^4\bigl(b\rho_k+\tfrac12 b(b-1)\lambda_k^2+b^4 t_k\lambda_k^4\bigr)
\end{equation}
for a positive sequence $\{t_k\}_{k=0}^{n-1}$ that is bounded uniformly in~$n\ge1$. Dividing this by \eqref{E:3.52} then gives the desired claim.
\end{proofsect}

We are now finally ready to give:

\begin{proofsect}{Proof of Theorem~\ref{thm-4}} 
Assume~$\beta=\beta_\cc$. Let us start with the bounds \twoeqref{E:3.15}{E:3.17}. First, using \eqref{E:3.33} in \eqref{E:3.25} and iterating yields
\begin{equation}
\label{E:3.26a}
\frac{a_k(q)}{a_k(0)}\le\theta_k^{q^2}\Bigl(\,\frac{bC}{\sqrt{1+k}}\Bigr)^{|q|}
\end{equation}
 whenever~$k\ge1$. (For~$k=0$ we note that  $a_0(q)/a_0(0)\le\theta_0^{q^2} \hat c$, for $\hat c$ as in Lemma~\ref{lemma-3.6}.)  
Since Assumption~\ref{ass-1} implies that $\theta_k^{q^2/2}(bC)^{|q|}$ is bounded uniformly in~$n\ge k\ge0$ and $q\in\Z$, this  is sufficient for  \eqref{E:3.15}.  The bound \eqref{E:3.16} is then proved using the same argument as in~\cite[Theorem 3.4]{BHu}. For the bound \eqref{E:3.17}, we plug \eqref{E:3.26a} in \eqref{E:3.30} using the same argument as for the subcritical case.

The main point of the proof is the asymptotic  \eqref{Eq:3.18u} and  \eqref{Eq:1.44}.  For the latter we isolate  the terms $|q| = 1$ from the rest of the sum in \eqref{E:3.30}  to get 
\begin{equation}
v_k'(z) = 4\pi \texte^{v_k(z)} a_k(0) \lambda_k \sin(2\pi z) - \sum_{|q| \ge 2} \texte^{v_k(z)} a_k(0) 2\pi i q \, \texte^{2\pi i qz} \,\frac{a_k(q)}{a_k(0)}.
\end{equation}
 Proceeding  as in \eqref{E:3.29} using \eqref{E:3.15} instead of \eqref{E:3.12}, we get uniform convergence of $a_k(0)^{-1} \texte^{-v_k(z)}$ to 1 with decay rate $(1+k)^{-1/2}$, which leads to 

\begin{equation}
\label{Eq:1.47}
\bigl| v_k'(z) -  4\pi\lambda_k\sin(2\pi z)\bigr|\le\frac Ck
\end{equation}
for some constant~$C\in(0,\infty)$ uniformly in~$z\in\R$.
 To get \eqref{Eq:1.44}  it thus suffices to prove the asymptotic~\eqref{Eq:3.18u}. 

 We will prove \eqref{Eq:3.18u} by iterating the top line in \eqref{Eq:1.47b} but for that we  first have to show that $\rho_k/\lambda_k^2$ is, for~$k$ large, close to a constant.
Here \eqref{Eq:1.47b}  expresses $\rho_{k+1}/\lambda_{k+1}^2$ as
\begin{equation}
(b\theta_{k+1})^2\Bigl(b^{-3}\frac{\rho_k}{\lambda_k^2}+\frac12 b^{-3}(b-1)+t_k\lambda_k\Bigr)\frac{1+b(b-1)\lambda_k^2+s_k\lambda_k^3}{\left(1+(b-1)\rho_k+\frac12(b-1)(b-2)\lambda_k^2+r_k\lambda_k^3\right)^2} 
\end{equation}
which in light of $\rho_k/\lambda_k^2$ being bounded from above gives 
\begin{equation}
\label{E:3.57}
\frac{\rho_{k+1}}{\lambda_{k+1}^2}
=\bigl(  (b\theta_{k+1})^2  +t'_k\lambda_k\bigr)\biggl(b^{-3}\frac{\rho_k}{\lambda_k^2}+\frac12 b^{-3}(b-1)\biggr)
\end{equation}
for a sequence~$\{t_k'\}_{k=0}^{n-1}$ that is bounded uniformly in~$n\ge1$.  Since  the prefactor is close to~$1$  for $\min\{k,n-k\}$ large, set  $\delta_k:=(b\theta_{k+1})^2+t'_{k}\lambda_{k}-1$  and observe that then
\begin{equation}
\biggl|\frac{\rho_{k+1}}{\lambda_{k+1}^2}-\frac12\frac{b-1}{b^3-1}\biggr|
\le \frac{1+\delta_k}{b^3}\biggl|\frac{\rho_{k}}{\lambda_{k}^2}-\frac12\frac{b-1}{b^3-1}\biggr|+\delta_k.
\end{equation}
Iteration shows
\begin{equation}
\label{E:3.58}
\biggl|\frac{\rho_{k}}{\lambda_{k}^2}-\frac12\frac{b-1}{b^3-1}\biggr|
\le\sum_{j=1}^k \delta_{k-j}  \prod_{i=1}^{j-1}  \frac{1+\delta_{k-i}}{b^3}+  \biggl(\,\prod_{j=0}^k\frac{1+\delta_{j}}{b^3}\biggr)  \biggl| \frac{\rho_{0}}{\lambda_{0}^2} -\frac12\frac{b-1}{b^3-1}\biggr| .
\end{equation}
 Using  $\lambda_k=O(k^{-1/2})$ and  $b\theta_k-1=O(|\sigma_k^2-1|)$  we get $\delta_k = O(k^{-1/2})+ O(|\sigma_{ k+1}^2 - 1|)$  and so $1+\delta_k\le\exp\{c k^{-1/2}+c|\sigma_{k+1}^2 - 1|\}$ for some constant~$c>0$. It follows that, for some constants~$c',c''>0$,
\begin{equation}
\prod_{i=1}^{j-1} \frac{1+\delta_{k-i}}{b^3}\le  b^{-3(j-1)}  \exp\biggl\{ c\sum_{i=0}^{j-1}(k-i)^{-1/2}+c\sum_{i=0}^n|\sigma_i^2-1|\biggr\}
\le c' b^{-3(j-1)} \exp\{c'' j/k^{1/2}\},  
\end{equation}
where in the second inequality we invoked  Assumption~\ref{ass-1} for the second sum and noticed that the first sum is bounded by a constant times~$j/k^{1/2}$.  The product is thus checked to decay at least as $O(b^{-2j})$ which then allows us to simplify \eqref{E:3.58} as 
\begin{equation}
\label{E:3.58a}
\frac{\rho_{k}}{\lambda_{k}^2}=\frac12\frac{b-1}{b^3-1}+O(k^{-1/2}) +O\biggl(\,\sum^k_{j= 0}  b^{-2j}  |\sigma_{k-j}^2-1|\biggr).
\end{equation}
Here the implicit constants are uniform in~$n\ge k\ge0$.

Moving to the proof of \eqref{Eq:3.18u}, we temporarily denote
\begin{equation}
\tilde\lambda_k:=\lambda_k\prod_{j=0}^k  (b\theta_j)^{-1} 
\end{equation}
and observe that the first line in \eqref{Eq:1.47b} can concisely be written as
\begin{equation}
\label{E:3.60}
\tilde\lambda_{k+1}= \frac{\tilde\lambda_k}{\sqrt{1+u_k\lambda_k^2}},
\end{equation}
where
\begin{equation}
\label{E:3.66k}
u_k :=\frac1{\lambda_k^2}\Biggl[\biggl(\frac{1+b(b-1)\lambda_k^2+s_k\lambda_k^3}{1+ (b-1)\rho_k+\frac12(b-1)(b-2)\lambda_k^2+r_k \lambda_k^3}\biggr)^2-1\Biggr].
\end{equation}
The reason for writing the iteration this way is because \eqref{E:3.60} can now be cast as
\begin{equation}
\frac1{\tilde\lambda_{k+1}^2} = \frac1{\tilde\lambda_k^2}+u_k\prod_{j=0}^k  (b\theta_j)^{2} .
\end{equation}
Iterating we then get
\begin{equation}
\label{E:3.63}
\frac1{\lambda_k^2}=\frac1{\lambda_0^2}\biggl(\,\prod_{j=1}^k  (b\theta_j)^{-2}  \biggr)+\sum_{\ell=0}^{k-1}\biggl(\,\prod_{j=\ell+1}^k  (b\theta_j)^{-2} \biggr)\rho_\ell,
\end{equation}
where we already returned to the original variables.

In order to extract the leading asymptotic of the sum in \eqref{E:3.63}, we note that \eqref{E:3.58a} along with $\lambda_k=O(k^{-1/2})$ show
\begin{equation}
\label{E:3.69k}
u_k = \frac{(b-1)^2(b+1)^3}{b^3-1}+O(k^{-1/2}) +O\biggl(\,\sum^k_{j=0}  b^{-2j} |\sigma_{k-j}^2-1|\biggr). 
\end{equation}
 Denote by~$u_\star$ the first quantity on the right. Plugging \eqref{E:3.69k} in \eqref{E:3.63} while noting that, thanks to Assumption~\ref{ass-1}, $\prod_{j=\ell+1}^k(b\theta_j)^{-2}$ is bounded uniformly in~$n\ge k>\ell\ge0$ we get
\begin{equation}
\label{E:3.83ii}
\lambda_k^{-2} = O(1)+ u_\star\sum_{\ell=0}^{k-1}\prod_{j=\ell+1}^k(b\theta_j)^{-2}  + O(k^{1/2})+O\biggl(\,\sum_{j=0}^{k-1}|\sigma_j^2-1|\biggr).
\end{equation}
The last term is again $O(1)$ by Assumption~\ref{ass-1} so we just need to control the middle term. Here we separate the terms with~$\ell\le\sqrt k$ at the cost of another~$O(k^{1/2})$ correction. In the remaining terms we note that
\begin{equation}
\prod_{j=\ell+1}^k(b\theta_j)^{-2} = \prod_{j=\ell+1}^kb^{2(\sigma_j^2-1)}
=1+O\biggl(\,\sum_{j\ge\sqrt k}\frakd_j\biggr)+O\biggl(\,\sum_{j\ge n-k}\frakd_j\biggr),
\end{equation}
where~$\{\frakd_j\}_{j\ge0}$ is the sequence from Assumption~\ref{ass-1}. Using this in \eqref{E:3.83ii} and inverting the two negative powers then proves \eqref{Eq:3.18u}.
Plugging this in \eqref{Eq:1.47} gives also \eqref{Eq:1.44}.
\end{proofsect}

Unlike Theorems~\ref{thm-3}--\ref{thm-4} whose proofs borrowed from results proved in our earlier work, Theorem~\ref{thm-5} will have to be proved from ``scratch'' using different ideas. The details will be given in  Section~\ref{sec-6}. 

\section{Asymptotic covariance structure}
\label{sec-4}\noindent
We are now ready to commence the actual proof of Theorem~\ref{thm-1}. We rely heavily
on the convergence of the renormalization group iterations established in the earlier sections along with the representation of the field  as a  tree-indexed  Markov chain. We again suppose that  Assumptions~\ref{ass-1}--\ref{ass-2}  hold throughout.

\subsection{Markov-chain representation} 
In Lemma~\ref{lemma-3.2} we showed that $P_{n,\beta}$ is the law on the leaves of a tree-indexed Markov chain. Along any branch of the tree, that tree-indexed Markov chain is just an ordinary Markov chain with transition probabilities \eqref{E:1.23}. As it turns out, the proof of Theorem~\ref{thm-1} can be reduced to properties of this chain.

In order to match the labeling of the~$v_k$'s, we will label the Markov chain backwards; i.e., from the leaves to the root. A run of the chain is thus a sequence of real-valued random variables $\{\phi_k\}_{k=n,\dots,0}$  such that, for all Borel~$A\subseteq\R$, 
\begin{equation}
\label{Eq:1.2b}
P(\phi_n\in A) = \frakp_n(A\,|\,0\,)
\end{equation}
and
\begin{equation}
\label{Eq:1.2a}
P(\phi_k\in A\,|\FF_{k+1}) = \frakp_{k}(A|\phi_{k+1}),
\end{equation}
where $\frakp_k(\cdot\,|\,\cdot)$  is  as in \eqref{E:1.23} and
\begin{equation}
\label{Eq:4.1}
\FF_k:=\sigma(\phi_i\colon i=k,\dots,n).
\end{equation}
 The above is consistent with setting the ``initial'' state of the chain to $\phi_{n+1}:=0$ and using~$\FF_{n+1}$ to denote  the trivial $\sigma$-algebra.

We will use~$E$ to denote expectation with respect to~$P$  and, for~$x,y\in\Lambda_n$, write $k(x,y)$ for the depth of the nearest common ancestor of~$x$ and~$y$ in the underlying $b$-ary tree structure on~$\T_n$.  Explicitly, for $x=(x_1,\dots,x_n)$ and $y=(y_1,\dots,y_n)$ in~$\Lambda_n$ we set 
\begin{equation}
\label{Eq:4.3i}
 k(x,y):=  \min\{\,j\in \{0,\dots, n\}: \,x_i = y_i \,\forall\, i = 1,\dots, n- j\,\}. 
\end{equation}
 Note that $k(x,y) = \log_{b^{1/2}} d(x,y)$ for distinct $x,y \in \Lambda_n$.  To reduce clutter we sometimes use the same letter for both the field values and the states of the chain, as the precise meaning will always be clear from context. 

A key step in the proof of Theorem~\ref{thm-1} then comes in:

\begin{proposition}
\label{prop-4.1}
Let~$\beta>0$ and assume that~$\{v_k'\}_{k=0}^n$ and~$\{v_k''\}_{k=0}^n$ are bounded uniformly in~$n\ge1$.
Then for all~$n\ge1$ and all $x,y\in\Lambda_n$,
\begin{equation}
\label{Eq:1.4a}
\langle\phi_x\phi_y\rangle_{n,\beta}
=\sum_{i=k(x,y)}^n\biggl(\frac1\beta+\frac1{\beta^2}\frac{b+1}{b-1}\,E\bigl[v_{i}'(\phi_{i+1})^2-v_{i}''(\phi_{i+1})\bigr]\biggr) +O(1)
\end{equation}
holds with~$O(1)$ that is bounded uniformly in~$n\ge1$ and~$x,y\in\Lambda_n$.
\end{proposition}

To prove this proposition note that, for~$k:=k(x,y)$, the tree-indexed Markov chain representation gives
\begin{equation}
\label{Eq:4.5i}
\langle\phi_x\phi_y\rangle_{n,\beta} = E\bigl(E(\phi_0|\FF_k)^2\bigr)
\end{equation}
and  so it suffices to extract the asymptotic form of $E(\phi_0|\FF_k)$, uniformly in $n\ge k\ge0$.
This leads to somewhat lengthy calculations for the underlying Markov chain, part of which we relegate to:

\begin{lemma}
Let~$\beta>0$  and $n\ge1$. Then for  all  $k=0,\dots,n$,
\begin{equation}
\label{Eq:4.5}
E(\phi_k|\FF_{k+1}) = \phi_{k+1}-\frac{\sigma_{k}^2}\beta v'_{k}(\phi_{k+1})
\end{equation}
and, for all $k=1,\dots,n$,
\begin{equation}
\label{Eq:4.6}
E\bigl(v_{k-1}'(\phi_k)\big|\,\FF_{k+1}\bigr)=\frac1b v'_{k}(\phi_{k+1}).
\end{equation}
\end{lemma}

\begin{proofsect}{Proof}
Abbreviate $\beta_k:=\beta\sigma_k^{-2}$ and recall our notation~$\mu_{\sigma^2}$ for the law of~$\NN(0,\sigma^2)$.
We start by computing some relevant Gaussian integrals. The first one of these is
\begin{equation}
\label{Eq:1.3}
\begin{aligned}
\int\texte^{-b v_{k-1}(\phi+\zeta)}
\zeta \,\mu_{1/\beta_{k}}(\textd\zeta)&=\frac1{\beta_{k}}\frac{\textd}{\textd\phi}\int\texte^{-b v_{k-1}(\phi+\zeta)}\mu_{1/\beta_{k}}(\textd\zeta)
\\
&=\frac1{\beta_{k}} \frac{\textd}{\textd\phi}\texte^{-v_{k}(\phi)} = -\frac1{\beta_k} v'_{k}(\phi)\texte^{-v_{k}(\phi)}
\end{aligned}
\end{equation}
for~$k=1,\dots,n$. Here we first performed  a  Gaussian integration by parts, then swapped differentiation with respect to~$\zeta$ by that with respect to~$\phi$ and finally applied the recursive relation between~$v_k$ and~$v_{k-1}$. The second integral to compute uses the relation between~$v_k$ and~$v_{k-1}$ with the result
\begin{equation}
\label{Eq:1.5a}
\int\texte^{-b v_{k-1}(\phi+\zeta)}
v_{k-1}'(\phi+\zeta) \,\mu_{1/\beta_{k}}(\textd\zeta) = -\frac1b \frac{\textd}{\textd\phi}\texte^{-v_{k}(\phi)} = \frac1b v_{k}'(\phi)\texte^{-v_{k}(\phi)}
\end{equation}
for all~$k=1,\dots,n$.
 The  above manipulations are justified by the fact that the $v_k$'s are periodic $C^\infty$-functions and so no issues arise from swapping derivatives and integrals and no boundary terms pop up during integration by parts.

Moving to statements above, for \eqref{Eq:4.6} we only need to use the top line in the definition~\eqref{E:1.23} of the transition probability  to get  
\begin{equation}
\label{Eq:1.5}
E\bigl(v_{k-1}'(\phi_k)\,\big|\,\FF_{k+1}\bigr)
= \texte^{v_{k}(\phi_{k+1})}\int\texte^{-b v_{k-1}(\phi_{k+1}+\zeta)}v_{k-1}'(\phi_{k+1}+\zeta)\mu_{1/\beta_{k}}(\textd\zeta)
\end{equation}
Using \eqref{Eq:1.5a} we now  obtain  \eqref{Eq:4.6}. For \eqref{Eq:4.5} we first treat~$k\ge1$ where \eqref{Eq:1.3}  combined  with the top line in~\eqref{E:1.23} show
\begin{equation}
\begin{aligned}
E(\phi_k|\FF_{k+1}) 
&= \texte^{v_{k}(\phi_{k+1})}\int\texte^{-b v_{k-1}(\phi_{k+1}+\zeta)}(\phi_{k+1}+\zeta)\mu_{1/\beta_{k}}(\textd\zeta)
\\
&=\phi_{k+1}-\frac1{\beta_{k}} v_{k}'(\phi_{k+1}).
\end{aligned}
\end{equation}
For~$k=0$ we instead use the bottom line in \eqref{E:1.23} to get
\begin{equation}
\begin{aligned}
E(\phi_0|\FF_1) 
&= \phi_1+\texte^{v_0(\phi_1)}\int \texte^{-\frac{\beta_0}2(\phi-\phi_1)^2}(\phi-\phi_1)\nu(\textd\phi)
\\
&=\phi_1+\texte^{v_0(\phi_1)}\frac1{\beta_0}\frac{\textd}{\textd\phi_1}\int \texte^{-\frac{\beta_0}2(\phi-\phi_1)^2}\nu(\textd\phi)
\\
&=\phi_1-\frac1{\beta_0}v_0'(\phi_1)
\end{aligned}
\end{equation}
by invoking the definition of~$v_0$ in the last step. 
\end{proofsect}

\begin{remark}
We note that \eqref{Eq:4.6} shows that, for the model corresponding to the renormalization fixed point, the associated potential is an eigenvector of the Markov transition kernel restricted to the space of $1$-periodic functions. This fact was key for the derivations in Benfatto and Renn; see~\cite[Eq.~4.22]{BR} and thereafter.
\end{remark}

As a consequence of the identities \twoeqref{Eq:4.5}{Eq:4.6} we then get:

\begin{corollary}
The sequence $\{M_k\}_{k=0}^{ n+1}$ defined by $M_0:=\phi_0$ and by
\begin{equation}
\label{Eq:4.12}
M_k:=\phi_k-\frac1\beta \biggl(\,\sum_{i=1}^{k} b^{1-i}\sigma_{k-i}^2\biggr)v_{k-1}'(\phi_k), \quad k=1,\dots, n+1, 
\end{equation}
is a reverse martingale; i.e., $E(M_k|\FF_{k+1})=M_{k+1}$ holds true a.s. for all~$k=0,\dots,n$.
\end{corollary}

\begin{proofsect}{Proof}
Write \eqref{Eq:4.12} as $M_k:=\phi_k-s_k\beta^{-1}v_{k-1}'(\phi_k)$ where~$s_0:=0$ to avoid having to interpret~$v_{-1}'$. For $\{M_k\}_{k=0}^{ n+1}$ to be a reverse martingale, the identities \twoeqref{Eq:4.5}{Eq:4.6} dictate that $s_{k+1}=b^{-1}s_k+\sigma_k^2$ for all~$k=0,\dots,n$. This is satisfied by $s_k:=\sum_{i=1}^{k} b^{1-i}\sigma_{k-i}^2$.
\end{proofsect}

Let us keep writing $s_k:=\sum_{i=1}^{k} b^{1-i}\sigma_{k-i}^2$. The fact that~$\{M_k\}_{k=0}^{ n+1}$ is a martingale with $M_0=\phi_0$ allows us to continue \eqref{Eq:4.5i} as
\begin{equation}
\label{Eq:1.7}
\begin{aligned}
\langle\phi_x\phi_y\rangle_{n,\beta} 
&= E\bigl(E(\phi_0|\FF_k)^2\bigr)
\\&= E(M_k^2)= E(\phi_k^2)-\frac{2s_k}\beta E\bigl(\phi_k v_{k-1}'(\phi_k)\bigr)+\frac{s_k^2}{\beta^2}E(v_{k-1}'(\phi_k)^2\bigr).
\end{aligned}
\end{equation}
To compute the expectation on the  left,  we  need to  iteratively compute the  expectations  arising on the right.  This is done in: 

\begin{lemma}
Abbreviate $\beta_k:=\beta\sigma_k^{-2}$. For all~$k=0,\dots,n$,
\begin{equation}
\label{Eq:1.12}
\begin{aligned}
E\bigl(\phi_k^2\,\big|\,\FF_{k+1}\bigr) =\phi_{k+1}^2+\frac1{\beta_k}+\frac1{\beta_k^2}\bigl[v_{k}'(\phi_{k+1})^2-v_{k}''(\phi_{k+1})\bigr]-\frac2{\beta_k}\phi_{k+1}v_{k}'(\phi_{k+1})
\end{aligned}
\end{equation}
and, for all $k=1,\dots,n$, also
\begin{equation}
\label{Eq:4.16u}
E\bigl(\phi_kv_{k-1}'(\phi_k)\,\big|\,\FF_{k+1}\bigr) = \frac1b\phi_{k+1}v_{k}'(\phi_{k+1})-\frac1{\beta_k}\,\frac1b\bigl[v_{k}'(\phi_{k+1})^2-v_{k}''(\phi_{k+1})\bigr].
\end{equation}
\end{lemma}

\begin{proofsect}{Proof}
We again start by computing some relevant integrals. The first of these uses similar ideas as \eqref{Eq:1.3} with the result
\begin{equation}
\label{Eq:1.4}
\begin{aligned}
\int\texte^{-b v_{k-1}(\phi+\zeta)}
\zeta^2 \,\mu_{1/\beta_{k}}(\textd\zeta)&=\frac1{\beta_{k}} \texte^{-v_{k}(\phi)}+\frac1{\beta_{k}}\frac{\textd}{\textd\phi}\int\texte^{-b v_{k-1}(\phi+\zeta)}\zeta\,\mu_{1/\beta_{k}}(\textd\zeta)
\\
&=\frac1{\beta_{k}} \texte^{-v_{k}(\phi)}+\frac1{\beta_k^2}\frac{\textd^2}{\textd\phi^2}\texte^{- v_{k}(\phi)} \\
&= \Bigl(\frac1{\beta_{k}}+ \frac1{\beta_{k}^2}\bigl[v'_{k}(\phi)^2-v''_{k}(\phi)\bigr]\Bigl)\texte^{-v_{k}(\phi)}
\end{aligned}
\end{equation}
for all $k=1,\dots,n$. Here we first interpreted one of the~$\zeta$'s as a term coming from the derivative of the probability density of~$\mu_{1/\beta_{k}}$ and then used that to integrate by parts, which reduces the computation to the integral in \eqref{Eq:1.3}.
 Differentiating twice the formula for~$\texte^{-v_0}$ in turn shows 
\begin{equation}
\label{Eq:4.16}
\begin{aligned}
\int \texte^{-\frac{\beta_0}2(\phi-\phi_1)^2}&(\phi-\phi_1)^2\nu(\textd\phi)
\\
&=\frac1{\beta_0}
\int \texte^{-\frac{\beta_0}2(\phi-\phi_1)^2}\nu(\textd\phi)
+\frac1{\beta_0^2}\frac{\textd^2}{\textd\phi^2}\int \texte^{-\frac{\beta_0}2(\phi-\phi_1)^2}\nu(\textd\phi)
\\
&=\Bigl(\frac1{\beta_0}+\frac1{\beta_0^2}\bigl[v'_{0}(\phi)^2-v''_{0}(\phi)\bigr]\Bigl)\texte^{-v_{0}(\phi)}.
\end{aligned}
\end{equation}
To get \eqref{Eq:1.12} we  now  invoke $\phi_k^2 = \phi_{k+1}^2+(\phi_{k}-\phi_{k+1})^2+2\phi_{k+1}(\phi_k-\phi_{k+1})$ and then apply \twoeqref{Eq:1.4}{Eq:4.16} to the second term and \eqref{Eq:4.5} to the third term.

The proof of \eqref{Eq:4.16u} again starts by  computing an integral; namely, 
\begin{equation}
\label{Eq:1.5b}
\begin{aligned}
\int&\texte^{-b v_{k-1}(\phi+\zeta)}
v_{k-1}'(\phi+\zeta) \zeta\,\mu_{1/\beta_{k}}(\textd\zeta) 
\\
&= -\frac1b \frac{\textd}{\textd\phi}\int\texte^{-b v_{k-1}(\phi+\zeta)}
\zeta\,\mu_{1/\beta_{k}}(\textd\zeta)  = -\frac1{\beta_{k}}\,\frac1b \bigl[v_{k}'(\phi)^2-v_{k}''(\phi)\bigr]\texte^{-v_{k}(\phi)},
\end{aligned}
\end{equation}
 where the last equality  follows by plugging in an intermediate step from \eqref{Eq:1.4}. This shows that, for all $k=1,\dots,n$,
\begin{equation}
\label{Eq:1.13}
\begin{aligned}
E\bigl(\phi_kv_{k-1}'(\phi_k)\,\big|\,\FF_{k+1}\bigr) &= \texte^{v_k(\phi_{k+1})} \int\texte^{-b v_{k-1}(\phi_{k+1}+\zeta)}
\,(\phi_{k+1}+\zeta)\,v_{k-1}'(\phi_{k+1}+\zeta)\,\mu_{1/\beta_{k}}(\textd\zeta) \\
&= \frac1b\phi_{k+1}v_{k}'(\phi_{k+1})-\frac1{\beta_k}\,\frac1b\bigl[v_{k}'(\phi_{k+1})^2-v_{k}''(\phi_{k+1})\bigr],
\end{aligned}
\end{equation}
where the first term arises via \eqref{Eq:1.5a} and the second term via \eqref{Eq:1.5b}.
\end{proofsect}

We are now ready to give:

\begin{proofsect}{Proof of Proposition~\ref{prop-4.1}}
We will keep using the above shorthands~$s_k$ and~$\beta_k$ whenever convenient.  The argument aims  directly at an iterative computation of~$E(M_k^2)$.  Here 
the identities \twoeqref{Eq:1.12}{Eq:4.16u} give
\begin{equation}
\label{E:4.22ui}
\begin{aligned}
E\bigl(M_k^2\,\big|\,\FF_{k+1}\bigr)&=E\biggl(\Bigl[\phi_k-\frac{s_k}\beta v_{k-1}'(\phi_k)\Bigr]^2\,\bigg|\,\FF_{k+1}\biggr)
\\
&=\phi_{k+1}^2+\frac1{\beta_k}+\frac1{\beta_k^2}\bigl[v_{k}'(\phi_{k+1})^2-v_{k}''(\phi_{k+1})\bigr]-\frac2{\beta_k}\phi_{k+1}v_{k}'(\phi_{k+1})
\\
&\quad-\frac{2s_k}\beta\biggl(\frac1b\phi_{k+1}v_{k}'(\phi_{k+1})-\frac1{\beta_k}\,\frac1b\bigl[v_{k}'(\phi_{k+1})^2-v_{k}''(\phi_{k+1})\bigr]\biggr)
\\
&\qquad+\frac{s_k^2}{\beta^2}E\bigl(v_{k-1}'(\phi_k)^2\,\big|\,\FF_{k+1}\bigr).
\end{aligned}
\end{equation}
Next note that, by the definition of~$\beta_k$ and the recursion $s_{k+1}=b^{-1}s_k+\sigma_k^2$, the terms containing $\phi_{k+1}v_k'(\phi_{k+1})$ combine into the cross-term that arises from squaring the expression $M_{k+1}=\phi_{k+1}-s_{k+1}\beta^{-1}v_k'(\phi_{k+1})$. This  wraps \eqref{E:4.22ui} into 
\begin{equation}
\label{Eq:1.15}
\begin{aligned}
E\bigl(M_k^2\,\big|\,\FF_{k+1}\bigr) = M_{k+1}^2 +&\frac{\sigma_k^2}\beta+\frac{\sigma_k^2}{\beta^2}\Bigl(\sigma_k^2+\frac{2s_k}b\Bigr)\bigl[v_{k}'(\phi_{k+1})^2-v_{k}''(\phi_{k+1})\bigr]
\\
&+\frac1{\beta^2}\Bigl( s_k^2 E\bigl(v_{k-1}'(\phi_k)^2\,\big|\,\FF_{k+1}\bigr) -s_{k+1}^2 v_{k}'(\phi_{k+1})^2\Bigr).
\end{aligned}
\end{equation}
Taking expectation and invoking the assumed uniform boundedness of~$\{v_k'\}_{k=0}^n$ and~$\{v_k''\}_{k=0}^n$, we may replace~$\sigma_k^2$ by~$1$ and~$s_k$ by~$\frac b{b-1}$ to get
\begin{equation}
\label{Eq:4.23}
\begin{aligned}
E(M_k^2) &= E(M_{k+1}^2) + \frac1\beta+\frac1{\beta^2}\frac{b+1}{b-1}E\bigl[v_{k}'(\phi_{k+1})^2-v_{k}''(\phi_{k+1})\bigr]
\\
&\quad+\frac1{\beta^2}\Bigl(\frac b{b-1}\Bigr)^2\Bigl( E\bigl(v_{k-1}'(\phi_k)^2\bigr) -E\bigl( v_{k}'(\phi_{k+1})^2\bigr)\Bigr)
\\
&\quad+O\biggl(\,\sum_{i=0}^k b^{-i}|\sigma_{k-i}^2-1|\biggr)+O(b^{-k}),
\end{aligned}
\end{equation}
where $O(b^{-k})$ arises from the bound on the tail of the infinite series for the asymptotic value of~$s_k$.
Under Assumption~\ref{ass-1},  adding up  the error terms  over~$k=0,\dots, n$ shows that these produce  an~$O(1)$ correction under iteration. The same applies to the middle term as it leads to a telescoping sum which is bounded by the assumed  uniform  boundedness of~$\{v_k'\}_{k=0}^n$. Iterations of \eqref{Eq:4.23} then prove the desired claim.
\end{proofsect}

\subsection{Technical lemmas}
In order to process the formula in Proposition~\ref{prop-4.1} further we need a couple of technical lemmas. First we note a way to simplify the right-hand side of \eqref{Eq:1.4a}.

\begin{lemma}
\label{lemma-4.6}
Let~$\beta>0$ and assume that~$\{v_k'\}_{k=0}^n$ and~$\{v_k''\}_{k=0}^n$ are bounded uniformly in~$n\ge1$. Then for all $n\ge1$ and all $k=0,\dots,n$,
\begin{equation}
\label{Eq:1.20a}
\sum_{i=k}^n E\bigl[v_{i}'(\phi_{i+1})^2-v_{i}''(\phi_{i+1})\bigr]
=-b\sum_{i=k}^n E\bigl(v_{i}'(\phi_{i+1})^2\bigr) + O(1),
\end{equation}
where~$O(1)$ is bounded uniformly in~$n\ge k\ge 0$.
\end{lemma}

\begin{proofsect}{Proof}
 By the assumed boundedness of~$\{v_k'\}_{k=0}^n$ and~$\{v_k''\}_{k=0}^n$ it suffices to prove this for $k\ge1$. Here the definition of~$v_{k}$ from~$v_{k-1}$ gives 
\begin{equation}
\begin{aligned}
E\bigl(b^2 &v_{k-1}'(\phi_k)^2-b v_{k-1}''(\phi_k)\,\big|\,\FF_{k+1}\bigr)
\\
&=\texte^{v_{k}(\phi_{k+1})}\int\texte^{-b v_{k-1}(\phi_{k+1}+\zeta)}\bigl[b^2 v_{k-1}'(\phi_{k+1}+\zeta)^2-b v_{k-1}''(\phi_{k+1}+\zeta)\bigr]\mu_{1/\beta_k}(\textd\zeta)
\\
&=\texte^{v_{k}(\phi_{k+1})}\frac{\textd^2}{\textd\phi_{k+1}^2}\int\texte^{-b v_{k-1}(\phi_{k+1}+\zeta)}\mu_{1/\beta_k}(\textd\zeta)=v_{k}'(\phi_{k+1})^2 - v_{k}''(\phi_{k+1}).
\end{aligned}
\end{equation}
Taking expectation then shows
\begin{equation}
\sum_{i=k+1}^{n} E\bigl[v_{i}'(\phi_{i+1})^2-v_{i}''(\phi_{i+1})\bigr]
=\sum_{i=k}^{n-1}E\bigl[b^2 v_{i}'(\phi_{i+1})^2-b v_{i}''(\phi_{i+1})\bigr].
\end{equation}
 Relying again on the boundedness of $\{v_k'\}_{k=0}^n$ and~$\{v_k''\}_{k=0}^n$,  rearranging terms  yields 
\begin{equation}
(b^2-1)\sum_{i=k}^n E\bigl(v_{i}'(\phi_{i+1})^2\bigr) = (b-1)\sum_{i=k}^n E\bigl(v_{i}''(\phi_{i+1})\bigr) +O(1).
\end{equation}
Canceling $b-1$ on both sides then  shows  the claim.
\end{proofsect}

Next, in order to determine the asymptotic behavior of the sum on the right of \eqref{Eq:1.20a}, we need to control the law of~$\phi_k$ for large~$k$. While this law does not converge by itself due to the fact that the variance of~$\phi_k$ stays of order~$k$, the law of its fractional part (which is all we need to compute expectations of $1$-periodic functions) does converge as long as~$v_k$ tends to a limit. In quantitative form, this is the subject of:

\begin{lemma}
\label{lemma-1.4}
Let~$\beta>0$ and  let $v_\star$ be a $1$-periodic continuous function such that  \eqref{Eq:3.33} holds.  Let~$\nu_\star$ be the Borel measure  on~$[0,1)$  defined by
\begin{equation}
\nu_\star(\textd z):=\Bigl[\int_{[0,1]}\texte^{-(b+1)v_\star(z')}\textd z'\Bigr]^{-1}\texte^{-(b+1)v_\star(z)}1_{[0,1)}(z)\textd z.
\end{equation}
For~$k=1,\dots,n$, abbreviate
\begin{equation}
\label{Eq:4.29i}
\begin{aligned}
A_{n,k}:=(b+1)&\sup_{z\in\R}\bigl|v_\star(z)\bigr|+\sup_{z,z'\in\R}\bigl|v_k(z)-b v_{k-1}(z')\bigr|
\\&+\frac\beta2\max\{\sigma_k^{-2},1\} 
+ \frac12\log(2\pi/\beta)+\log\max\{\sigma_k,1\}
\end{aligned}
\end{equation}
and
\begin{equation}
\label{Eq:4.30}
\delta_k:=\sup_{z,z'\in\R}\Bigl| \bigl[v_k(z)-bv_{k-1}(z')\bigr] - \bigl[v_\star(z)-bv_\star(z')\bigr]\Bigr|+\log\max\{\sigma_k,\sigma_k^{-1}\}.
\end{equation}
Then for all $k=1,\dots,n$,
\begin{equation}
\label{Eq:1.28}
\begin{aligned}
\bigl\Vert P(\phi_k\text{\rm\ mod }1\in\cdot) - \nu_\star\bigr\Vert_{\text{\rm TV}}
\le \prod_{i=k}^n (1-\texte^{-A_{n,i}}) +\sum_{j=k+1}^n (1-\texte^{-\delta_{j-1}})\prod_{i=k}^{j-2}(1-\texte^{-A_{n,i}}),
\end{aligned}
\end{equation}
 where the last product is interpreted as~$1$ when~$j=k+1$. 
\end{lemma}

 Before we delve into the proof, let us explain how this lemma will be used later. Our results about the RG flow imply that the quantities~$A_{n,k}$ are uniformly bounded for all~$n\ge k\ge0$. This, in turn, ensures that the products on the right-hand side of~\eqref{Eq:1.28} decay exponentially in the number of factors. Meanwhile, the term~$\delta_k$ is small whenever $\min\{k,n-k\}$ is large. As a result, the bound~\eqref{Eq:1.28} implies that the total variation distance between the law of~$\phi_k \bmod 1$ and the stationary measure~$\nu_\star$ --- which corresponds to the law of $\phi \bmod 1$ at the RG fixed point --- is small as well.

\begin{proofsect}{Proof of Lemma~\ref{lemma-1.4}}
The proof proceeds by a coupling argument. First note that the assumptions about~$v_\star$ ensure that
\begin{equation}
\cmss P_\star(B\,|\,\phi) = \texte^{v_{\star}(\phi)}\int\texte^{-b v_\star(\phi+\zeta)}
1_B(\phi+\zeta\text{ mod }1)\mu_{1/\beta}(\textd\zeta)
\end{equation}
is a transition kernel on~$[0,1)$. Let~$(\phi_n^\star,\dots,\phi_0^\star)$ denote a run of Markov chain with transition probability $\cmss P_\star$ and~$\phi_n^\star$ drawn from~$\nu_\star$ above. Write~$P_\star$ for the distribution of the chain.  Using that $v_\star$ obeys \eqref{Eq:3.33} we now check that $\nu_\star$ is stationary for~$\cmss P_\star$. In particular, we have  $P_\star(\phi_k^\star\in\cdot)=\nu_\star$ for all~$k=1,\dots,n$.

Recall the following standard coupling of random variables~$X$ and~$Y$ taking values in~$[0,1)$ with probability densities denoted as $f$, resp., $g$:
\begin{equation}
\label{Eq:1.30}
\begin{aligned}
P\bigl((X,Y)\in B\bigr) =& \int 1_B(x,x) f\wedge g(x)\textd x 
\\
&+\int_B\frac{[f(x)-f\wedge g(x)][g(y)-f\wedge g(y)]}{1-\int f\wedge g(z)\,\textd z}\,\textd x\,\textd y,
\end{aligned}
\end{equation}
where~$f\wedge g(x):=\min\{f(x),g(x)\}$. If the pair is drawn using the first term, we will say that~$X$ and~$Y$ get ``coupled,'' while if the pair is drawn  using  the second term, we say that they get ``uncoupled.''

 Denoting $\phi_k' := \phi_k \text{\rm\ mod }1,$  we will now apply the above coupling recursively to generate a sequence
\begin{equation}
(\phi_{ n+1}',\phi_{ n+1}^\star),\dots, (\phi_1',\phi_1^\star)
\end{equation}
 of $[0,1)\times[0,1)$-valued pairs of random variables  as follows: Draw $(\phi_{ n+1}',\phi_{ n+1}^\star)$  from $\delta_0\otimes\nu_\star$. Then, given a sample of $(\phi_{k+1}',\phi_{k+1}^\star)$  for some~$k=1,\dots,n$,  draw $(\phi_{k}',\phi_{k}^\star)$ from the above coupling measure with
\begin{equation}
\label{Eq:1.31}
f(z):=\texte^{v_{k}(\phi_{k+1})-b v_{k-1}(z)}\frac1{\sqrt{2\pi/\beta_k}}\sum_{j\in\Z}\texte^{-\frac{\beta_k}2(z-\phi_{k+1}'+j)^2},
\end{equation}
where~$\beta_k:=\beta\sigma_k^{-2}$, and
\begin{equation}
\label{Eq:1.32}
g(z):=\texte^{ v_\star(\phi_{k+1}^\star)-b v_\star(z)}\frac1{\sqrt{2\pi/\beta}}\sum_{j\in\Z}\texte^{-\frac\beta2(z-\phi_{k+1}^\star+j)^2}.
\end{equation}
 As~$f$ is the probability density of $\frakp_k(\cdot\text{ mod }1|\phi_{k+1}')$ and~$g$ the  probability density $\cmss P_\star(\cdot|\phi_{k+1}^\star)$,  this gives us a (Markovian) coupling of $(\phi_{n+1}\text{ mod }1,\dots,\phi_1\text{ mod }1)$ and a run of the Markov chain with transition probability~$\cmss P_\star$ with initial law~$\nu_\star$. 

We will now use the explicit expression \eqref{Eq:1.30} to control the probabilities that $(\phi_{k}',\phi_{k}^\star)$ get ``coupled'' or get ``uncoupled.''
We start by deriving bounds on the terms entering on the right of \eqref{Eq:1.30}. First note that, for~$f$ and~$g$ as in \twoeqref{Eq:1.31}{Eq:1.32}, retaining only the~$j=0$ term in the sums shows 
\begin{equation}
f\wedge g(z)\ge \texte^{-A_{n,k}}
\end{equation}
for all~$z\in[0,1)$,  regardless of the values of~$\phi_{k+1}'$ and~$\phi_{k+1}^\star$. 
 On the other hand, assuming that $\phi_{k+1}'=\phi_{k+1}^\star$, if~$\sigma_k^2\ge1$, then $\beta_k\le\beta$ and from \eqref{Eq:4.30} we get $f(z)\ge\texte^{-\delta_k} g(z)$  for all~$z\in[0,1)$. This implies $f\wedge g(z)\ge\texte^{-\delta_k}g(z)$ leading to 
\begin{equation}
 g(z)-f\wedge g(z)\le (1-\texte^{-\delta_k})g(z) 
\end{equation}
for all $z\in[0,1)$.
Similarly,  still under $\phi_{k+1}'=\phi_{k+1}^\star$, if~$\sigma_k^2\le 1$ then $\beta\le\beta_k$ and so  $f(z)\le\texte^{\delta_k}g(z)$ for all~$z\in[0,1)$. This now gives $f\wedge g(z)\ge \texte^{-\delta_k} f(z)$ and so
\begin{equation}
f(z)-f\wedge g(z)\le(1-\texte^{-\delta_k}) f(z)
\end{equation}
holds for all~$z\in[0,1)$. 
Using the fact that~$(\phi_{k}',\phi_{k}^\star)$ get ``coupled'' with probability equal to the total ``mass'' of the first term on the right of \eqref{Eq:1.30} and get ``uncoupled'' with probability equal to the total ``mass'' of  the second part, the above observations readily translate into the inequalities
\begin{equation}
\label{Eq:1.35}
\wt P\bigl(\phi_k' = \phi_k^\star\,\big|\FF_{k+1}^\star\bigr)
\ge \texte^{-A_{n,k}}\quad\text{\rm a.s.}
\end{equation}
and
\begin{equation}
\label{Eq:1.36}
\wt P\bigl(\phi_k' \ne \phi_k^\star\,\big|\FF_{k+1}^\star\bigr)\le 
 1-\texte^{-\delta_k} \quad\text{\rm a.s.\, on }\{\phi_{k+1}'= \phi_{k+1}^\star\},
\end{equation}
where we set $\FF_k^\star:=\sigma(\phi_i',\phi_i^\star\colon i=k,\dots,n)$  and wrote~$\wt P$ for the coupling measure. 

We now observe that the event that $\phi_k' \ne \phi_k^\star$ entails one of two possibilities: either the chains never got ``coupled'' up to and including the state indexed by~$k$, or they did get ``coupled'' at some index $j=k+1,\dots, n$  but then got ``uncoupled'' and stayed so until and including the state indexed by~$k$. This means
\begin{equation}
\wt P(\phi_k'\ne\phi_k^\star) \le \wt P\biggl(\,\bigcap_{i=k}^{ n}\{\phi_i'\ne\phi_i^\star\}\biggr)
+\sum_{j=k+1}^n  \wt P\biggl(\,\{\phi_j'=\phi_j^\star\}\,\cap\,\bigcap_{i=k}^{j-1}\{\phi_i'\ne\phi_i^\star\}\biggr).
\end{equation}
Using \eqref{Eq:1.35}, the first probability is bounded inductively by  the product of $1-\texte^{-A_{n,i}}$ for~$i$ ranging from~$k$ to~$n$ while \twoeqref{Eq:1.35}{Eq:1.36} bounds the second probability by $1-\texte^{-\delta_{j-1}}$ times the product of $1-\texte^{-A_{n,i}}$ for~$i$ ranging from~$k$ to~$j-2$.  Since $P(\phi_k'\ne\phi_k^\star)$  dominates  the total variation on  the left of \eqref{Eq:1.28}, the claim follows. 
\end{proofsect}

\subsection{Proof of Theorem~\ref{thm-1}}
We are now in a position to prove the conclusion concerning the covariance of the field. We start with sub and supercritical cases that can be handled concurrently:

\begin{proofsect}{Proof of Theorem~\ref{thm-1}, $\beta\ne\beta_\cc$}
Suppose~$b\ge2$ and~$\beta>0$ are such that either Theorem~\ref{thm-3} or Theorem~\ref{thm-5} applies, whichever is relevant. This in particular means existence of a $1$-periodic continuously differentiable~$v_\star\colon\R\to\R$ such that \twoeqref{Eq:3.32a}{Eq:3.32}  hold  with some~$C',\eta'>0$. (For~$\beta<\beta_\cc$ we put  $v_\star=0$ and invoke \twoeqref{E:3.13}{E:3.14}.) Moreover, in both cases $\{v_k'\}_{k=0}^n$ and $\{v_k''\}_{k=0}^n$ are  bounded uniformly in~$n\ge1$, which enables Proposition~\ref{prop-4.1} and Lemma~\ref{lemma-4.6}. 

In the notation of Lemma~\ref{lemma-1.4} the bound \eqref{Eq:3.32a} implies 
\begin{equation}
\delta_k\le C'\texte^{-\eta'\min\{k,n-k\}}+\frac12\max\bigl\{|\sigma_k^2-1|,|\sigma_k^{-2}-1|\bigr\}
\end{equation}
while the uniform boundedness of~$\{v_k\}_{k=0}^n$ and $\{v_\star\}_{k=0}^n$ gives $\sup_{n\ge k\ge0}A_{n,k}<\infty$. Relying on \eqref{E:3.13} and, for~$\beta>\beta_\cc$, on \eqref{Eq:3.32a} and the assumption that~$\{\frakd_j\}_{j\ge0}$ decays exponentially, $k\mapsto\delta_k$ decays exponentially with $\min\{k,n-k\}$. The coupling inequality \eqref{Eq:1.28} then shows that $\Vert P(\phi_k\text{ mod }1\in\cdot)-\nu_\star\Vert_{\text{TV}}$ decays exponentially in~$\min\{k,n-k\}$. With the help of \eqref{Eq:3.32} this gives
\begin{equation}
\label{E:4.45i}
\sup_{n\ge1}\sum_{k=0}^n\Bigl|E\bigl(v_k'(\phi_{k+1})^2\bigr) - E_\star\bigl(v_\star'(\phi)^2\bigr)\Bigr|<\infty,
\end{equation}
where~$E_\star$ denotes expectation with respect to~$\nu_\star$ and where the $k=0$ term is handled using the uniform bound on~$v_0'$ and~$v_\star'$. 
 Since $E_\star(v_\star'(\phi)^2)$ does not depend on~$k$, it  follows that the sum on the right of \eqref{Eq:1.20a} equals
\begin{equation}
(n-k)E_\star(v_\star'(\phi)^2)+O(1).
\end{equation}
This implies the claim with
\begin{equation}
\label{E:4.39}
\sigma^2(\beta) := \frac1\beta-b\frac1{\beta^2}\frac{b+1}{b-1}\frac{\int_0^1 \texte^{-(b+1)v_\star(\phi)} v_\star'(\phi)^2\,\textd\phi}{\int_0^1\texte^{-(b+1) v_\star(\phi)}\,\textd\phi}
\end{equation}
which equals~$1/\beta$ when~$\beta\le\beta_\cc$ as~$v_\star=0$ but is strictly less than that whenever~$v_\star$ is non-trivial, as is the case for~$\beta>\beta_\cc$. The continuity in~$\beta$ follows from the continuity of~$v_\star$ in~$\beta$, which is itself implied by uniqueness of the fixed point \eqref{Eq:3.33}. The monotonicity for~$\beta>\beta_\cc$ is a consequence of Ginibre inequalities proved by Fr\"ohlich and Park~\cite{FP}.
\end{proofsect}

\begin{remark}
While \eqref{E:4.39} makes it perfectly clear that~$\sigma^2(\beta)<1/\beta$ once~$v_\star$ is non-trivial, the reader may wonder why (since~$\sigma^2(\beta)$ is a variance) the second term never exceeds the first one. This follows from the formula
\begin{equation}
\label{E:4.48}
\begin{aligned}
&\sigma^2(\beta) \int_0^1\texte^{-(b+1) v_\star(\phi)}\textd\phi
\\
&\quad= \int\biggl(\zeta-\frac1\beta\frac b{b-1}\bigl[v_\star'(\phi+\zeta)-v_\star'(\phi)\bigr]\biggr)^2\texte^{-bv_\star(\phi)-bv_\star(\phi+\zeta)}\mu_{1/\beta}(\textd\zeta)1_{[0,1)}(\phi)\textd\phi
\end{aligned}
\end{equation}
which can be checked using the calculations underpinning the proof of Proposition~\ref{prop-4.1} and the fact that~$v_\star$ solves \eqref{Eq:3.33}. 
\end{remark}

\begin{remark}
\label{rem-sigma-asymp}
The computations in Section~\ref{sec-6}  (see, e.g., Lemma~\ref{lemma-6.8})  show that
\begin{equation}
\lambda_\star(1)=\sqrt{\frac{ 2(b^3-1)}{(b-1)^2(b+1)^3}}\sqrt{b\theta-1}+O(b\theta-1).
\end{equation}
This, along with the bounds \eqref{Eq:3.30w} implies
\begin{equation}
v_\star'(\phi)=4\pi\,\sqrt{\frac{ 2 (b^3-1)}{(b-1)^2(b+1)^3}}\sqrt{b\theta-1}\, \sin(2\pi\phi)+O(b\theta-1).
\end{equation}
It follows that the oscillation of~$v_\star$ is at most order~$\sqrt{b\theta-1}$ and so, up to corrections of order~$(b\theta-1)^{3/2}$,   the ratio  in \eqref{E:4.39} reduces to $\int_0^1 v_\star'(\phi)^2\textd\phi$. As $\int_0^1\sin(2\pi\phi)^2\textd\phi=1/2$, we~get
\begin{equation}
\sigma^2(\beta)=\frac1\beta - (4\pi)^2\frac1{\beta_\cc^2} b\frac{b+1}{b-1}\frac{ 2(b^3-1)}{(b-1)^2(b+1)^3}(b\theta-1)\frac12+O(|b\theta-1|^{3/2}).
\end{equation}
A computation shows $ b\theta-1 = \frac{2\pi^2}{\beta_\cc^2}(\beta-\beta_\cc) +O((\beta-\beta_\cc)^2)$ and gives \eqref{E:1.16ui}.
\end{remark}

With the above stated, the critical case requires only minor changes:

\begin{proofsect}{Proof of Theorem~\ref{thm-1}, $\beta=\beta_\cc$}
Note that~$v_\star=0$ in this case and so~$\nu_\star$ is just the Lebesgue measure on~$[0,1)$.  Define $A_n := \sup_{0 \le k \le n} A_{n,k}$.  Lemma~\ref{lemma-1.4} then gives
\begin{equation}
\begin{aligned}
\sum_{k=1}^{n+1}\biggl| &E\bigl(v_{k-1}'(\phi_{k})^2\bigr) - \int_{[0,1)}v_{k-1}'(\phi)^2\textd\phi\biggr|
\le \sum_{k=1}^{n+1} \Vert v_{k-1}'\Vert_\infty^2\bigl\Vert P(\phi_k\text{ mod }1\in\cdot)-\nu_\star\bigr\Vert_{\text{TV}}
\\
& \le \sum_{k=1}^{n+1} \Vert v_{k-1}'\Vert_\infty^2 (1-\texte^{-A_n})^{n-k+1}+\sum_{j=2}^{n+1}\delta_{j-1}\sum_{k=1}^{j-1}\Vert v_{k-1}'\Vert_\infty^2(1-\texte^{-A_n})^{j-1-k}. 
\end{aligned}
\end{equation}
Assumption~\ref{ass-1} and \eqref{E:3.16} imply $\sup_{n\ge1}A_n<\infty$ and $\delta_k=O(k^{-1/2})+\frakd_{\min\{k,n-k\}}$ uniformly in~$n\ge k\ge0$. Since $\Vert v'_{k-1}\Vert_\infty=O(k^{-1/2})$, the second sum on the right is order~$j^{-1}$ and so the sum over~$j$ is finite uniformly in~$n$.

Writing the first term on the right of \eqref{Eq:1.44}  as $[4\pi A+o(1)]\sin(2\pi z)/\sqrt{k}$ we have
\begin{equation}
v_k'(z)^2 = (4\pi A)^2\sin(2\pi z)^2 k^{-1}+O(k^{-3/2})+O\biggl(\,\,k^{-1}\!\!\!\!\!\!\!\sum_{j\ge\min\{\sqrt k,n-k\}}\frakd_j\biggr).
\end{equation}
Using $\int_0^1\sin(2\pi\phi)^2\textd\phi=1/2$, the sum in \eqref{Eq:1.4a} thus equals
\begin{equation}
\label{E:4.53}
\frac1{\beta_\cc} (n-k) + \biggl[\frac1{\beta_\cc^2}b\frac{b+1}{b-1}8\pi^2 A^2\biggr]\log\frac nk
\end{equation}
plus a quantity of order
\begin{equation}
\sum_{k=1}^n k^{-3/2}+\sum_{k=1}^n k^{-1}\sum_{j\ge\min\{\sqrt k,n-k\}}\frakd_j\le O(1)+\sum_{j\ge1}\frakd_j\biggl(1+\sum_{k=1}^{j^2}\frac1k\biggr).
\end{equation}
Under the assumption that~$\sum_{j\ge1}\frakd_j\log(j)<\infty$, this is bounded uniformly in~$n\ge1$. Plugging $n:=\log_{b^{1/2}}(\diam(\Lambda_n))$ and $k:=\log_{b^{1/2}}(2+d(x,y))+O(1)$ in \eqref{E:4.53} and substituting for~$A$ yields the claim. 
\end{proofsect}

\begin{remark}
\label{rem-crit-near-crit1}
The previous proofs demonstrate the reason for the numerical closeness of the coefficients in front of the second order term of the covariance at~$\beta_\cc$ and the near-critical expansion of~$\sigma^2(\beta)$ in~\eqref{E:1.16ui}. Indeed, both rely on the expansion $v_k'(\phi)=(4\pi A \epsilon)\sin(2\pi\phi)+ O(\epsilon^2)$,  where~$\epsilon:=k^{-1/2}$ in the critical case and~$\epsilon:=\sqrt{ 2(b\theta-1)}$ in the near-critical case and~$A$ is the constant on the right of \eqref{Eq:3.18u}. Plugging this  into  \eqref{Eq:1.4a} with  the help of  \eqref{Eq:1.20a}, we also need that $E(\sin(2\pi\phi_k)^2)=1/2+ O(\epsilon)$  once~$k$ and~$n-k$ are sufficiently large.
\end{remark}

\begin{remark}
\label{rem-max}
The covariance computation reveals a log-correlated structure that has, in recent years, allowed  control of the limit law of the maximum and the full extremal process of the underlying field in a number of specific models of interest. Two examples most relevant for the present work where this has been done are the Branching Random Walk (A\"idekon~\cite{Aidekon}, Madaule~\cite{Madaule2}) and the GFF on subdomains of~$\Z^2$ (Bramson, Ding and Zeitou\-ni~\cite{BDingZ}, Biskup and Louidor~\cite{BL1,BL2,BL3}; see also~\cite{B-notes}). In~\cite{BHu} the present authors used the tree-indexed Markov chain representation of the hierarchical DG~model to control the maximum and the extremal process throughout the subcritical regime. This leaves the question of what happens at, and beyond,~$\beta_\cc$.

Unfortunately, the covariance calculation is not sufficient to make reliable predictions about the law of the maximum. Indeed, we need  a  sharp asymptotic of the probability that~$\phi_x$ exceeds quantities of order~$n$, and likely quite a bit more. We presently do not see how to extract this information from our calculations.
\end{remark}

\section{Fractional charge asymptotic}
\label{sec-5}\noindent
We now move to the computation of the asymptotic of the fractional charge  correlation  stated in Theorem~\ref{thm-2}. We start with some general observations that apply  to  all~$\beta>0$ and then prove the statement for subcritical, critical and supercritical~$\beta$.

\subsection{General considerations}
The proof of Theorem~\ref{thm-2} again relies on the Markov chain representation of the field; see Section~\ref{sec-3.1}. As before, we will write $\phi_{n+1}, \phi_n,\dots,\phi_0$ for a run of the chain along the path from the root to a generic leaf-vertex indexed so that~$\phi_n$ is the value at the root  and $\phi_{n+1}=0$.  We start by introducing the key iteration for the fractional-charge setting:

\begin{lemma}
Let~$\alpha\in\R$ and, given $\{w_1(q)\}_{q\in\Z}\in\ell^1(\Z)$, define
$\{w_k(\cdot)\}_{k=1}^{ n+1}$ by the recursion
\begin{equation}
\label{Eq:5.2}
w_{k+1}(q):=\frac{a_{k-1}(0)^b}{a_{k}(0)}\sum_{\begin{subarray}{c}
\ell_1,\dots,\ell_b\in\Z\\\ell_1+\dots+\ell_b=q
\end{subarray}}
\biggl(\,\prod_{i=1}^{b-1}\frac{a_{k-1}(\ell_i)}{a_{k-1}(0)}\biggr)\,w_k(\ell_b)\,\theta_{k}^{(q+\alpha)^2}.
\end{equation}
Then $\{w_k(q)\}_{q\in\Z}\in\ell^1(\Z)$ for all~$k=1,\dots, n+1$  and so we can set 
\begin{equation}
\label{Eq:5.2b}
f_k(z) := a_{k-1}(0)\texte^{v_{k-1}(z)}\sum_{q\in\Z}w_k(q)\texte^{2\pi\texti qz}.
\end{equation}
Moreover, with $\FF_k$ as in \eqref{Eq:4.1},
\begin{equation}
\label{Eq:2.4}
E\bigl(\texte^{2\pi\texti\alpha \phi_k}f_k(\phi_k)\,\big|\,\FF_{k+1}\bigr) =
\texte^{2\pi\texti\alpha \phi_{k+1}}f_{k+1}(\phi_{k+1})
\end{equation}
then holds for all~$k=1,\dots,n$. In short, $\{\texte^{2\pi\texti\alpha \phi_k}f_k(\phi_k)\}_{k=1}^{ n+1}$ is a reverse martingale.
\end{lemma}

\begin{proofsect}{Proof}
Using that $a_k \in\ell^1(\Z)$ for all  $k=0,\dots,n$  we check that~$w_{k+1}\in\ell^1(\Z)$ whenever $w_k\in\ell^1(\Z)$, so the main part to prove is \eqref{Eq:2.4}.
Continuing to abbreviate $\beta_k:=\beta\sigma_k^{-2}$, let  $k=1,\dots,n$  and use \eqref{Eq:5.2b} to write 
\begin{equation}
\label{E:5.4w}
\begin{aligned}
E\bigl(\texte^{2\pi\texti\alpha \phi_k}&f_k(\phi_k)\,\big|\,\FF_{k+1}\bigr) 
\\
&= \texte^{v_{k}(\phi_{k+1})}\int \texte^{-b v_{k-1}(\phi_{k+1}+\zeta)}\texte^{2\pi\texti\alpha (\phi_{k+1}+\zeta)}  f_k(\phi_{k+1}+\zeta)\mu_{1/\beta_{k}}(\textd\zeta)
\\
&= a_{k-1}(0)\texte^{2\pi\texti \alpha\phi_{k+1}}\texte^{v_{k}(\phi_{k+1})}
\\
&\qquad\qquad\times\int \texte^{-(b-1) v_{k-1}(\phi_{k+1}+\zeta)}\texte^{2\pi\texti\alpha \zeta}  \sum_{\ell_b\in\Z}  w_k(\ell_b)\,\texte^{2\pi\texti \ell_b(\phi_{k+1}+\zeta)}\mu_{1/\beta_k}(\textd\zeta).
\end{aligned}
\end{equation}
The sum over~$\ell_b$ can be exchanged with the integral thanks to $\{w_k(q)\}_{q\in\Z}\in\ell^1(\Z)$. Invoking $\texte^{-v_{k-1}(z)}=\sum_{\ell\in\Z}a_{k-1}(\ell)\texte^{2\pi\texti\ell z}$ along with the fact that $\int\texte^{2\pi\texti r\zeta}\mu_{1/\beta_k}(\textd\zeta)=\theta_k^{r^2}$ for any~$r\in\R$, we then compute the last integral in \eqref{E:5.4w} to be
\begin{equation}
\begin{aligned}
\sum_{\ell_b\in\Z}w_k(\ell_b)\int &\texte^{-(b-1) v_{k-1}(\phi_{k+1}+\zeta)}\texte^{2\pi\texti\alpha \zeta}\texte^{2\pi\texti \ell_b(\phi_{k+1}+\zeta)}
\mu_{1/\beta_k}(\textd\zeta)
\\
&=\sum_{\ell_1,\dots,\ell_{b}\in\Z}
\biggl(\,\prod_{i=1}^{b-1}a_{k-1}(\ell_i)\biggr)\,w_k(\ell_b)\texte^{2\pi\texti(\ell_1+\dots+\ell_b)\phi_{k+1}}\,\theta_k^{(\ell_1+\dots+\ell_b+\alpha)^2}.
\end{aligned}
\end{equation}
Writing the sum as two sums, one over~$q\in\Z$ and the other over~$\ell_1,\dots,\ell_b\in\Z$ subject to~$\ell_1+\dots+\ell_b=q$,  the definition of~$w_{k+1}$ reduces the expectation of interest to
\begin{equation}
E\bigl(\texte^{2\pi\texti\alpha \phi_k}f_k(\phi_k)\,\big|\,\FF_{k+1}\bigr) 
=\texte^{2\pi\texti \alpha\phi_{k+1}}a_k(0)\texte^{v_{k}(\phi_{k+1})}\sum_{q\in\Z} w_{k+1}(q)\texte^{2\pi\texti q\phi_{k+1}}.
\end{equation}
This equals $\texte^{2\pi\texti\alpha \phi_{k+1}}f_{k+1}(\phi_{k+1})$ thus proving \eqref{Eq:2.4} as desired.
\end{proofsect}

Note that the case~$k=0$ is excluded from the previous lemma in light of the conditional law of~$\phi_0$ given~$\FF_1$ being of a different form than the law of $\phi_k$ given~$\FF_{k+1}$ for~$k\ge1$. To get the iteration started, we thus need:

\begin{lemma}
\label{lemma-5.2}
Suppose that~$f\colon\R\to\R$ admits the Fourier representation
\begin{equation}
f(z) = \sum_{q\in\Z}w(q)\texte^{2\pi\texti q z}
\end{equation}
such that $\{w(q)\}_{q\in\Z}\in\ell^1(\Z)$. Then for all $\alpha\in\R$,
\begin{equation}
\label{E:5.8}
E\bigl(\,\texte^{2\pi\texti\alpha\phi_0}f(\phi_0)\,\big|\,\FF_1\bigr) = \texte^{v_0(\phi_1)+2\pi\texti\alpha\phi_1}
\sqrt{\frac{2\pi\sigma_0^2}\beta}
\sum_{q\in\Z} w\ast a(q)\,\theta_0^{(q+\alpha)^2}\texte^{2\pi\texti q\phi_1},
\end{equation}
where~$\{a(q)\}_{q\in\Z}$ are the Fourier coefficients of~$\nu$; i.e., $a(q):=\int_{[0,1)}\texte^{- 2\pi\texti qz}\nu(\textd z)$, and $w\ast a$ is the usual convolution,
\begin{equation}
w\ast a(q):=\sum_{\ell\in\Z} w(q-\ell)a(\ell),\quad q\in\Z.
\end{equation}
\end{lemma}

\begin{proofsect}{Proof}
In light of $\{w( q)\}_{ q\in\Z}\in\ell^1(\Z)$  it suffices to focus on  $f(z): = \texte^{2\pi\texti qz}$ for some~$q\in\Z$. With $\beta_0:=\beta\sigma_0^{-2}$, the bottom line in \eqref{E:1.23} then gives
\begin{equation}
E\bigl(\,\texte^{2\pi\texti(q+\alpha)\phi_0}\,\big|\,\FF_1\bigr)
=\texte^{v_0(\phi_1)+2\pi\texti(q+\alpha)\phi_1}\int \texte^{2\pi\texti(q+\alpha)(\phi-\phi_1) - \frac{\beta_0}2(\phi-\phi_1)^2}\nu(\textd\phi).
\end{equation}
Denote the last integral as~$h_q(\phi_1)$. 
The $1$-periodicity of~$\nu$ implies that  $h_q$ is a $1$-periodic function. The Gaussian decay of the integrand in turn permits us to swap any number of derivatives with respect to~$\phi_1$ with the integral which means that $h_q\in C^\infty(\R)$. It follows that one can express $h_q$ as a uniformly convergent Fourier series. Writing~$\nu$ as the sum of shifts of its restriction to~$[0,1)$, for the~$m$-th Fourier coefficient of~$h_q$ we get
\begin{equation}
\begin{aligned}
\hat h_q(m):=&\int_{[0,1)}  h_q(\phi_1) \, \texte^{-2\pi\texti m\phi_1}\textd\phi_1
\\
&=\sum_{n\in\Z}\int_{[0,1)\times[0,1)}\texte^{-2\pi\texti m\phi_1}\,
\texte^{2\pi\texti(q+\alpha)(\phi+n-\phi_1) - \frac{\beta_0}2(\phi+n-\phi_1)^2}\,\textd\phi_1\nu(\textd\phi)
\\
&=\sum_{n\in\Z}\int_{[0,1)\times[0,1)}\texte^{-2\pi\texti m \phi}\,
\texte^{2\pi\texti(q+\alpha+m)(\phi+n-\phi_1) - \frac{\beta_0}2(\phi+n-\phi_1)^2}\,\textd\phi_1\nu(\textd\phi),
\end{aligned}
\end{equation}
where we used that~$mn\in\Z$. Now combine the sum over~$n\in\Z$ with the integral over~$\phi_1\in[0,1)$ to get an integral over~$\phi_1\in\R$. The shift-invariance of the Lebesgue measure allow us to replace~$\phi_1-\phi$ by~$\phi_1$ which brings the integral into a product form. A computation shows $\hat h_q(m)=a(m)\theta_0^{(q+m+\alpha)^2}\sqrt{2\pi/\beta_0}$. 

From the previous calculations we conclude
\begin{equation}
E\bigl(\,\texte^{2\pi\texti(q+\alpha)\phi_0}\,\big|\,\FF_1\bigr)
=\texte^{v_0(\phi_1)+2\pi\texti(q+\alpha)\phi_1}\sqrt{\frac{2\pi}{\beta_0}}\sum_{m\in\Z} a(m)\theta^{(q+m+\alpha)^2}\texte^{2\pi\texti m\phi_1}
\end{equation}
which is \eqref{E:5.8} when $f(z) = \texte^{2\pi\texti qz}$.
\end{proofsect}

To connect the above to the main objective of Theorem~\ref{thm-2}, we note:

\begin{corollary}
\label{cor-5.3}
Let~$\beta>0$ and~$\alpha\in\R$ and let $\{w_k\}_{k=1}^{ n+1}$ be defined by \eqref{Eq:5.2} with
\begin{equation}
\label{Eq:2.8}
w_1(q):=
\sqrt{\frac{2\pi\sigma_0^2}\beta}
\,\frac{a(q)}{a_0(0)}\,\theta_0^{(q+\alpha)^2},\quad q\in\Z,
\end{equation}
where $\{a(q)\}_{q\in\Z}$  denote the Fourier coefficients of~$\nu$.
Let~$\{f_k\}_{k=1}^{ n+1}$ be defined from $\{w_k\}_{k=1}^{ n+1}$ as in \eqref{Eq:5.2b}. Then $w_k(\cdot)$ is strictly positive for all~$k=1,\dots,n+1$  and
\begin{equation}
\label{Eq:2.9}
\langle\texte^{2\pi\texti\alpha(\phi_x-\phi_y)}\rangle_{n,\beta} = E\bigl(\,|f_k(\phi_k)|^2\bigr)
\end{equation}
holds with~$k:=k(x,y)$ for all distinct $x,y\in\Lambda_n$.  (Here~$k(x,y)$ is as in \eqref{Eq:4.3i}.)
\end{corollary}

\begin{proofsect}{Proof}
The strict positivity of~$w_k(\cdot)$ is immediate from \eqref{Eq:2.8}, \eqref{Eq:5.2}  and the positivity of~$a_k$'s so we only need to prove \eqref{Eq:2.9}.
Let~$x,y\in\Lambda_n$ be distinct and set~$k:=k(x,y)$. The tree-indexed Markov chain representation of the field~$\{\phi_x\}_{x\in\Lambda_n}$ from Section~\ref{sec-3.1} yields
\begin{equation}
\label{Eq:5.13i}
\langle\texte^{2\pi\texti\alpha(\phi_x-\phi_y)}\rangle_{n,\beta}=E\Bigl(E(\texte^{2\pi\texti\alpha\phi_0}\,|\,\FF_k)E(\texte^{-2\pi\texti\alpha\phi_0}\,|\,\FF_k)\Bigr).
\end{equation}
The first conditional expectation corresponds to the choice $f:=1$ in \eqref{E:5.8} for which $w(q):=\delta_{q,0}$. The identity  \eqref{E:5.8} then gives~$E(\texte^{2\pi\texti\alpha\phi_0}\,|\,\FF_1)=\texte^{2\pi\texti\alpha\phi_1}f_1(\phi_1)$ for 
\begin{equation}
f_1(\phi_1):=\texte^{v_0(\phi_1)}\sqrt{\frac{2\pi\sigma_0^2}\beta}\sum_{q\in\Z}a(q)\,\theta^{(q+\alpha)^2}\,\texte^{2\pi\texti q\phi_1}.
\end{equation}
Since this coincides with the function~$f_1$ from \eqref{Eq:5.2b} for~$w_1$ as in \eqref{Eq:2.8}, using \eqref{Eq:2.4} we obtain $E(\texte^{2\pi\texti\alpha\phi_0}\,|\,\FF_k)=\texte^{2\pi\texti\alpha\phi_k}f_k(\phi_k)$. Noting that the second expectation in \eqref{Eq:5.13i} is  the  complex conjugate of the first, plugging this in \eqref{Eq:5.13i} yields \eqref{Eq:2.9}.
\end{proofsect}

We remark that, besides~$k$, the sequence $\{w_k(q)\}_{q\in\Z}$ depends also on~$n$ through the $n$-dependence of the sequence~$\{\sigma_k^2\}_{k=0}^n$. While we keep that dependence implicit, it may need to be noted in statements where uniformity in~$n$ is required.

\subsection{Below criticality}
Corollary~\ref{cor-5.3} tells us that, for the asymptotic of the fractional charge at large separations of~$x$ and~$y$, we need to track the large-$k$ asymptotic form of~$f_k$ initiated from \eqref{Eq:2.8}. As it turns out, the cases of~$\beta\le\beta_\cc$ are united by the fact that~$f_k$ is dominated by the coefficient~$w_k(0)$ for $k$ large. 

\begin{lemma}
\label{cor-5.5}
Let~$\beta\in(0,\beta_\cc]$, $\alpha\in(-\ffrac12,\ffrac12)$ and let $w_1$ be as in \eqref{Eq:2.8}. Then
\begin{equation}
\label{Eq:5.14w}
 \lim_{m\to\infty} \sup_{n\ge k\ge m} \,\sum_{q\ne0}\frac{w_k(q)}{w_k(0)}=0.
\end{equation}
\end{lemma}

\begin{proofsect}{Proof}
We will prove that, under the stated conditions on~$\beta$ and~$\alpha$, there exist~$\wt C>0$,~$\wt\eta>0$ and $k_0\ge2$ such that
\begin{equation}
\label{Eq:5.10}
\sum_{q\ne0}\frac{w_k(q)}{w_k(0)}\le\wt C\sum_{j=2}^{k-k_0} \texte^{-\wt\eta j}\sum_{q\ne0}\frac{a_{k-j}(q)}{a_{k-j}(0)}
\end{equation}
holds for all~$k\ge k_0+2$, uniformly in~$n$. For this we note that, reducing the sum in \eqref{Eq:5.2} to the term with $\ell_1=\dots=\ell_{b-1}=0$ shows
\begin{equation}
\label{E:5.16w}
w_{k+1}(0)\ge \frac{a_{k-1}(0)^b}{a_{k}(0)} \,\theta_{k}^{\alpha^2}\,w_k(0).
\end{equation}
Dividing each side of \eqref{Eq:5.2} by the corresponding side of this bound then gives
\begin{equation}
\label{Eq:5.12}
\frac{w_{k+1}(q)}{w_{k+1}(0)}\le\sum_{\begin{subarray}{c}
\ell_1,\dots,\ell_b\in\Z\\\ell_1+\dots+\ell_b=q
\end{subarray}}
\biggl(\,\prod_{i=1}^{b-1}\frac{a_{k-1}(\ell_i)}{a_{k-1}(0)}\biggr)\,\frac{w_k(\ell_b)}{w_k(0)}\,\theta_{k}^{(q+\alpha)^2-\alpha^2}.
\end{equation}
Now observe that for~$q\ne0$ we have $(q+\alpha)^2-\alpha^2= q^2+2q\alpha\ge 1-2|\alpha|$ and abbreviate
\begin{equation}
u_k:=\sum_{q\ne0}\frac{a_{k-1}(q)}{a_{k-1}(0)}.
\end{equation}
Summing \eqref{Eq:5.12} over~$q\ne0$, we have two cases  to consider  on the right-hand side: either~$\ell_b=0$, in which case at least one of~$\ell_1,\dots,\ell_{b-1}$ must be non-zero, or $\ell_b\ne0$ in which case the remaining indices can, as an upper bound, be summed over as free variables. This yields
\begin{equation}
\label{Eq:5.14}
\sum_{q\ne0}\frac{w_{k+1}(q)}{w_{k+1}(0)}
\le (1+u_k)^{b-2}b\theta_{k}^{1-2|\alpha|}u_k+ (1+u_k)^{b-1}\theta_{k}^{1-2|\alpha|}\sum_{q\ne0}\frac{w_k(q)}{w_k(0)},
\end{equation}
where the factor~$b$ in the first term dominates the number of choices of the first non-zero term among $\ell_1,\dots,\ell_{b-1}$. 

Now observe that Theorems~\ref{thm-3} and~\ref{thm-4} imply $\sup_{n\ge k \ge m}u_k \to0$ as~$m\to\infty$. In light of~$|\alpha|<1/2$, Assumption~\ref{ass-1} and~$\sup_{k\le n}\theta_k<1$, there exists~$k_0\ge1$ such that
\begin{equation}
\texte^{-\wt\eta}:=\sup_{n\ge k\ge k_0}(1+u_k)^{b-1}\theta_{k}^{1-2|\alpha|}<1.
\end{equation}
Using this in \eqref{Eq:5.14} and iterating we get \eqref{Eq:5.10} by invoking~$\{w_{k_0}(q)\}_{q\in\Z}\in\ell^1(\Z)$ which is checked from \eqref{Eq:5.2}, \eqref{Eq:2.8} and the boundedness of $\{a(q)\}_{q\in\Z}$. With \eqref{Eq:5.10} in hand, we then get \eqref{Eq:5.14w} by using again $\sup_{n\ge k\ge  m}u_k\to0$ as~$m\to\infty$.
\end{proofsect}

Using very similar arguments, we also get:

\begin{lemma}
\label{lemma-5.5b}
Let~$\beta\in(0,\beta_\cc]$, $\alpha\in(-\ffrac12,\ffrac12)$ and let $w_1$ be as in \eqref{Eq:2.8}. Then
\begin{equation}
\label{E:5.21i}
0<\inf_{n\ge k\ge1}\frac{w_{k+1}(0)}{w_k(0)\theta_k^{\alpha^2}}\le \sup_{n\ge k\ge1}\frac{w_{k+1}(0)}{w_k(0)\theta_k^{\alpha^2}}<\infty.
\end{equation}
\end{lemma}

\begin{proofsect}{Proof}
For the lower bound in \eqref{E:5.21i} we use \eqref{E:5.16w} to get
\begin{equation}
\frac{w_{k+1}(0)}{w_k(0)\theta_k^{\alpha^2}}\ge\frac{a_{k-1}(0)^b}{a_k(0)}.
\end{equation}
The quantity on the right is positive uniformly in $n\ge k\ge1$ thanks to \eqref{E:3.12} (for~$\beta<\beta_\cc$) and \eqref{E:3.15} (for $\beta=\beta_\cc$) and the bound~$a_k(0)\le (\sum_{q\in\Z} a_{k-1}(q))^b$. 

For the upper bound in \eqref{E:5.21i} we drop the condition $\ell_1+\dots+\ell_b=0$ in the definition of $w_{k+1}(0)$ in \eqref{Eq:5.2} to get
\begin{equation}
\label{E:5.23i}  
\frac{w_{k+1}(0)}{w_k(0)\theta_k^{\alpha^2}}
\le \frac{a_{k-1}(0)^b}{a_k(0)} \, \biggl(\,\sum_{\ell\in\Z}\frac{a_{k-1}(\ell)}{a_{k-1}(0)}\biggr)^{b-1}\sum_{q\in\Z}\frac{w_k(q)}{w_k(0)}.
\end{equation}
 Now observe that $a_{k-1}(0)^b \le a_k(0)$, and the first sum is bounded uniformly in $ n + 1\ge k\ge1$ thanks to \eqref{E:3.12} and \eqref{E:3.15}, while the second sum is bounded thanks to \eqref{Eq:5.14w}.
\end{proofsect}

We are now ready for:

\begin{proofsect}{Proof of Theorem~\ref{thm-2}, $\beta<\beta_\cc$}
 We will derive an approximate recursive equation for the sequence $\{w_k(0)\}_{k=1}^{n+1}$.  For this we separate the term with $\ell_1=\dots=\ell_b=0$ in \eqref{Eq:5.2} and use arguments similar to those underlying the proof of \eqref{E:5.23i} to bound the remaining terms. This yields
\begin{equation}
\begin{aligned}
\biggl|\,w_{k+1}(0)-&\frac{a_{k-1}(0)^b}{a_{k}(0)}\theta_{k}^{\alpha^2}w_k(0)\biggr|
\\
&\le \Biggl[b\biggl(\,\sum_{\begin{subarray}{c}
\ell\in\Z\\\ell\ne0
\end{subarray}}\frac{a_{k-1}(\ell)}{a_{k-1}(0)}\biggr)\biggl(\sum_{\ell\in\Z}\frac{a_{k-1}(\ell)}{a_{k-1}(0)}\biggr)^{b-2}\sum_{\ell\in\Z}\frac{w_k(\ell)}{w_k(0)}\Biggr] \theta_{k}^{\alpha^2} w_k(0),
\end{aligned}
\end{equation}
where the prefactor~$b$ accounts for the choice of the smallest index with~$\ell_i\ne0$. The  exponential decay \eqref{E:3.12} in Theorem~\ref{thm-3} shows that the first sum on the right decays exponentially in~$k$ and that the second sum is bounded uniformly in~$n\ge k\ge1$. Lemma~\ref{cor-5.5} gives the same for the last sum. On the left-hand side, Theorem~\ref{thm-3} gives that $a_{k-1}(0)^b/a_k(0)$ differs from~$1$ by a factor that decays exponentially in $k$. It follows that for each $\beta\in(0,\beta_\cc)$ there exists~$\eta>0$ such that
\begin{equation}
\label{E:5.21}
w_{k+1}(0) = \texte^{O(\texte^{-\eta k})} \theta_{k}^{\alpha^2}w_k(0),
\end{equation}
with the implicit constant in $O(\texte^{-\eta k})$ bounded uniformly in~$n\ge k\ge1$. Here we relied on Lemma~\ref{lemma-5.5b} in moving the error to  the  exponent for small~$k$.

We now set  $n':=\lfloor n/2\rfloor$ and let  $C_n:= (w_{ n'}(0)\theta^{-\alpha^2 (n'-1)})^2$.   In light of $\theta_i^{\alpha^2} = \theta^{\alpha^2}\theta^{\alpha^2(\sigma_i^2-1)}$,  for~$k\le n'$ we then have 
\begin{equation}
\label{E:5.26ui}
w_k(0)\theta^{-\alpha^2 k} = \sqrt{C_n}\,\exp\Biggl\{-\sum_{i=k}^{ n'-1}\log\Bigl(\frac{w_{i+1}(0)}{w_i(0)\theta_i^{\alpha^2}}\Bigr)+\,\alpha^2\log(1/\theta)\sum_{i=k}^{n'-1}(\sigma_i^2-1)\Biggr\}.
\end{equation}
 The same identity with the limits of the sums interchanged applies to~$k\ge n'$ (with $k\le n$).  Invoking \eqref{E:5.21} and Assumption~\ref{ass-1}, we get $w_k(0)\theta^{-\alpha^2 k}=\sqrt{C_n}\,[1+o(1)]$, where~$o(1)\to0$ as $\min\{k,n-k\}\to\infty$. To see that~$\{C_n\}_{n\ge1}$ is bounded and uniformly positive, take  $k:=1$ in \eqref{E:5.26ui} and  note  that, by Assumption~\ref{ass-1},~$w_1(0)$ is positive and  bounded  uniformly in~$n\ge1$. 
With the help of Lemma~\ref{cor-5.5} and some elementary facts from Fourier analysis, we now conclude that
\begin{equation}
 \bigl|f_k(z)\bigr|^2=\bigl[C_n+o(1)\bigr]\theta^{2\alpha^2k}
\end{equation}
holds with~$o(1)\to0$ as~ $\min\{k,n-k\}\to\infty$, uniformly in~$z\in\R$. Plugging this in \eqref{Eq:2.9} then proves \eqref{E:1.17} for~$\beta<\beta_\cc$.
\end{proofsect}

\renewcommand{\gg}{\gamma}

\subsection{At criticality}
\label{sec-5.2}\noindent
 In the critical and supercritical situations, a simple (albeit approximate) recursion linking~$w_k(0)$ to~$w_{k-1}(0)$ is not enough to capture the actual decay. Indeed, as we will show,~$w_k(0)$ will receive non-trivial contributions from $w_{k-j}(0)$ with~$j\ge2$ as well. To prepare the needed formulas,  we first condense \eqref{Eq:5.2} as
\begin{equation}
\label{Eq:5.19}
w_{k+1}(q) = \theta_{k}^{(q+\alpha)^2} \sum_{\ell\in\Z} \gg_k(q-\ell)w_k(\ell),
\end{equation}
where
\begin{equation}
\label{E:5.29ui}
\gg_k(q):=\frac{a_{k-1}(0)^b}{a_{k}(0)}\sum_{\begin{subarray}{c}
\ell_1,\dots,\ell_{b-1}\in\Z\\\ell_1+\dots+\ell_{b-1}=q
\end{subarray}}
\prod_{i=1}^{b-1}\frac{a_{k-1}(\ell_i)}{a_{k-1}(0)}.
\end{equation}
 Note that, as many of the above objects,  $\gg_k$ depends also on~$n$ but we keep that dependence implicit.
We now  rewrite \eqref{Eq:5.19} as follows:

\begin{lemma}
\label{lemma-5.6a}
For each  $ n+1\ge  k\ge2$ and~$p,q\in\Z$ let
\begin{equation}
\Gamma_{k,1}(p,q):=\theta_{k-1}^{(p+\alpha)^2}\theta^{-\alpha^2}\gg_{k-1}( p-q)
\end{equation}
and for  $ n+1\ge k\ge3$ and~$j=2,\dots,k-1$, set
\begin{equation}
\label{Eq:5.21}
\Gamma_{k,j}(p,q):=\sum_{\begin{subarray}{c}
q_1,\dots,q_{j-1}\in\Z\smallsetminus\{0\}\\ q_0=p,\,q_j= q
\end{subarray}}
\biggl(\,\prod_{i=0}^{j-1}\theta_{k-i-1}^{(q_i+\alpha)^2}\theta^{-\alpha^2}\biggr)\prod_{i=1}^{j}\gg_{k-i}( q_{i-1}-q_i).
\end{equation}
Then for all  $ n+1\ge k\ge2$ and all $p\in\Z$, 
\begin{equation}
\label{Eq:5.23}
w_{k}(p) = \sum_{j=1}^{k-1}\Gamma_{k,j}(p,0)\theta^{\alpha^2j}w_{k-j}(0)+\sum_{q\in\Z\smallsetminus\{0\}}\Gamma_{k,k-1}(p,q)\theta^{\alpha^2(k-1)}w_1(q).
\end{equation}
Moreover, for all~$ n+1\ge k\ge2$,  all  $1\le j< k$ and  all  $p\in\Z$,
\begin{equation}
\label{Eq:5.24}
\sum_{q\in\Z}\Gamma_{k,j}(p,q)\le \theta^{-\tilde c+|p|(|p|-1)\sigma^2_{\text{\rm min}}+(1-2|\alpha|)(j-1)}\prod_{i=1}^{j}\sum_{\ell\in\Z}\gg_{k-i}(\ell)
\end{equation}
holds with~$\tilde c:=\max\{\alpha^2,(1-|\alpha|)^2\}\sum_{i=0}^n|\sigma_i^2-1|$ and $\sigma^2_{\text{\rm min}}:=\inf_{i=0,\dots,n}\sigma_i^2$.
\end{lemma}

\begin{proofsect}{Proof}
We start by proving \eqref{Eq:5.23}. As is checked directly from \eqref{Eq:5.19}, this holds for~$k=2$. For general~$k\ge3$, we plug \eqref{Eq:5.23} for~$w_{k-1}(\ell)$ in 
\begin{equation}
w_k(p) = \theta_{k-1}^{(p+\alpha)^2}\gg_{k-1}(p)w_{k-1}(0)+\sum_{\ell\ne0}\theta_{k-1}^{(p+\alpha)^2}\gg_{k-1}(p-\ell)w_{k-1}(\ell)
\end{equation}
which itself follows from \eqref{Eq:5.19}.  With the help of 
\begin{equation}
\theta_{k-1}^{(p+\alpha)^2}\gg_{k-1}( p)=\Gamma_{k,1}(p,0)\theta^{\alpha^2}
\end{equation}
and
\begin{equation}
\sum_{\ell\ne0}\theta_{k-1}^{(p+\alpha)^2}\gg_{k-1}(p-\ell)\Gamma_{k-1,j}(\ell,q) = \theta^{\alpha^2}\Gamma_{k,j+1}(p,q)
\end{equation}
 we thus get  \eqref{Eq:5.23} for~$k$ from \eqref{Eq:5.23} for~$k-1$. By induction, this proves \eqref{Eq:5.23} to be valid for all~$ n+1\ge k\ge2$. 

For the bound \eqref{Eq:5.24} we first note that
\begin{equation}
\label{E:5.35ui}
(q+\alpha)^2\ge |q|^2-2|\alpha||q|+\alpha^2 = |q|(|q|-1)+|q|(1-2|\alpha|)+\alpha^2.
\end{equation}
For~$q_0:=p$ and any $q_1,\dots,q_{j-1}\ne0$, this implies
\begin{equation}
\label{E:5.37w}
\begin{aligned}
\prod_{i=0}^{j-1}&\theta_{k-i-1}^{(q_i+\alpha)^2}\theta^{-\alpha^2}
\\
&\le \biggl(\,\theta^{\alpha^2(\sigma_{k-1}^2-1)}\prod_{i=1}^{j-1}\theta^{(\sigma_{k-i-1}^2-1)(1-2|\alpha|+\alpha^2)}\biggr)\theta^{|p|(|p|-1)\sigma^2_{k-1}+(1-2|\alpha|)(j-1)}.
\end{aligned}
\end{equation}
 Plugging this in \eqref{Eq:5.21}, using $\sigma_{k-1}^2\ge\sigma^2_{\text{min}}$ and noting that the quantity in the large parentheses is bounded by $\theta^{-\tilde c}$, we get \eqref{Eq:5.24} by performing the sums over $q_1,\dots,q_{j}$. 
\end{proofsect}

The main point of the bound \eqref{Eq:5.24} is that, for $\alpha\in(-\ffrac12,\ffrac12)$  and with Assumption~\ref{ass-2} in force,  the second sum on the right of \eqref{Eq:5.23} decays exponentially faster than the first and can thus be regarded as an error  whenever
\begin{equation}
\theta^{(1-2|\alpha|)}\sum_{\ell\in\Z}\gg_k(\ell)<1.
\end{equation}
This is true for $\min\{k,n-k\}$ large  up to, and even slightly above~$\beta_\cc$.

Formula \eqref{Eq:5.23} reduces the iterations to those of the sequence $\{w_k(0)\}_{k=1}^{n+1}$. (That this suffices at criticality has been shown in Lemma~\ref{cor-5.5}. We wrote \eqref{Eq:5.23} in more generality to prepare for supercritical situations.) Next we need good control of the asymptotic behavior of the coefficients~$\Gamma_{k,j}(0,0)$. This comes in:

\begin{lemma}
\label{lemma-5.7e}
Let $\beta=\beta_\cc$, $\alpha\in(-\ffrac12,\ffrac12)$ and let $A:=[\frac{b^3-1}{(b-1)^2(b+1)^3}]^{1/2}$ be the constant on the right of \eqref{Eq:3.18u}.  For $0\le k\le n$ set
\begin{equation}
\frake_k:=\sum_{j\ge\min\{\sqrt k, n-\sqrt k\}}\frakd_j,
\end{equation}
where~$\{\frakd_j\}_{j\ge0}$ is the sequence from Assumption~\ref{ass-1}. Then 
\begin{equation}
\label{E:5.39w}
\Gamma_{k,1}(0,0)=\exp\Bigl\{- 2(b-1)A^2 k^{-1}+O(k^{-3/2})+O(k^{-1}\frake_{k-1})+O(\frakd_{\min\{k-1,n-k+1\}})\Bigr\}
\end{equation}
and, for $j=2,\dots,k-1$ and~$\eta$ defined by~$\texte^{-\eta}:=b^{-\frac12(1-2|\alpha|)}$,
\begin{equation}
\label{E:5.40w}
\begin{aligned}
\Gamma_{k,j}(0,0) = \bigl(b^{-(1+2\alpha)(j-1)}+&b^{-(1-2\alpha)(j-1)}\bigr)(b-1)^2 A^2 k^{-1}
\\
&+O(\texte^{-\eta j}k^{-3/2})+O\bigl(\texte^{-\eta j}k^{-1}(\frake_{k-1}+\frake_{k-j})\bigr),
\end{aligned}
\end{equation}
where the error bounds are uniform in $j$, $k$ and~$n$ such that $2\le j< k\le n+1$ and in~$\alpha$ on compact subsets of $(-\ffrac12,\ffrac12)$. 
\end{lemma}

\begin{proofsect}{Proof}
The proof will use the conclusions of Theorem~\ref{thm-4} to extract the asymptotic behavior of the relevant coefficients $\gg_k$. For $\Gamma_{k,1}(0,0)$ we only need $\gg_{k-1}(0)$.  Here we treat explicitly the term corresponding to permutations of $(1,-1,0,\dots,0)$ in the sum in \eqref{E:5.29ui} as well as in the sum representing~$a_{k}(0)$ via $a_{k-1}(\cdot)$ and note that, by \eqref{E:3.15}, the remaining terms give contributions of order~$k^{-2}$. Abbreviating $\lambda_k:=a_k(1)/a_k(0)$, recall that \eqref{Eq:3.18u} in Theorem~\ref{thm-4} gives $\lambda_k=A k^{-1/2}+O(k^{-1})+O(k^{-1/2}\frake_k)$. This yields
\begin{equation}
\label{E:5.40q}
\begin{aligned}
\gg_k(0)&=\frac{1+(b-1)(b-2) \lambda_{k-1}^2+O(\lambda_{k-1}^4)}
{1+b(b-1) \lambda_{k-1}^2+O(\lambda_{k-1}^4)}
\\
&=1-2(b-1)\lambda_{k-1}^2+O(\lambda_{k-1}^4) 
\\
&=1-2(b-1)A^2 k^{-1} + O(k^{-3/2})+O(k^{-1}\frake_{k-1})
\\
&=\exp\bigl\{- 2(b-1)A^2 k^{-1}+O(k^{-3/2}) +O(k^{-1}\frake_{k-1})\bigr\},
\end{aligned}
\end{equation}
where we also used that~$\frake_k^2=O(\frake_k)$ and noted that the exponential form is legit thanks to $\gg_{k-1}(0)$ being positive uniformly in $n\ge k\ge1$.  From \eqref{E:5.40q} and $\theta_{k-1} = \theta^{\sigma_{k-1}^2}$ we then readily get 
\begin{equation}
\label{Eq:5.28}
\Gamma_{k,1}(0,0) = \texte^{O(\sigma_{ k-1}^2-1)
}\exp\bigl\{- 2(b-1)A^2 k^{-1}+O(k^{-3/2}) +O(k^{-1}\frake_{k-1})\bigr\}.
\end{equation}
With the help of \eqref{E:1.4u} this shows \eqref{E:5.39w}.

For $\Gamma_{k,j}(0,0)$ with $j=2,\dots,k-1$ we will need the leading-order asymptotic of~$\gg_k(1)$ and suitable bounds on~$\gamma_k(q)$ with~$|q|>1$. For~$\gamma_k(1)$ we treat explicitly the terms corresponding to permutations of $(1,0,\dots,0)$ in \eqref{E:5.29ui} and apply \eqref{Eq:3.18u} to get 
\begin{equation}
\label{E:5.41q}
\begin{aligned}
\gg_k(1) &= (b-1)A\frac 1{\sqrt k}+O(k^{-1})  + O(k^{-1/2}\frake_k)
\\
&=(b-1)A\frac1{\sqrt k}\,\texte^{O(k^{-1/2})+O(\frake_k)}.
\end{aligned}
\end{equation}
For the remaining coefficients we invoke \eqref{E:3.15} with the result
\begin{equation}
\label{E:5.43ui}
\gg_k(q)\le c\Bigl(\frac1{1+\sqrt k}\Bigr)^{|q|},
\end{equation}
where~$c$ is a constant independent of~$q\in\Z$ and~$n\ge k\ge1$.

With \twoeqref{E:5.41q}{E:5.43ui} in hand, we now treat the  $(j-1)$-tuples with $q_1=\dots=q_{j-1}=\pm 1$ in \eqref{Eq:5.21} (with $p=q=0$) and, for the remaining terms, bound the prefactors using \eqref{E:5.37w} while noting that the $q_0,\dots,q_j$ under the sum necessarily satisfy $|q_1-q_0|+\dots+|q_j-q_{j-1}|\ge4$. Denoting the first product in \eqref{Eq:5.21} for~$\theta:=b^{-1}$ and $q_1,\dots,q_{j-1}=1$ by
\begin{equation}
\label{E:5.47q}
h_{k,j}(\alpha):=b^{-(1+2\alpha)(j-1)-\alpha^2(\sigma_{k-1}^2-1)}\prod_{i=1}^{j-1}b^{-(1+\alpha)^2(\sigma_{k-i-1}^2-1)},
\end{equation}
this yields
\begin{equation}
\label{Eq:5.35u}
\begin{aligned}
\Gamma_{k,j}(0,0) = 
\bigl[h_{k,j}(\alpha)+&h_{k,j}(-\alpha)\bigr]
\gamma_{k-1}(1)\gamma_{k-j}(1)\prod_{i=2}^{j-1}\gamma_{k-i}(0)
\\
&+O(1)b^{-(1-2|\alpha|)j}\sum_{\begin{subarray}{c}
p_1,\dots,p_{j-1}\in\Z\\
|p_1|+\dots+|p_{j-1}|\ge4
\end{subarray}}
\prod_{i=1}^{j}\gg_{k-i}(p_i),
\end{aligned}
\end{equation}
where $p_i$ represents~$q_i-q_{i-1}$ and where the symmetry $\gamma_k(-\ell)=\gamma_k(\ell)$ was invoked to simplify the first term. 

Let $\eta>0$ be as in the statement. Using \eqref{E:5.43ui} we now check that the second term on the right of \eqref{Eq:5.35u} is at most $O(\texte^{-\eta j}k^{-2})$, uniformly in~$j=2,\dots,k-1$. As to the first term, collecting the error terms in \eqref{E:5.40q}, \eqref{E:5.41q} and \eqref{E:5.47q} shows
\begin{equation}
\label{E:5.48q}
h_{k,j}(\alpha)\gamma_{k-1}(1)\gamma_{k-j}(1)\prod_{i=2}^{j-1}\gamma_{k-i}(0) = b^{-(1+2\alpha)(j-1)}(b-1)^2\frac{A^2}k \texte^{O(\tilde\fraku_{k,j})},
\end{equation}
where
\begin{equation}
\begin{aligned}
\tilde \fraku_{k,j}:=\log\frac k{k-j}+&
\sum_{i=1}^j|\sigma_{k-i}^2-1|+\frake_{k-1}+\frake_{k-j}+(k-1)^{-1/2}+(k-j)^{-1/2}
\\
&+\sum_{i=1}^j\frac1{k-i}+\sum_{i=1}^j\frac1{(k-i)^{3/2}}+\sum_{i=1}^j\frac{\frake_{k-i}}{k-i}.
\end{aligned}
\end{equation}
Now observe that the first, seventh and eighth term are order $\log(\frac{k}{k-j})$, while the second, third, fourth and the last term are bounded by~$(\frake_{k-1}+\frake_{k-j})[1+\log(\frac{k}{k-j})]$. Using that $\log(\frac{k}{k-j})=O(j/k)$ for~$j\le k/2$ and $(k-j)^{-1/2}\texte^{-\eta j}=O(k^{-1/2})$, it follows that
\begin{equation}
\texte^{-2\eta j}\tilde\fraku_{k,j}=O(\frake_{k-1}+\frake_{k-j})\texte^{-\eta j}+O(k^{-1/2}) \texte^{-\eta j}.
\end{equation}
This shows that the quantity in \eqref{E:5.48q} equals
\begin{equation}
b^{-(1+2\alpha)(j-1)}(b-1)^2\frac{A^2}k
+O(\frake_{k-1}+\frake_{k-j})k^{-1}\texte^{-\eta j} + O(k^{-3/2}) \texte^{-\eta j}.
\end{equation}
Plugging this in \eqref{Eq:5.35u}, we get \eqref{E:5.40w}.
\end{proofsect}

We will now show how to use the above to control the fractional-charge asymptotic when~$\beta=\beta_\cc$.  Set  $\tilde\tau:=\frac12\tau(\alpha)$ for $\tau(\alpha)$ from \eqref{E:tau} and observe that we then have 
\begin{equation}
\label{Eq:5.32i}
\tilde\tau+ 2(b-1)A^2 = \sum_{j=2}^\infty \bigl(b^{-(1+2\alpha)(j-1)}+b^{-(1-2\alpha)(j-1)}\bigr)(b-1)^2A^2
\end{equation}
for  $A$ as above.  Setting
\begin{equation}
\label{E:5.39u}
r_k:=k^{-\tilde\tau}\theta^{-\alpha^2 k}w_k(0),
\end{equation}
the proof of Theorem~\ref{thm-4} will boil down to showing that~$r_k$ is close to a positive and finite $n$-dependent constant once~$k$ and~$n-k$ are large. For this we first prove:

\begin{lemma}
\label{lemma-5.6}
Let~$\beta=\beta_\cc$  and suppose that the sequence $\{\frakd_j\}_{j\ge0}$ from Assumption~\ref{ass-1} obeys $\sum_{j\ge1}\frakd_j\log(j)<\infty$.  Then
\begin{equation}
\label{E:5.40u}
0<\inf_{n\ge k\ge1}r_k\le\sup_{n\ge k\ge1}r_k<\infty.
\end{equation}
\end{lemma}

\begin{proofsect}{Proof}
 For~$A$ as above, denote the quantity in the leading-order term in \eqref{E:5.40w} as 
\begin{equation}
h_j:= (b^{-(1+2\alpha)(j-1)}+b^{-(1-2\alpha)(j-1)})(b-1)^2A^2
\end{equation}
 and, recalling $\texte^{-\eta}:=b^{-\frac12(1-2|\alpha|)}$, abbreviate the errors in \twoeqref{E:5.39w}{E:5.40w} as 
\begin{equation}
\begin{aligned}
\fraks_k&:=k^{-3/2}+k^{-1}\frake_k+\frakd_{\min\{k-1,n-k+1\}},
 \\
\fraku_{k,j}&:=\texte^{-\eta j}\bigl[k^{-3/2}+k^{-1}(\frake_{k-1}+\frake_{k-j})\bigr].
\end{aligned}
\end{equation}
Then use the asymptotic forms from Lemma~\ref{lemma-5.7e} to cast \eqref{Eq:5.23} as
\begin{equation}
\label{Eq:5.36i}
\begin{aligned}
r_k = \Bigl(1-\frac1k\Bigr)^{\tilde\tau}&\texte^{-2(b-1)A^2 k^{-1}+O(\fraks_k)}r_{k-1}
\\
&+\sum_{j=2}^{k-1}h_j\frac1k\Bigl(1-\frac jk\Bigr)^{\tilde\tau} r_{k-j}
+O\biggl(\,\sum_{j=2}^{k-1}\fraku_{k,j} r_{k-j}\biggr)+O(\texte^{-\eta k}),
\end{aligned}
\end{equation}
where the last error term arises  from applying the bound \eqref{Eq:5.24} to the second sum in \eqref{Eq:5.23}. Next we denote~$M_k:=\max_{j=1,\dots,k}r_j$ and use it to dominate the terms in \eqref{Eq:5.36i} with the help of $1-\frac1k\le\texte^{-1/k}$ and $1-\frac jk\le1$ to get
\begin{equation}
M_k\le\biggl[\texte^{-[\tilde\tau+ 2(b-1)A^2]k^{-1}}+\frac1k\sum_{j\ge2}h_j+c\Bigl(\fraks_k+\sum_{j=2}^{k-1}\fraku_{k,j}\Bigr)\biggr]M_{k-1}+c\texte^{-\eta k}
\end{equation}
for some constant~$c>0$ independent of~$n\ge k>j\ge2$.
The choice of~${\tilde\tau}$ ensures (via \eqref{Eq:5.32i}) that the sum of the first two terms in the square brackets equals $1+O(k^{-2})$. To prove uniform boundedness of~$\{M_k\}_{k=1}^n$ it thus suffices to show that $\sum_{k=1}^n\fraks_k$ and $\sum_{k=1}^n\sum_{j=2}^{k-1}\fraku_{k,j}$ are bounded uniformly in~$n$. 
In light of Assumption~\ref{ass-1}, this reduces to uniform boundedness of~$\sum_{k=1}^n k^{-1}\frake_k$. Here we compute
\begin{equation}
\sum_{k=1}^n k^{-1}\frake_k \le 2\sum_{k\ge1}\frac1k\sum_{j\ge\sqrt k}\frakd_j
\le2\sum_{j\ge1}\frakd_j\sum_{1\le k\le j^2}\frac1k \le 4\sum_{j\ge1}\frakd_j\log(j^2)
\end{equation}
which we assumed to be finite.

Concerning the lower bound, we denote $m_k:=\min_{j=1,\dots,k}r_j$ and apply the inequality $w_{k-j}(0)\ge (k-j)^{\tilde\tau}\theta^{\alpha^2(k-j)}m_{k-1}$ to the first sum on the right of \eqref{Eq:5.23}. Dropping the second sum and invoking Lemma~\ref{lemma-5.7e} along with the bound $(1-\frac jk)^{\tilde\tau}\ge 1-{\tilde\tau} j/k$ then shows
\begin{equation}
\label{E:5.49u}
m_k\ge \biggl[\texte^{-[\tilde\tau+ 2(b-1)A^2]k^{-1}-c'\fraks_k}+\frac1k\sum_{j=2}^k h_j-\frac{\tilde\tau}{k^2}\sum_{j\ge2}j h_j -c'\sum_{j=2}^{k-1}\fraku_{k,j}\biggr]m_{k-1}
\end{equation}
for some constant~$c'>0$. Using \eqref{Eq:5.32i} we now check that there is~$k_0\ge1$ such that the term in the square bracket is uniformly positive for~$n\ge k\ge k_0$ and differs from~$1$ by a quantity that is uniformly summable on~$k=1,\dots, n$. It follows that~$m_k\ge c'' m_{k_0}$ for all~$n\ge k\ge k_0$. To extend the bound to~$k\le k_0$ we call upon Lemma~\ref{lemma-5.5b} which gives~$m_k\ge c'''k^{-\tilde \tau}$ for~$c'''>0$ independent of~$n\ge k\ge1$.
\end{proofsect}

We are now ready for:

\begin{proofsect}{Proof of Theorem~\ref{thm-2}, $\beta=\beta_\cc$}
We start by deriving a recursive bound on the difference $r_k-r_{k-1}$.
For this we only need to expand a bit on the arguments from the proof of Lemma~\ref{lemma-5.6}.
Indeed, using that~$\{r_k\}_{k=1}^n$ is bounded, we can trim \eqref{Eq:5.36i} to the form
\begin{equation}
r_k = r_{k-1}-\frac1 k \bigl[{\tilde\tau}+ 2 (b-1)A^2\bigr] r_{k-1}+\frac1k\sum_{j=2}^{k-1}h_jr_{k-j} +O(\frakv_k),
\end{equation}
 where~$\frakv_k :=k^{-2}+\texte^{-\eta k}+\fraks_k+\sum_{j=2}^{k-1}\fraku_{k,j}$.
For our choice of~${\tilde\tau}$ (and~$\eta$), this gives
\begin{equation}
r_k-r_{k-1} = \frac1k\sum_{j=2}^{k-1} h_j (r_{k-j}-r_{k-1})   + O(\frakv_k).
\end{equation}
With the help of the triangle inequality and a simple interchange of two sums we get
\begin{equation}
\label{Eq:5.42}
|r_k-r_{k-1}| \le \frac1k\sum_{i=1}^{k-2} \Bigl(\,\sum_{j\ge i+1} h_j\Bigr) |r_{k-i}-r_{k-i-1}|+ a\frakv_k
\end{equation}
for a constant~$a\ge0$ independent of~$n\ge k\ge1$. Setting
\begin{equation}
k_0:=\left\lceil\sum_{i\ge1}\texte^{\eta i}\Bigl(\,\sum_{j\ge i+1} h_j\Bigr)\right\rceil,
\end{equation}
where the sum over~$i$ is finite for~$\eta$ as above, we adjust~$a$ so that also
\begin{equation}
\label{Eq:5.35}
|r_j-r_{j-1}|\le a\sum_{i=0}^{j-2}\texte^{-\eta i}\frakv_{j-i}
\end{equation}
holds for all~$j=2,\dots,k_0$.

We now claim that \eqref{Eq:5.35} is true for all~$j=1,\dots,n$. Indeed, suppose $k\ge k_0$ is such that \eqref{Eq:5.35} holds for~$j=2,\dots,k-1$. Then plugging the bound in \eqref{Eq:5.42} yields
\begin{equation}
\begin{aligned}
|r_k-r_{k-1}|&\le a \frakv_k +\frac ak\sum_{i=1}^{k-2} \Bigl(\,\sum_{j\ge i+1} h_j\Bigr)\sum_{\ell=i}^{k-2}\texte^{-\eta(\ell-i)}\frakv_{k-\ell}
\\&= a \frakv_k +\frac ak\sum_{\ell=1}^{k-2}\biggl[\,\sum_{i=1}^\ell \texte^{\eta i}\Bigl(\,\sum_{j\ge i+1} h_j\Bigr)\biggr] \,\texte^{-\eta\ell}\frakv_{k-\ell}.
\end{aligned}
\end{equation}
Noting that the quantity in the square bracket is bounded by~$k_0$, the fact that~$k\ge k_0$ then implies \eqref{Eq:5.35} for~$j:=k$. This proves \eqref{Eq:5.35} for all~$j=2,\dots,n$ by induction.

With \eqref{Eq:5.35} in hand, we proceed similarly as in the subcritical situations. Indeed, abbreviate $n':=\lfloor n/2\rfloor$, set~$C_n:=r_{n'}^2$ and note that, by Lemma~\ref{lemma-5.6}, $C_n$ is positive and finite uniformly in~$n\ge1$. The inequality \eqref{Eq:5.35} implies
\begin{equation}
\bigl|r_k-\sqrt{C_n}\bigr|\le \frac a{1-\texte^{-\eta}}\sum_{j=\min\{k,n-k\}}^{\max\{k,n-k\}}\frakv_j+ \frac a{1-\texte^{-\eta}}\sum_{j=1}^{k-1}\texte^{-\eta j}\frakv_{k-j}.
\end{equation}
The assumption $\sum_{j\ge1}\frakd_j\log(j)<\infty$ along with the bounds in the proof of Lemma~\ref{lemma-5.6} then give~$r_k-\sqrt{C_n}\to0$ as~$\min\{k,n-k\}\to \infty$.
Using this in \eqref{E:5.39u} along with $2\tilde\tau=\tau(\alpha)$, Lemma~\ref{cor-5.5} and standard facts about Fourier series show
\begin{equation}
\bigl|f_k(z)\bigr|^2= \bigl[C_n+o(1)\bigr] k^{\tau(\alpha)}\theta^{2\alpha^2 k} 
\end{equation}
with~$o(1)\to0$ as~$\min\{k,n-k\}\to\infty$ uniformly in~$z\in\R$.  Plugging~$k:=k(x,y)$ we get \eqref{E:1.17} for~$\beta=\beta_\cc$.
\end{proofsect}

\subsection{Above criticality}
Our last item of business in this section is the asymptotic of the fractional charge for~$\beta$ slightly above~$\beta_\cc$. Throughout we assume that the sequence~$\{\frakd_k\}_{k\ge0}$ in Assumption~\ref{ass-1} exhibits exponential decay. 

 We will again rely on the representation from Lemma~\ref{lemma-5.6a} for which we need to identify the asymptotic values of the coefficients $\Gamma_{k,j}(p,0)$. Recall the notation $\Xi_b(q)$ from \eqref{E:5.73i} and, with~$\{\lambda_\star(q)\}_{n\in\Z}$  as in Theorem~\ref{thm-5}, set
\begin{equation}
\gg_\star(q):=\frac{\sum_{\bar\ell\in\Xi_{b-1}(q)}\prod_{i=1}^{b-1}\lambda_\star(\ell_i)}{\sum_{\bar\ell'\in\Xi_b(0)}\prod_{i=1}^{b}\lambda_\star(\ell_i')},
\end{equation}
where $\bar\ell$ stands for $(\ell_1,\dots,\ell_{b-1})$ and $\bar\ell'$ stands for $(\ell_1',\dots,\ell_{b}')$.
Note that, by \eqref{E:3.28u}, the sums are finite for~$b\theta-1$ small and $\gg_\star(p)=\gg_\star(-p)$ by $\lambda_\star(q)=\lambda_\star(-q)$. Now let
\begin{equation}
\Gamma^\star_1(p):=\gg_\star(p)
\end{equation}
and, for~$j\ge2$, let
\begin{equation}
\label{Eq:5.21a}
\Gamma^\star_{j}(p):=\sum_{\begin{subarray}{c}
q_1,\dots,q_{j-1}\in\Z\smallsetminus\{0\}\\ q_0=p,\,q_j=0
\end{subarray}}
\biggl(\,\prod_{i=0}^{j-1}\theta^{(q_i+\alpha)^2-\alpha^2}\biggr)\prod_{i=1}^{j}\gg_\star(q_i-q_{i-1}).
\end{equation}
Note that these are the quantities in \eqref{E:5.29ui} and \eqref{Eq:5.19}, respectively, corresponding to the renormalization-group flow fixed point. 
Theorem~\ref{thm-5} then shows:

\begin{lemma}
\label{lemma-5.7}
For all~$\alpha\in(-\ffrac12,\ffrac12)$ and~$\delta\in(0,1-2|\alpha|)$ there exist~$\epsilon>0$ and $C'>0$ and, for all~$\beta>\beta_\cc$ with~$1/\beta>1/\beta_\cc-\epsilon$, there exist $C>0$  and  $\eta>0$ such that
\begin{equation}
\label{E:5.49ui}
\Gamma^\star_j(p)\le C'\sqrt{b\theta-1}\,\,\,b^{-(1-2|\alpha|-\delta)j}\texte^{-\eta|p|(|p|-1)}
\end{equation}
holds for all~$j\ge2$ and $p\in\Z$ and
\begin{equation}
\label{E:5.50ui}
\bigl|\Gamma_{k,j}(p,0)-\Gamma^\star_j(p)\bigr|\le Cb^{-(1-2|\alpha|-\delta)j}\texte^{-\eta|p|(|p|-1)}\texte^{-\eta \min\{k,n-k\}}
\end{equation}
holds for $n\ge k> j\ge1$ and~$p\in\Z$.
\end{lemma}

\begin{proofsect}{Proof}
We start by estimates for the weights~$\gamma_k$ and~$\gamma_\star$. For~$b\theta-1$ sufficiently small, the bounds \twoeqref{E:3.28u}{Eq:3.30w} imply that, for some constant~$C'>0$,
\begin{equation}
\label{E:5.76t}
\max\biggl\{\sum_{\ell\in\Z}\gg_k(\ell), \sum_{\ell\in\Z}\gg_\star(\ell)\biggr\}\le 1+C'\sqrt{b\theta-1}
\end{equation}
holds uniformly in~$1\le k\le n$. The assumption of exponential decay of $\{\frakd_k\}_{k\ge0}$ in turn allows us to summarize the inequalities 
\twoeqref{E:3.28u}{Eq:3.31} as
\begin{equation}
\label{E:5.75t}
\Bigl|\frac{a_k(q)}{a_k(0)}-\lambda_\star(q)\Bigr|\le C''\texte^{-2\eta' \max\{\min\{k,n-k\},|q|\}},
\end{equation}
for some~$\eta'>0$ that may depend on~$\beta$. Assuming~$\eta'$ is small that the quantity on the right of \twoeqref{E:3.28u}{Eq:3.30w} obeys $2b^{1/2}\sqrt{b\theta-1}\le\texte^{-\eta'}$, we have the telescopic bound
\begin{equation}
\begin{aligned}
\sum_{q\in\Z}\biggl|\,&\sum_{\bar\ell\in\Xi_{b}(q)}
\prod_{i=1}^{b}\frac{a_k(\ell_i)}{a_k(0)}-\sum_{\bar\ell\in\Xi_{b}(q)}\prod_{i=1}^{b}\lambda_\star(\ell_i)\biggr|
\\
&\le\sum_{q\in\Z}\sum_{\bar\ell\in\Sigma_{b}(q)}\sum_{j=1}^{b}
\biggl(\,\prod_{i=1}^{j-1}\lambda_\star(\ell_i)\biggr)\biggl(\prod_{i=j+1}^b\frac{a_k(\ell_i)}{a_k(0)}\biggr)\biggl|\frac{a_k(\ell_i)}{a_k(0)}-\lambda_\star(\ell_i)\biggr|
\\
&\le b\biggl(\frac{1+\texte^{-\eta'}}{1-\texte^{-\eta'}}\biggr)^{b-1}C''\biggl(1+2\min\{k,n-k\}+\frac2{1-\texte^{-2\eta'}}\biggr)\texte^{-2\eta'\min\{k,n-k\}}.
\end{aligned}
\end{equation}
Applying this in the numerator and denominator of the ratios defining~$\gamma_k(q)$ and~$\gamma_\star(q)$ along with \eqref{E:5.76t} gives
\begin{equation}
\label{E:5.63e}
\sum_{q\in\Z}
\bigl|\gg_k(q)-\gg_\star(q)\bigr|\le C'''\texte^{-\eta' \min\{k,n-k\}}.
\end{equation}
for some~$C'''<\infty$. We will assume that~$\eta'\le\frac12\delta\log b$.

Now assume that~$\epsilon>0$ is so small that, for~$\beta$ with $1/\beta>1/\beta_\cc-\epsilon$ we have
\begin{equation}
(b\theta)^{1-2|\alpha|}\sum_{\ell\in\Z}\gg_\star(\ell)\le b^{\delta}
\end{equation}
which is possible in light of \eqref{E:5.76t}. 
The reasoning underlying \eqref{Eq:5.24} gives
\begin{equation}
\Gamma^\star_j(p)\le \theta^{\sigma_{\text{min}}^2|p|(|p|-1)+(1-2|\alpha|)(j-1)}\biggl(\,\sum_{q\in\Z\smallsetminus\{0\}}\gg_\star(q)\biggr)\biggl(\,\sum_{\ell\in\Z}\gg_\star(\ell)\biggr)^{j-1},
\end{equation}
where we noted that~$q_j-q_{j-1}\ne0$ in \eqref{Eq:5.21a}. Since the first sum on the right is order $\sqrt{b\theta-1}$ by the same argument that proved \eqref{E:5.76t}, we get \eqref{E:5.49ui} with~$\texte^{-\eta}:=\theta^{\sigma_{\text{min}}^2}$. 

In order to prove \eqref{E:5.50ui}, we telescopically swap the~$k$-dependent terms in the  expression  for~$\Gamma_{k,j}(0)$ for the corresponding terms in~$\Gamma^\star_j$. Using  \eqref{E:5.35ui} this gives
\begin{equation}
\label{E:5.62ua}
\begin{aligned}
\bigl|\Gamma_{k,1}&(p,0)-\Gamma^\star_1(p)\bigr|
\\
&\le \theta^{-\alpha^2+\sigma_{\text{min}}^2|p|(|p|-1)}\bigl|\gg_{k-1}(p)-\gg_\star(p)\bigr|+\bigl|\theta_k^{(p+\alpha)^2}\theta^{-\alpha^2}-\theta^{(p+\alpha)^2}\theta^{-\alpha^2}\bigr|\,\gg_\star( p).
\end{aligned}
\end{equation}
For $j\ge2$, we in turn invoke the bound \eqref{E:5.37w} and recall~$\tilde c$ defined below \eqref{Eq:5.24} to get
\begin{equation}
\label{E:5.63ua}
\bigl|\Gamma_{k,j}(p,0)-\Gamma^\star_j(p)\bigr|\le \theta^{-\tilde c+  (1-2|\alpha|)(j-1)+\sigma_{\text{min}}^2|p|(|p|-1)} B_1+\theta^{ -\tilde c+(1-2|\alpha|)(j-2)} B_2,
\end{equation}
where
\begin{equation}
\label{E:5.64ua}
 B_1:=\sum_{m=1}^j\biggl(\,\prod_{i=1}^{m-1}\sum_{\ell\in\Z}\gg_{k-i}(\ell)\biggr)\biggl(\,\sum_{\ell\in\Z}\bigl|\gg_{k- m}(\ell)-\gg_\star(\ell)\bigr|\biggr)\biggl(\,\sum_{\ell\in\Z}\gg_\star(\ell)\biggr)^{j- m}
\end{equation}
and
\begin{equation}
\label{E:5.65ua}
\begin{aligned}
 B_2:=\biggl(\,\sum_{\ell\in\Z}\gg_\star(\ell)\biggr)^{j}\Biggl[\bigl|\theta_{k-1}^{(p+\alpha)^2}&\theta^{-\alpha^2}-\theta^{(p+\alpha)^2}\theta^{-\alpha^2}\bigr|
\\&+\theta^{(p+\alpha)^2-\alpha^2}
\sum_{i=1}^{j-1}\sum_{q\in\Z}\bigl|\theta_{k-i-1}^{(q_i+\alpha)^2}\theta^{-\alpha^2}-\theta^{(q+\alpha)^2}\theta^{-\alpha^2}\bigr|\Biggr].
\end{aligned}
\end{equation}
The bounds \twoeqref{E:5.76t}{E:5.63e} give 
\begin{equation}
B_1\le C'''\texte^{\eta' j}\texte^{-\eta' \min\{k,n-k\}} j [1+C'\sqrt{b\theta-1}]^{j-1}.
\end{equation}
 If~$b\theta-1$ is so small that~$ 1+C'\sqrt{b\theta-1}\le b^{\delta/3}$, then the assumption $\eta'\le\frac12\delta\log b$ along with the fact that~$\sup_{j\ge1} j b^{-j\delta/6}<\infty$ bounds the first term on the right of \eqref{E:5.63ua} by a quantity proportional to the right-hand side of \eqref{E:5.50ui}. 

The first term on the right of \eqref{E:5.62ua} is bounded via \eqref{E:5.63e} so it remains to bound the second terms in \eqref{E:5.62ua} and \eqref{E:5.63ua}. Here the term involving~$\gamma_\star$ is controlled with the help of \eqref{E:5.76t} so it remains to estimate the quantity in absolute value. The elementary inequality $|\texte^{-a}-\texte^{-\tilde a}|\le \texte^{-\min\{a,\tilde a\}}|a-\tilde a|$ combine into 
\begin{equation}
\label{E:5.83u}
\bigl|\theta_\ell^{(q+\alpha)^2}\theta^{-\alpha^2}-\theta^{(q+\alpha)^2}\theta^{-\alpha^2}\bigr|
\le \log(1/\theta)\,\theta^{\min\{1,\sigma_\ell^2\}(q+\alpha)^2 -\alpha^2}(q+\alpha)^2|\sigma_\ell^2-1|.
\end{equation}
Assumption~\ref{ass-1} and the exponential decay of~$\{\frakd_j\}_{j\ge0}$ give
\begin{equation}
\label{E:5.84u}
|\sigma_{\ell}^2-1|\le \wt C\texte^{-\wt \eta\min\{\ell,n-\ell\}}
\end{equation}
 for some~$\wt C,\wt\eta>0$. In light of $\inf_{n\ge \ell\ge 1}\sigma_\ell^2>0$, this bounds the second term on the right of \eqref{E:5.62ua} by a quantity proportional to the right-hand side of \eqref{E:5.50ui}.

The bounds \twoeqref{E:5.83u}{E:5.84u} dominate the first term in the square bracket in \eqref{E:5.65ua} by a constant times $\theta^{\frac12\min\{\sigma_{\text{min}}^2,1\}|p|(|p|-1)}\texte^{-\wt\eta\min\{k,n-k\}}$, where ``half'' of the exponential decay in~$|p|$ was used to control the term~$(p+\alpha)^2$. Similarly, the second term in the square bracket in \eqref{E:5.65ua} is bounded by the same quantity as the first times~$j\texte^{\wt\eta j}$. Summarizing,
\begin{equation}
B_2\le \wt C'\theta^{\frac12\min\{\sigma_{\text{min}}^2,1\}|p|(|p|-1)}\texte^{-\wt\eta\min\{k,n-k\}}[1+C'\sqrt{b\theta-1}]^j\bigl[1+j\texte^{\wt\eta j}\bigr],
\end{equation}
where \eqref{E:5.76t} was used for the~$\gamma_\star$-dependent prefactor.
Assuming that~$\wt\eta<\delta/2$ and $b\theta-1$ so small that $1+C'\sqrt{b\theta-1}\le b^{\delta/6}$, the last two terms on the right are at most a constant times $b^{\delta j}$. Inserting this on the right of \eqref{E:5.63ua}, we get the claim.
\end{proofsect}

Given~$\beta>\beta_\cc$ for which Theorem~\ref{thm-5} applies and $\alpha\in(-\ffrac12,\ffrac12)$, define $t_\star = t_\star(\alpha,\beta)$ by
\begin{equation}
t_\star:=\inf\biggl\{t>0\colon \sum_{j\ge1}\Gamma^\star_j(0) t^{-j}\le1\biggr\}.
\end{equation}
Clearly, $t_\star\in(0,\infty)$ and, by Fatou's lemma, $\sum_{j\ge1}\Gamma^\star_j(0) t_\star^{-j}\le1$. Key for our use of this quantity is the fact that equality holds. Indeed, we have:

\begin{lemma}
\label{lemma-5.8}
 For each~$\alpha_0\in(0,\ffrac12)$ there exists~$\tilde\epsilon>0$ such that 
\begin{equation}
\label{Eq:5.47u}
\sum_{j\ge1}\Gamma^\star_j(0) t_\star^{-j}=1
\end{equation}
holds true for all~$\beta>\beta_\cc$ with~$1/\beta>1/\beta_\cc-\tilde\epsilon$ and all~$\alpha\in[0,\alpha_0]$. Moreover, for such~$\alpha$ and~$\beta$ we have $t_\star\ge1$ with $t_\star>1$ when~$\alpha\ne0$.
\end{lemma}

\begin{proofsect}{Proof}
Plugging the explicit form of~$\Gamma^\star_j(0)$ yields
\begin{equation}
\label{E:5.95}
\begin{aligned}
\sum_{j\ge1}&\,\Gamma^\star_j(0)t^{-j}=\gamma_\star(0)t^{-1}
\\
&\quad+\sum_{j\ge2}\,t^{-j}\sum_{\begin{subarray}{c}
q_1,\dots,q_{j-1}\in\Z\smallsetminus\{0\}\\ q_0=0,\,q_j=0
\end{subarray}}
\cosh\biggl(2\alpha  \log(\theta) \sum_{i=1}^{j-1}q_i\biggr)
\biggl(\,\prod_{i=1}^{j-1}\theta^{q_i^2}\biggr)\prod_{i=1}^{j}\gg_\star(q_i-q_{i-1}),
\end{aligned}
\end{equation}
where we also used~$\gamma_\star(-\ell)=\gamma_\star(\ell)$ to symmetrize the expression with respect to the simultaneous reflection of all~$q_0,\dots,q_j$. Next we observe that the same reasoning as used for \eqref{E:5.76t} implies
\begin{equation}
\label{E:5.96}
\sum_{q\in\Z}\frac{\gamma_\star(q)}{\gamma_\star(0)}\le 1+C\sqrt{b\theta-1}
\end{equation}
for some~$C>0$. Proceeding as in the derivation of \eqref{E:5.37w} to bound the $\theta$-dependent factor in \eqref{E:5.95} by $\theta^{(1-2|\alpha|)(j-1)}$, the sum in \eqref{E:5.95} converges for~$t:=\gamma_\star(0)$ once $\beta-\beta_\cc$ is so small that
\begin{equation}
\label{E:5.97}
\theta^{1-2|\alpha_0|}(1+C\sqrt{b\theta-1})<1.
\end{equation}
 Under this condition the sum converges absolutely, and is thus continuous in $t\ge\gamma_\star(0)$. Since the whole quantity on the right also trivally exceeds~$1$ at~$t:=\gamma_\star(0)$, equality \eqref{Eq:5.47u} holds at~$t=t_\star$.

It remains to show that~$t_\star\ge1$ with~$t_\star>1$ when~$\alpha\ne0$. For this we observe that~$\lambda_\star$ solves the recursion \eqref{lambda-rec}
which by partitioning according to whether~$\ell_b=0$ or not rewrites into
\begin{equation}
u_\star(p) = \gamma_\star(p)+\sum_{q\ne 0}\theta^{q^2}\gamma_\star(p-q)u_\star(q)
\end{equation}
for the quantity $u_\star(p):=\theta^{-p^2}\lambda_\star(p)$. Iterating $j$-times and taking $j\to\infty$ while controlling the error using the argument \twoeqref{E:5.96}{E:5.97} and $u_\star(q)\le1$ we readily get
\begin{equation}
1 = u_\star(0)=\gamma_\star(0)+ 
\sum_{j\ge2}  \,\sum_{\begin{subarray}{c}
q_1,\dots,q_{j-1}\in\Z\smallsetminus\{0\}\\ q_0=0,\,q_j=0
\end{subarray}}
\biggl(\,\prod_{i=1}^{j-1}\theta^{q_i^2}\biggr)\prod_{i=1}^{j}\gg_\star(q_i-q_{i-1}),
\end{equation}
where we used that~$u_\star(0)=1$ on the right as well. This shows that the right-hand side of \eqref{E:5.95} is at least~$1$ when~$t=1$ and it exceeds~$1$ strictly when also~$\alpha\ne0$. Hence we must have~$t_\star\ge1$ with~$t_\star>1$ when~$\alpha\ne0$. 
\end{proofsect}
\eMB

We now follow the same blueprint as in the critical case. Define $\{\tilde r_k\}_{k=1}^n$ by
\begin{equation}
\label{E:5.58w}
\tilde r_k:= (t_\star\theta^{\alpha^2})^{-k}w_k(0).
\end{equation}
This quantity depends also on~$n$ but we keep that dependence implicit.
Our aim is to show that $\tilde r_k$ is close to an $n$-dependent constant once~$\min\{k,n-k\}$ is sufficiently large. As in the critical regime, we first prove:

\begin{lemma}
\label{lemma-5.9}
For each~$\alpha\in(-\ffrac12,\ffrac12)$ there exists~$\tilde \epsilon>0$ such that for all~$\beta>\beta_\cc$ satisfying $1/\beta>1/\beta_\cc-\tilde\epsilon$ we have
\begin{equation}
\label{E:5.69ee}
0<\inf_{n\ge k\ge1}\tilde r_k\le\sup_{n\ge k\ge1}\tilde r_k<\infty.
\end{equation}
\end{lemma}

\begin{proofsect}{Proof}
Applying \eqref{E:5.58w} in \eqref{Eq:5.23} yields
\begin{equation}
\label{E:5.91q}
\tilde r_k = \sum_{j=1}^{k-1}\Gamma_{k,j}(0,0)t_\star^{-j}\tilde r_{k-j}+\sum_{q\in\Z\smallsetminus\{0\}}\Gamma_{k,k-1}(0,q)t_\star^{-k}\theta^{-\alpha^2}w_1(q).
\end{equation}
An inspection of \eqref{Eq:2.8} shows that~$w_1$ is bounded. Since $t_\star\ge1$, the uniform exponential decay \eqref{Eq:5.24} implies that the last term is~$O(\texte^{-\eta' k})$ for some~$\eta'>0$. 
Invoking Lemma~\ref{lemma-5.7} to swap $\Gamma_{k,j}(0,0)$ for~$\Gamma_j^\star(0)$ gives
\begin{equation}
\label{Eq:5.51}
\tilde r_k \le O(\texte^{-\eta' k})+ \sum_{j=1}^{k-1}\Gamma_j^\star(0) t_\star^{-j}\tilde r_{k-j} +  C \texte^{-\eta\min\{k,n-k\}}\sum_{j=1}^{k-1} b^{-(1-2|\alpha|-\delta)j}t_\star^{-j}\tilde r_{k-j}.
\end{equation}
Setting $M_k:=\max_{1\le j\le k}\tilde r_j$, the same reasoning applied to the second sum along with the bound $\sum_{j=1}^{k-1}\Gamma^\star_j(0) t_\star^{-j}\le1$ yields
$M_k\le  O(\texte^{-\eta' k})+(1+O(\texte^{-\eta\min\{k,n-k\}}))M_{k-1}$
 with the implicit constants uniform in~$n$. This now gives the upper bound in \eqref{E:5.69ee}.

For the complementary direction we first prove that~$\{\tilde r_k\}_{k\ge0}$ cannot decay exponentially fast. Indeed, for this we pick~$\delta>0$ and note that Lemma~\ref{lemma-5.7} along with \eqref{Eq:5.47u} show that, for some~$\ell\ge1$ and $k_0>\ell$, 
\begin{equation}
\sum_{j=1}^\ell \Gamma_{k,j}(0,0)t_\star^{-j}\ge\texte^{-\delta}
\end{equation}
holds once~$\min\{k,n-k\}\ge k_0$.
Now observe that plugging \eqref{E:5.58w} for~$w_j(0)$ in \eqref{Eq:5.23} and retaining only the terms with~$j\le\ell$ from the first sum (and dropping the second sum)~yields
\begin{equation}
\label{E:5.63}
\tilde r_k \ge \sum_{j=1}^{\ell}\Gamma_{k,j}(0,0) t_\star^{-j}\tilde r_{k-j}.
\end{equation}
Setting $m_k:=\min_{1\le j\le k}\tilde r_j$ we get $m_k\ge \texte^{-\delta}m_{k-1}$ once~$\min\{k,n-k\}\ge k_0$. For~$k$ violating this inequality, we in turn use that $\tilde r_k\ge \Gamma_{k,1} (0)m_{k-1}$ and observe that the product of $\min\{\Gamma_{k,1} (0,0),1\}$ for~$k=2,\dots,k_0$ and $k=n-k_0,\dots,n$ is positive uniformly in~$n\ge1$. Writing~$c$ for this product we conclude that $m_k\ge c\texte^{-\delta(k-k_0)}\tilde r_1$. As~$\delta>0$ is arbitrary, we have ruled out exponential decay.

We now redo the argument leading to \eqref{E:5.63} while invoking Lemma~\ref{lemma-5.7} and the boundedness of~$\{\tilde r_k\}_{k=1}^n$ proved earlier to get
\begin{equation}
\tilde r_k \ge \sum_{j=1}^{k-1}\Gamma_j^\star(0) t_\star^{-j}\tilde r_{k-j} -O(\texte^{-\eta \min\{k,n-k\}})
\end{equation}
The boundedness of~$\{\tilde r_k\}_{k\ge1}$ along with the exponential decay \eqref{E:5.49ui} allows us to extend the sum all the way to infinity at the cost of an $O(\texte^{-\eta' k})$ error. From \eqref{Eq:5.47u} we then get $m_k\ge m_{k-1}-a\texte^{-\eta \min\{k,n-k\}}$ for some constant~$a>0$. But the fact that $\{m_k\}_{k\ge1}$ does not decay exponentially means that we can wrap this as $m_k\ge (1-a\texte^{-\eta k/2})m_{k-1}$ once~$k$ is sufficiently large. Any positive sequence satisfying this recursive bound is necessarily bounded away from zero.
\end{proofsect}

Next we prove an iterative bound on the increments of $\{\tilde r_k\}_{k=1}^n$:

\begin{lemma}
For each~$\alpha\in(-\ffrac12,\ffrac12)$ there exists~$\tilde \epsilon>0$ and~$\wt C,\wt\eta>0$ such that for all~$\beta>\beta_\cc$ satisfying $1/\beta>1/\beta_\cc-\tilde\epsilon$,
\begin{equation}
\label{E:5.67a}
|\tilde r_k-\tilde r_{k-1}|\le \wt C\texte^{-\wt\eta \min\{k,n-k\}}
\end{equation}
holds for all~$n\ge k\ge2$.
\end{lemma}

\begin{proofsect}{Proof}
We start by using \eqref{E:5.50ui} in \eqref{E:5.91q} along with \eqref{Eq:5.24}, the boundedness of~$\{\tilde r_k\}_{k=1}^{n+1}$ and the fact that $t_\star\ge1$ to get
\begin{equation}
\tilde r_k = O(\texte^{-\eta \min\{k,n-k\}})+\sum_{j=1}^{k-1}\Gamma_j^\star(0) t_\star^{-j}\tilde r_{k-j}.
\end{equation}
Next we invoke \eqref{Eq:5.47u} to rewrite this as
\begin{equation}
\tilde r_k-\tilde r_{k-1} = O(\texte^{-\eta \min\{k,n-k\}})+\sum_{j=2}^{k-1}\Gamma_j^\star(0)  t_\star^{-j}(\tilde r_{k-j}- \tilde r_{k-1}),
\end{equation}
where we also noted that the~$j=1$ term cancels on the right-hand side.
Using the same argument as in \eqref{Eq:5.42}, this yields
\begin{equation}
\label{E:5.66}
|\tilde r_k-\tilde r_{k-1}|\le \sum_{i=1}^{k-2}\biggl(\,\sum_{j=i+1}^{k-1} \Gamma^\star_j(0)t_\star^{-j} \biggr) | \tilde r_{k-i}-\tilde r_{k-i-1}| + a\texte^{-\eta \min\{k,n-k\}}
\end{equation}
for some constant~$a>0$ which we will for convenience assume exceeds $\texte^{2\eta}|\tilde r_2-\tilde r_1|$. 

Next observe that, since $t_\star\ge1$, for all~$j\ge2$ we have $\Gamma^\star_j(0) t_\star^{-j}\le C'\sqrt{b\theta-1}\,\texte^{-\eta' j}$
with~$C'$ a constant and~$\eta'$ close to $(1-2|\alpha|)\log b$ by \eqref{E:5.49ui}. Abbreviating $\eta'':=\min\{\eta,\eta'\}$, the bound
\begin{equation}
\label{E:5.67}
|\tilde r_{\ell}-\tilde r_{\ell-1}|\le 2a\texte^{-\frac12\eta''\min\{\ell,n-\ell\}},\quad \ell=2,\dots,k-1,
\end{equation}
then iterates via \eqref{E:5.66} to
\begin{equation}
\begin{aligned}
|\tilde r_k-\tilde r_{k-1}|
&\le 2a C'\sqrt{b\theta-1} \sum_{i=1}^{k-2}\sum_{j=i+1}^{k-1}  \texte^{-\eta' j}  \texte^{-\frac12\eta''\min\{k-i,n-k+i\}} +a\texte^{-\eta''\min\{k,n-k\}}
\\
&\le\biggl(\frac{2C'\sqrt{b\theta-1}}{ (1-\texte^{-\eta'/2})^2}+1\biggr) a\texte^{-\frac12\eta'' \min\{k,n-k\}},
\end{aligned}
\end{equation}
where we also used 
$\min\{k-i,n-k+i\}\ge\min\{k,n-k\}-i$ and applied $\eta'\ge\eta''$.
Noting that~$C'$ and~$\eta'$ do not depend on~$\beta$, for $b\theta-1$ so small that $2C'\sqrt{b\theta-1}\le (1-\texte^{-\eta'/2})^2$ we proved \eqref{E:5.67} for~$\ell:=k$ from \eqref{E:5.67} for~$\ell<k$. Since \eqref{E:5.67} holds for~$\ell:=2$ by our assumption on~$a$, it holds for all~$\ell\ge1$ by induction.
\end{proofsect}

We are now ready for: 

\begin{proofsect}{Proof of Theorem~\ref{thm-2}, $\beta>\beta_\cc$}
For each~$n\ge1$ abbreviate $n':=\lfloor n/2\rfloor$. We start by noting that the exponential decay \eqref{E:5.67a} implies
\begin{equation}
\label{E:5.83ui}
|\tilde r_k-\tilde r_{n'}|\le\frac{\wt C\texte^{\wt \eta}}{1-\texte^{-\wt \eta}}\,\texte^{-\wt\eta\min\{k,n-k\}},\quad  1\le k\le n.
\end{equation}
Next we note that the identity \eqref{Eq:5.23} rewrites via \eqref{E:5.58w} as
\begin{equation}
(t_\star\theta^{\alpha^2})^{-k}w_k(p)= \sum_{j=1}^{k-1}  \Gamma_{k,j}(p,0)t_\star^{-j}\tilde r_{k-j}+\sum_{q\in\Z\smallsetminus\{0\}}\Gamma_{k,k-1}(p,q)t_\star^{-k}\theta^{-\alpha^2}w_1(q).
\end{equation}
Denoting
\begin{equation}
 \tilde r_{n'}(p):=\tilde r_{n'}\sum_{j\ge1}\Gamma^\star_j(p)t_\star^{-j}  
\end{equation}
we thus get
\begin{equation}
\begin{aligned}
\bigl| (t_\star\theta^{\alpha^2})^{-k}&w_k(p)-\tilde r_{n'}(p)\bigr|
\le 
 \sum_{j=1}^{k-1}\bigl|\Gamma_{k,j}(p,0)-\Gamma^\star_j(p)\bigr|t_\star^{-j}\tilde r_{k-j}
+\tilde r_{n'}\sum_{j\ge k} \eHY \Gamma_j^\star(p) t_\star^{-j}
\\
&+  \sum_{j=1}^{k-1}  \Gamma_j^\star(p)t_\star^{-j}|\tilde r_{k-j}-\tilde r_{n'}|
+\sum_{q\in\Z\smallsetminus\{0\}}\Gamma_{k,k-1}(p,q)t_\star^{-k}\theta^{-\alpha^2}w_1(q).
\end{aligned}
\end{equation}
Using that $\{\tilde r_k\}_{k=1}^n$ is bounded uniformly in~$n$, invoking the bounds \twoeqref{E:5.49ui}{E:5.50ui} in the  first three  terms, the decay \eqref{E:5.83ui} in the third term and the bound \eqref{Eq:5.24} in the last term along with the fact that $t_\star\ge1$ shows
\begin{equation}
\label{E:5.88a}
\bigl| (t_\star\theta^{\alpha^2})^{-k}w_k(p)-\tilde r_{n'}(p)\bigr|\le C'\texte^{-\eta'|p|(|p|-1)}\texte^{-\eta'\min\{k,n-k\}}
\end{equation}
for some~$C',\eta'>0$ independent of~$n$ provided~$b\theta-1$ is sufficiently small.
 
We now define
\begin{equation}
\label{E:5.107}
f_\star(z):=\texte^{\tilde v_\star(z)}\sum_{p\in\Z}\sum_{j\ge1}\Gamma^\star_j(p)t_\star^{-j}\texte^{2\pi\texti p z},
\end{equation}
where 
\begin{equation}
\label{E:tilde-v_star}
\texte^{-\tilde v_\star(z)}:=\sum_{q\in\Z}\lambda_\star(q)\texte^{2\pi\texti qz}
\end{equation}
and where the sum in \eqref{E:5.107}
converges absolutely thanks to \eqref{E:5.49ui} and the fact that $t_\star\ge1$. The definition of~$f_k$ in \eqref{Eq:5.2b} then shows
\begin{equation}
\label{E:5.115}
\begin{aligned}  
 \bigl| (t_\star\theta^{\alpha^2})^{-k}& f_k(z)- \tilde r_{n'} f_\star(z)\bigr|
\\
&\le  \tilde r_{n'}\bigl|\texte^{- \tilde v_\star(z)}-a_{k-1}(0)^{-1}\texte^{-v_{k-1}(z)}\bigr| a_{k-1}(0)\texte^{v_{k-1}(z)}\bigl| f_\star(z)\bigr| 
\\
&
\qquad+a_{k-1}(0)\texte^{v_{k-1}(z)}\sum_{p\in\Z}\bigl| (t_\star\theta^{\alpha^2})^{-k}w_k(p)-\tilde r_{n'}(p)\bigr|.
\end{aligned}
\end{equation}
 \eHY 
Now observe that
\begin{equation}
 \bigl|\texte^{-\tilde v_\star(z)}  -a_{k-1}(0)^{-1}\texte^{-v_{k-1}(z)}\bigr|\le\sum_{q\in\Z}\Bigl|\frac{a_{k-1}(q)}{a_{k-1}(0)}-\lambda_\star(q)\Bigr|
\end{equation}
With the help of \eqref{E:5.88a} and \eqref{E:5.75t}, the quantity in \eqref{E:5.115} is thus at most $C'\texte^{-\eta'\min\{k,n-k\}}$. 
 Recalling that $\Vert P(\phi_k\text{ mod }1\in\cdot)-\nu_\star\Vert_{\text{TV}}$ decays exponentially in~$\min\{k,n-k\}$ (see the text before \eqref{E:4.45i}), the fact that~$\tilde r_{n'}$ and~$f_\star$ are also bounded shows 
\begin{equation}
\label{E:5.109}
  \Bigl| (t_\star\theta^{\alpha^2})^{-2k} E\bigl(f_k(\phi_k)^2\bigr)-\tilde r_{n'}^2E_\star\bigl(f_\star(\phi)^2\bigr) \eHY
\Bigr|\le C''\texte^{-\eta\min\{k,n-k\}}
\end{equation}
for some constant~$C''>0$.

 Set $C_n:=\tilde r_{n'}^2E_\star\bigl(f_\star(\phi)^2)$. Plugging \eqref{E:5.109}  in Corollary~\ref{cor-5.3}, for all distinct $x,y\in\Lambda_n$ we thus obtain
\begin{equation}
\langle\texte^{2\pi\texti\alpha(\phi_x-\phi_y)}\rangle_{n,\beta} = \bigl[
C_n+O(\texte^{-\eta\min\{k,n-k\}})\bigr] (t_\star\theta^{\alpha^2})^{2k}
\end{equation}
where~$k:=k(x,y)$. Denoting
\begin{equation}
\label{Eq:5.56}
\kappa(\alpha,\beta) := 4\frac{2\pi^2}{\beta\log b}\alpha^2 - \frac4{\log b} \log t_\star(\alpha,\beta) 
\end{equation}
we have $(t_\star\theta^{\alpha^2})^{2k} = d(x,y)^{-\kappa(\alpha,\beta)}$, where the factors~$\log b$ arise from the conversion $d(x,y) = (b^{1/2})^{k(x,y)}$. The inequality in \eqref{E:1.16} holds by the fact that~$t_\star>1$ for $\alpha\ne0$. The sequence $\{C_n\}_{n\ge1}$ is uniformly positive and finite for $b\theta-1$ small thanks to Lemma~\ref{lemma-5.9} and the fact that~$f_\star$ is dominated by the $(p,j):=(0,1)$ term with the rest being at least order~$\sqrt{b\theta-1}$.
\end{proofsect}

\begin{figure}[t]
\refstepcounter{obrazek}
\label{fig4}
\vglue5mm
\centerline{\includegraphics[width=4in]{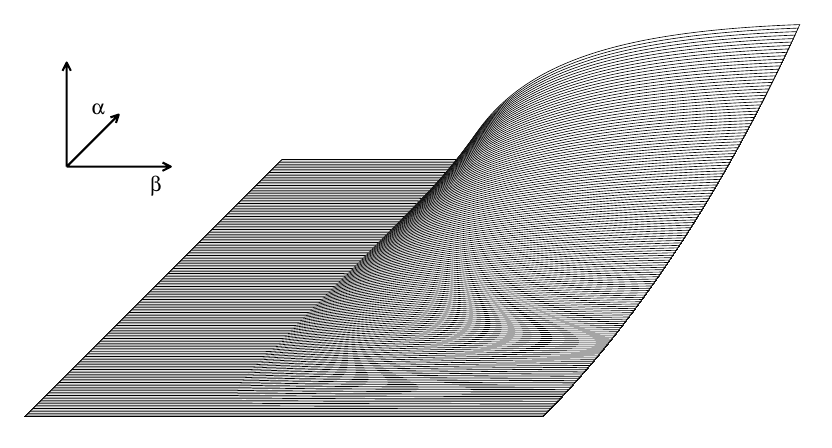}}
\vspace{.2in}
\begin{quote}
\fontsize{9}{5}\selectfont
{\sc Fig~\theobrazek.} A plot of the graph of $(\alpha,\beta)\mapsto-\log t_\star(\alpha,\beta)$ for $\beta\in[25,35]$ and~$\alpha\in[0,0.4]$ when $b = 2$. The flat region to the left corresponds to~$\beta<\beta_\cc$.
\normalsize
\end{quote}
\end{figure}

\begin{remark}
\label{rem-kappa}
It remains to show that~$\kappa(\alpha,\beta)$ obeys  \eqref{E:1.20i}. For this we will have to extract the leading-order asymptotic of~$t_\star$ in powers of~$\epsilon:=\sqrt{b\theta-1}$, for~$\epsilon$ small positive. We start by noting that
\begin{equation}
\gg_\star(0)=1- 2 (b-1)\lambda_\star(1)^2+O(\epsilon^4)
\end{equation}
while
\begin{equation}
\gg_\star(1) = (b-1)\lambda_\star(1)+O(\epsilon^3)
\end{equation}
and~$\sum_{n\ge2}\gg_\star(n) = O(\epsilon^2)$. This now gives
\begin{equation}
\Gamma^\star_1(0)=1- 2 (b-1)\lambda_\star(1)^2+O(\epsilon^4)
\end{equation}
and, for $j\ge2$,
\begin{equation}
\Gamma^\star_j(0) = (b-1)^2\bigl(\theta^{(1+2\alpha)(j-1)}+\theta^{(1-2\alpha)(j-1)}\bigr)\lambda_\star(1)^2+O(\epsilon^4)b^{-(1-2|\alpha|)j},
\end{equation}
where the implicit constant in $O(\epsilon^4)$ does not depend on~$j$.
Using this in \eqref{Eq:5.47u} while noting that~$\lambda_\star(1)=O(\sqrt{b\theta-1})$, a calculation shows
\begin{equation}
\label{E:5.118}
t_\star = 1+(b-1)\biggl[\frac{b-1}{b^{1+2\alpha}-1}+\frac{b-1}{b^{1-2\alpha}-1}- 2\biggr]\lambda_\star(1)^2+O\bigl((b\theta-1)^2\bigr).
\end{equation}
In the derivation of the asymptotic \twoeqref{E:6.38}{E:6.38a} we will show that
\begin{equation}
\lambda_\star(1) = \sqrt{\frac{2(b^3-1)}{(b-1)^2(b+1)^3}}\,\sqrt{b\theta-1}+O\bigl(|b\theta-1|^{3/2}\bigr).
\end{equation}
Inserting  this in \eqref{E:5.118} then gives
\begin{equation}
t_\star(\alpha,\beta)=1+\tau(\alpha)(b\theta-1)+O\bigl((b\theta-1)^2\bigr).
\end{equation}
Invoking \eqref{Eq:5.56} along with $b\theta-1=\frac{2\pi^2}{\beta_\cc^2}(\beta-\beta_\cc)+O((\beta-\beta_\cc)^2)$ we get \eqref{E:1.20i}. The numerical closeness of the critical and near-critical coefficients again stems from the structure of $\Gamma^\star_j(0)$ and the similarity of leading order term of~$\lambda_k(1)$ in the critical and slightly supercritical regimes.
\end{remark}

\begin{remark}
\label{rem-single-charge}
The above proofs were tailored to the asymptotic of the fractional dipole-charge correlation  but the structure also applies to its monopole counterpart, where the conclusion is in fact easier. Indeed, we get the identity
\begin{equation}
\langle \texte^{2\pi\texti\alpha\phi_x}\rangle_{n,\beta} = f_{n+1}(0)
\end{equation}
for~$f_{n+1}$ obtained by taking $k:=n+1$ in \eqref{Eq:5.2b}. (Here we used that the underlying Markov chain effectively takes value zero at the initial time so no expectation is needed.). Setting
\begin{equation}
C_n':=\begin{cases}
(t_\star\theta^{\alpha^2})^{-n}\langle \texte^{2\pi\texti\alpha\phi_x}\rangle_{n,\beta},\qquad&\text{if }\beta\ne\beta_\cc,
\\
k^{-\tilde\tau}\theta^{-\alpha^2n}\langle \texte^{2\pi\texti\alpha\phi_x}\rangle_{n,\beta},\qquad&\text{if }\beta=\beta_\cc,
\end{cases}
\end{equation}
where~$t_\star:=1$ when~$\beta<\beta_\cc$, the arguments used above show that~$C_n'$ is bounded away from zero and infinity.  Writing this using~$\kappa(\alpha,\beta)$ we get \eqref{E:1.21} as desired.
\end{remark}

\section{Supercritical iterations}
\label{sec-6}\noindent
The principal objective of this section is the proof of Theorem~\ref{thm-5}. This requires studying the flow of the iterations \eqref{E:3.8} under the conditions when these admit a ``non-trivial'' fixed point. The analysis carries a significant technical overhead that goes quite beyond what was sufficient for the subcritical and critical cases.

\newcommand{\blambda}{\bm\lambda}

\subsection{Renormalization group flow}
We start by casting the iterations in a more convenient and also more general form. 
Recall the definition of~$\theta_k$ from \eqref{Eq:3.14w} and~$\theta$ from \eqref{Eq:3.27w} and note that $\beta>\beta_\cc$ is equivalent to~$b\theta>1$. Let $\Sigma$ be the set of (doubly-infinite) positive sequences $\blambda=\{\lambda(q)\}_{q\in\Z}$ satisfying
the symmetry condition~$\lambda(-q)=\lambda(q)$ for all~$q\in\Z$ and such that $\lambda(0)=1$ and $\sup_{q\ge0}\lambda(q+1)/\lambda(q)<\infty$ hold true. For each~$\blambda\in\Sigma$ and~$q\in\Z$ set
\begin{equation}
\label{E:6.1}
G_q(\blambda) := \sum_{\begin{subarray}{c}
\ell_1,\dots,\ell_b\in\mathbb Z\\\ell_1+\dots+\ell_b=q
\end{subarray}}
\prod_{i=1}^b \lambda(\ell_i)
\end{equation}
and, for each $k\ge0$, let
\begin{equation}
F_q^{(k)}(\blambda) :=  \frac{G_q(\blambda)}{G_0(\blambda)}\theta_{k}^{q^2}.
\end{equation}
We will write~$F^{(k)}$ for the map assigning~$\blambda$ the sequence $\{F^{(k)}_q(\blambda)\}_{q\in\Z}$ and use the notation~$F$ for the corresponding map in which the~$\theta_k$'s have been replaced by~$\theta$.

In order to make the connection to the problem at hand, note that extending our earlier notation \eqref{E:3.46} to
\begin{equation}
\label{E:6.3}
\lambda_k(q) := \frac{a_k(q)}{a_k(0)},\quad q\in\Z,
\end{equation}
and denoting $\blambda_k:=\{\lambda_k(q)\}_{q\in\Z}$, the iterations \eqref{E:3.8} become
\begin{equation}
\label{E:1.5}
\blambda_{k} = F^{(k)}(\blambda_{k-1}),\quad k\ge1,
\end{equation}
where the parametrization reflects that $a_{0}(\cdot)/a_{0}(0)$ corresponds to~$\blambda_{0}$.
We are thus interested in the convergence/limit properties of the flow of compositions of functions $\{F^{(k)}\colon k\ge0\}$ evaluated on elements from~$\Sigma$. 

Our control of the iterations turns out to be  slightly  stronger when~$b$ is even. Indeed, for~$b$ even we can work with any starting~$\blambda_{0}\in\Sigma$ while for~$b$ odd we have to assume that the initial~$\blambda_{0}$ arises from the setting of the present work. We thus set $\Sigma':=\Sigma$ when~$b$ is even and let~$\Sigma'$ be the set of~$\blambda\in\Sigma$ that are the Fourier coefficients of a reflection-symmetric Borel probability measure on~$[0,1)$ when~$b$ is odd. Then we restate the key part of Theorem~\ref{thm-5} as:

\begin{theorem}
\label{thm-6.1}
Let~$b\ge2$. There exists~$\epsilon>0$ and, for each~$\beta>0$ satisfying $1<b\theta<1+\epsilon$, there exists a unique $\blambda_\star\in\Sigma'$ such that
\begin{equation}
F(\blambda_\star) = \blambda_\star.
\end{equation}
Moreover, under Assumption~\ref{ass-1}, for each~$\beta$ as above there exist~$\eta>0$ and~$C>0$ and, for each $\blambda_{0}\in\Sigma'$, there exists $k_0\ge0$ such that, for all~$n>2k_0$,
the sequence $\{\blambda_k\}_{k=0}^n$ defined from~$\blambda_{0}$ via \eqref{E:1.5} obeys 
\begin{equation}
\label{E:6.6iu}
\sum_{q\ge1}\bigl[4 b^{3}(b\theta-1) \bigr]^{\frac{q-1}2}\bigl|\lambda_k(q)-\lambda_\star(q)\bigr|\le C\biggl[\texte^{-\eta k}+\sum_{j=0}^k\texte^{-\eta(k-j)}\,\frakd_{\min\{j,n-j\}}\biggr]
\end{equation}
whenever~$\min\{k,n-k\}\ge k_0$. Here~$\{\frakd_j\}_{j\ge0}$ is the sequence from Assumption~\ref{ass-1}
\end{theorem}

While the above may appear to be a run-off-the-mill conclusion of the Banach Fixed Point Theorem, the proof is considerably more complicated. A key problem is that~$F$ is not contractive on~$\Sigma'$ when~$\beta>\beta_\cc$ due to the ``subcritical'' fixed point (corresponding to $\lambda(q)=\delta_{q,0}$) lingering on the ``boundary'' of~$\Sigma'$. This fixed point is unstable for~$\beta>\beta_\cc$ which mucks up uniform control of the iterations.

Our way to overcome this is by following the iterations until they reach a suitable subset $\Sigma_0\subseteq\Sigma'$ where contractivity can  be proved. We assume $b\theta-1$ small as, under this condition, the evolution of~$\blambda_k$ is completely controlled by $\lambda_k(1)$ and~$\lambda_k(2)$, just as we saw happens for the critical case in Lemma~\ref{lemma-3.12} and the proof of Theorem~\ref{thm-4}. Indeed, these two coordinates evolve autonomously (modulo error terms) according to \eqref{Eq:1.47b}  while the remaining ones are just ``swept along.'' 

Working near critical~$\beta$ unfortunately means that the convergence~$\blambda_k\to\blambda_\star$ is very slow; in fact, it gets the slower the closer~$b\theta$ is to~$1$. Indeed, our proof gives \eqref{E:6.6iu} with~$\eta$ proportional to~$b\theta-1$ which, as is easy to check, reflects also the true decay rate when~$\sigma_k^2=1$ for all~$k\ge0$. The inhomogeneity of~$\{\sigma_k^2\}_{k=0}^n$ causes further errors that are governed by the tails of the convergent series $\sum_{j\ge0}\frakd_j$.

\subsection{Preliminary observations}
We start by some preliminary technical estimates. As noted above, our control of the iterations is better when~$b$ is even. This is due to the availability of:

\begin{lemma}
\label{lemma-1.2}
Suppose that~$b\ge2$ is even and let~ $G_q$  be as in \eqref{E:6.1}. Then for all~$\blambda\in\Sigma$,
\begin{equation}
\label{E:2.15}
G_q(\blambda)\le G_0(\blambda),\quad q\in\Z.
\end{equation}
In particular, we have $F_q(\blambda)\le\theta^{q^2}\le1$ for all~$q\in\Z$.
\end{lemma}

\begin{proofsect}{Proof}
Assume~$b$ even and recall the definition of $\Xi_b(n)$ from \eqref{E:5.73i}. 
Our goal is to show that, for all $\blambda\in\Sigma$ and~$n\in\Z$,
\begin{equation}
\label{E:2.15a}
\sum_{\bar\ell\in\Xi_b(n)}
\prod_{i=1}^b \lambda(\ell_i)\le\sum_{\bar\ell\in\Xi_b(0)}
\prod_{i=1}^b \lambda(\ell_i),
\end{equation}
where $\bar\ell=(\ell_1,\dots,\ell_b)$. For this we note that, since~$b$ is even, $b':=b/2$  is a natural number and the distributive law yields
\begin{equation}
\label{E:2.16}
\sum_{\bar\ell\in\Xi_b(q)}
\prod_{i=1}^b \lambda(\ell_i)=
\sum_{j\in\Z}\biggl(\,\sum_{\bar\ell'\in\Xi_{b'}(j)}
\prod_{i=1}^{b'} \lambda(\ell_i')\biggr)\biggl(\,\sum_{\bar\ell''\in\Xi_{b'}(q-j)}
\prod_{i=1}^{b'} \lambda(\ell_i'')\biggr).
\end{equation}
On the other hand, the same argument and the symmetry condition $\lambda(-\ell)=\lambda(\ell)$ shows
\begin{equation}
\label{E:6.11ui}
\sum_{j\in\Z}\biggl(\,\sum_{\bar\ell'\in\Xi_{b'}(j)}
\prod_{i=1}^{b'} \lambda(\ell_i')\biggr)^2 = 
\sum_{\bar\ell\in\Xi_b(0)}
\prod_{i=1}^b \lambda(\ell_i).
\end{equation}
To get the desired claim \eqref{E:2.15a}, it suffices to invoke the Cauchy-Schwarz inequality in \eqref{E:2.16} and apply \eqref{E:6.11ui}.
\end{proofsect}

The distinction between~$b$ even and $b$ odd now enters solely through the following enhancement of Lemmas~\ref{lemma-3.5}--\ref{lemma-3.6}:

\begin{lemma}
\label{lemma-6.2}
For any $\blambda_{0}\in\Sigma'$ and all~$n\ge k\ge1$,
\begin{equation}
\label{E:6.7}
\sup_{q\ge0}\frac{\lambda_k(q+1)}{\lambda_k(q)}\le\biggl(\, \prod_{i=1}^k \max\{b\theta_i^3,1\}\biggr)\max\biggl\{1,\sup_{q\ge0}\frac{\lambda_{0}(q+1)}{\lambda_{0}(q)}\biggr\}.
\end{equation}
\end{lemma}

\begin{proofsect}{Proof}
Let $\blambda_{0}\in\Sigma'$ and, setting $a_{0}(0):=1$, let~$\{a_{0}(q)\}_{q\in\Z}$ be such that $a_{0}(q)/a_{0}(0)=\lambda_{0}(q)$ for each~$q\in\Z$.
Lemma~\ref{lemma-3.5} along with the fact that~ $(q + 1)^2 - q^2 \ge 3$  once $q\ge1$ imply
\begin{equation}
\label{E:6.13}
\sup_{q\ge1}\frac{a_{k}(q+1)}{a_{k}(q)}\le b\theta_{k}^{3}\sup_{q\ge0} \frac{a_{k-1}(q+1)}{a_{k-1}(q)} .
\end{equation}
Now note that, for~$b$ even, Lemma~\ref{lemma-1.2} shows that $a_k(n)\le a_k(0)$ for all~$n\in\Z$ and~$k\ge0$ while for~$b$ odd this holds by the fact that the Fourier coefficients of a probability measure on~$[0,1)$ are bounded by~$1$. Denoting, as before, the supremum on the left of \eqref{E:6.7} as~$c_k$, we are thus led to the inequality
\begin{equation}
c_{k}\le\max\{1,b\theta_{k}^3 c_{k-1}\}.
\end{equation}
To iterate this, set $\tilde c_k:=\max\{1,c_k\}$ are note that the above gives $\tilde c_{k}\le \max\{1,b\theta_{k}^3\}\tilde c_{k-1}$. This now readily implies \eqref{E:6.7}.
\end{proofsect}

As a consequence we get:

\begin{corollary}
\label{cor-6.4}
Suppose~$b\theta^3<1$. Then, under Assumption~\ref{ass-1}, for all~$\blambda_{0}\in\Sigma'$,
\begin{equation}
\label{E:6.11i}
A_1:=\biggl(\,\sup_{n\ge1}\prod_{i=0}^n\max\{b\theta_i^3,1\}\biggr)\max\biggl\{1,\sup_{q\ge0}\frac{\lambda_{0}(q+1)}{\lambda_{0}(q)}\biggr\}<\infty
\end{equation}
and, for each~$n\ge1$, the iterations~$\{\blambda_k\}_{k=0}^n$ generated from $\blambda_{0}$ via \eqref{E:1.5} obey
\begin{equation}
\label{E:6.8}
\sup_{n\ge k\ge0}\sup_{q\ge0}\frac{\lambda_k(q+1)}{\lambda_k(q)}\le A_1.
\end{equation}
In particular, $F^{(k)}(\Sigma')\subseteq\Sigma'$ for each~$k\ge0$.
\end{corollary}

\begin{proofsect}{Proof}
The condition~$b\theta^3<1$ along with Assumption~\ref{ass-1} imply that the product is uniformly bounded. Invoking \eqref{E:6.7}, we get \eqref{E:6.8} with~$A_1$ as in \eqref{E:6.11i}. Using \eqref{E:6.1} we check that~$F^{(k)}$ preserves the symmetry and, for~$b$ odd, a rewrite using the iteration of the potentials shows that $F^{(k)}(\blambda)$ is a Fourier transform of a probability measure on~$[0,1)$ if~$\blambda$ is. It follows that $F^{(k)}(\Sigma')\subseteq\Sigma'$ as desired. 
\end{proofsect}

\begin{remark}
To phrase the above in the vernacular of the renormalization group theory, under the condition $b\theta^3<1$, the estimate \eqref{E:6.13} says that $\lambda(q)$ with~$|q|\ge2$ are irrelevant (i.e., contracting) ``directions'' of  the renormalization group ``flow.'' Due to the normalization~$\lambda(0)=1$, the only possibly relevant (i.e., expanding) ``direction'' in this regime is thus~$\lambda(1)=\lambda(-1)$. The punchline of Corollary~\ref{cor-6.4} is that the expansion in this coordinate is still clamped down by a uniform bound.
\end{remark}

Before we move to the consequences of above observations, let us record the following general bound that shows up repeatedly in the sequel:

\subsection{Near-critical bounds}
We will now improve the above crude estimate on $\sup_{q\ge0}\lambda_k(q+1)/\lambda_k(q)$ to a bound that is small in the ``near-critical'' regime, i.e., for~$b\theta-1$ small positive.

\begin{lemma}
\label{lemma-6.7a}
Given~$\wt\alpha\in(0,1)$, suppose~$\beta>0$ is such that $b\theta^2<b-1$ and $(b-1+\wt\alpha)\theta<1$. Under Assumption~\ref{ass-1}, for each~$\blambda_{0}\in\Sigma'$ there exists~$k_1\ge0$ such that
\begin{equation}
\label{E:6.26b}
\sup_{q\ge0}\frac{\lambda_k(q+1)}{\lambda_k(q)}\le A_2\sqrt{b\theta-1}
\end{equation}
holds with
\begin{equation}
\label{E:6.19a}
A_2:= \frac1{\wt\alpha}\sqrt{\frac1{1-(b-1)\theta}}
\end{equation}
for all~$k\ge0$ satisfying $\min\{k,n-k\}\ge k_1$.
\end{lemma}

\begin{proofsect}{Proof} 
Let $c_k$ denote the supremum on the left-hand side of \eqref{E:6.26b}. The assumed bounds on~$\theta$ allow us to pick~$\delta>0$ so small that $b\theta^{2+2\delta}\le b-1$ and $(b-1+\wt\alpha \theta^{-3\delta})\theta\le \theta^{3\delta}$. Thanks to Assumption~\ref{ass-1} there exists~$k_1'\ge1$ with $2\sum_{j\ge k_1'}\frakd_j<\delta$ and~$\min\{k,n-k\}> k_1'$ implies
\begin{equation}
\label{E:6.19w}
b\theta_{k+1}^3\le (b-1)\theta_k  \quad\text{and}\quad
(b-1+\wt\alpha \theta^{-3\delta})\theta_{k+1}\le \theta^{2\delta}
\end{equation}
uniformly in the choice of~$\{\sigma_k^2\}_{k=0}^n$. 

The first inequality in \eqref{E:6.19w} enables the iterative bound \eqref{E:3.37}, where~$\alpha_k$ is defined in~\eqref{E:3.46iu}. 
Rewriting \eqref{E:3.37} as $c_{k+1}\le (b-1+\alpha_k)\theta_{k+1}c_k$, the second inequality in \eqref{E:6.19w} shows that, for all~$k$ with $\min\{k,n-k\}> k_1'$,
\begin{equation}
\label{E:6.20w}
\alpha_k\le\wt\alpha\theta^{-3\delta}\quad\Rightarrow\quad
c_{k+1}\le \theta^{2\delta}c_k.
\end{equation}
In the opposite regime of~$\alpha_k$, we introduce a variant of definition \eqref{E:3.48w} as
\begin{equation}
\tilde c_k:=  c_k  \exp\biggl\{\log(1/\theta)  \sum_{j=k_1'}^k  (\sigma_j^2-1)\biggr\}
\end{equation}
and observe that $2\sum_{j\ge k_1'}\frakd_j<\delta$ gives $c_k\ge\tilde c_k\theta^{\delta}$ for all~$k$ with $\min\{k,n-k\}> k_1'$ and thus
\begin{equation}
\alpha_k\ge\wt\alpha \theta^{-3\delta}\quad\Rightarrow\quad \tilde c_{k+1}\le h(\tilde c_k),
\end{equation}
where
\begin{equation}
h(u):=\frac{\theta u}{1+\wt\alpha^2 \theta^{-6\delta}  u^2 }+(b-1)\theta u.
\end{equation}
A calculation shows that~$h$ admits a unique positive fixed point at
\begin{equation}
\label{E:6.21b}
u_\star:=  \frac{\theta^{3\delta}}{\wt\alpha}\, 
\sqrt{\frac{1}{1-(b-1)\theta}}\,\sqrt{b\theta-1}
\end{equation}
with iterations of~$h$ started above~$u_\star$ decreasing geometrically fast to~$u_\star$ and those started below~$u_\star$  never exceeding~$u_\star$.  Since~$c_{k+1}\le \theta^{2\delta}c_k$ in \eqref{E:6.20w} gives $\tilde c_{k+1}\le \theta^{\delta}\tilde c_k$, whose iterations also decrease geometrically, we conclude the existence of~$k_1> k_1'$ such that $\min\{k,n-k\}\ge k_1$ implies $\tilde c_k\le u_\star\theta^{-\delta}$. Since~$c_k\le \tilde c_k\theta^{-\delta}$, this gives the claim.
\end{proofsect}

We now use the above to extract the asymptotic values of $\lambda_k(1)$ and~$\lambda_k(2)$ for~$b\theta-1$ small. This also shows that the supremum in \eqref{E:6.26b} is dominated by the $q=0$ term.

\begin{lemma}
\label{lemma-6.8}
For each~$\delta,\delta'\in(0,1)$ there exists~$\epsilon>0$ such that for all~$\blambda_{0}\in\Sigma'$ and all~$\beta>0$ with $1<b\theta<1+\epsilon$ and  there exists~$k_2\ge0$ for which
\begin{equation}
\label{E:6.38}
\lambda_k(1)\ge\biggl(\frac{(b-1)^2}2\frac{(b+1)^3}{b^3-1}+\delta\biggr)^{-1/2}\sqrt{b\theta-1}
\end{equation}
\begin{equation}
\label{E:6.38a}
\lambda_k(1)\le\biggl(\frac{(b-1)^2}2\frac{(b+1)^3}{b^3-1}-\delta\biggr)^{-1/2}\sqrt{b\theta-1}
\end{equation}
and
\begin{equation}
\label{E:6.26ua}
\biggl|\frac{\lambda_k(2)}{b\theta-1}- \frac1{(b-1)(b+1)^3}\Bigr|\le\delta'
\end{equation}
hold when $\min\{k,n-k\}\ge k_2$.
\end{lemma}

\begin{proofsect}{Proof}
As before, we will use the shorthands $\lambda_k:=\lambda_k(1)$ and $\rho_k:=\lambda_k(2)$  and, committing major abuse of notation,  abbreviate $\epsilon:=\sqrt{b\theta-1}$. Our first goal is to show that $\lambda_k$ will eventually be at least order~$\epsilon$. For this we invoke the inequality \eqref{Eq:3.55} from the proof of Lemma~\ref{lemma-3.9} which reads
\begin{equation}
\lambda_{k}\ge  \frac{b\theta_{k}}{1+\alpha' c_{k-1}^3}  \frac{\lambda_{k-1}}{1+b(b-1)\lambda_{k-1}^2}
\end{equation}
for some constant~$\alpha'>0$, whenever~$c_{k-1}\le1/2$. The latter is enabled by assuming $\min\{k,n-k\}> k_1$ and $A_2\epsilon\le1/2$, for~$k_1$ and~$A_2$ as in Lemma~\ref{lemma-6.7a} with, say, $\wt\alpha:=1/2$. Plugging in the bound $c_{k-1}\le A_2\epsilon$ while noting $b\theta_k = b\theta^{\sigma_k^2} = (1+\epsilon^2)\theta^{\sigma_k^2-1}$ we then get
\begin{equation}
\label{E:6.28}
 \frac{b\theta_{k}}{1+\alpha' c_{k-1}^3} \ge \frac{1+\epsilon^2}{1+\alpha' A_2^3\epsilon^3}\theta^{|1-\sigma_k^2|}.
\end{equation}
We now take $k_1'\ge k_1$ so large that $\frakd_j\le\epsilon^3$ once~$j\ge k_1'$ and~$\epsilon$ so small that right-hand side of \eqref{E:6.28} is at least~$1+\epsilon^2/2$. It follows that 
\begin{equation}
\lambda_k\ge\Bigl(1+\frac{\epsilon^2}2\Bigr)\frac{\lambda_{k-1}}{1+b(b-1)\lambda_{k-1}^2}
\end{equation}
once $\min\{k,n-k\}\ge k_1'$. Interpreting the right-hand side as $h(\lambda_{k-1})$ for a suitable function $h\colon\R_+\to\R_+$, we now check from the properties of~$h$ that iterations of~$h$ started from any positive value are attracted to its unique positive fixed point at $[2b(b-1)]^{-1/2}\epsilon$. Hence, there exists~$k_1''\ge k_1'$ such that
\begin{equation}
\label{E:6.30ii}
\lambda_k\ge\frac1{2\sqrt{b(b-1)}}\epsilon
\end{equation}
once~$\min\{k,n-k\}\ge k_1''$.

Next we observe that, with~$c_k$ bounded by a constant times~$\lambda_k$, whenever~$k$ obeys $\min\{k,n-k\}\ge k_1''$ the identities \eqref{Eq:1.47b} from Lemma~\ref{lemma-3.12} are still in force. The calculation leading up to \eqref{E:3.58} still applies. As~$\delta_k$ there is order~$\epsilon$ which by \eqref{E:6.30ii} is order~$\lambda_k$, instead of \eqref{E:3.58a} we then get
\begin{equation}
\label{E:6.31}
\frac{\rho_{k}}{\lambda_{k}^2}=\frac12\frac{b-1}{b^3-1}+ t_k''\lambda_k
\end{equation}
for some bounded sequence $\{t_k''\}_{k=0}^n$. Plugging this in the first line of \eqref{Eq:1.47b} yields
\begin{equation}
\lambda_{k+1}=b\theta_{k+1}\frac{\lambda_k+ \frac{b-1}2(b-2+\frac{b-1}{b^3-1})\lambda_k^3+r_k' \lambda_k^4}{1+b(b-1)\lambda_k^2+s_k\lambda_k^3}
\end{equation}
where~$r_k':=r_k+(b-1)t_k''$. 

In order to analyze this further, note that for $\min\{k,n-k\}$ so large that~$\frakd_k\le\epsilon$ the above leads to the inequalities
\begin{equation}
\label{E:6.33ue}
b\theta\frac{\lambda_k}{1+(A_\star+B\epsilon)\lambda_k^2}\le\lambda_{k+1}\le b\theta\frac{\lambda_k}{1+(A_\star-B\epsilon)\lambda_k^2},
\end{equation}
where
\begin{equation}
\label{E:6.34i}
\begin{aligned}
A_\star:=b(b-1)-\frac{b-1}2\Bigl(b-2+\frac{b-1}{b^3-1}\Bigr) = \frac{b-1}2\Bigl[b+2-\frac1{b^2+b+1}\Bigr]
\\=\frac{b-1}2\frac{(b+1)^3}{b^2+b+1} = \frac{(b-1)^2}2\frac{(b+1)^3}{b^3-1}.
\end{aligned}
\end{equation}
and~$B$ is a positive constant derived from the bounds on the sequences~$r_k'$ and~$s_k$. Noting iterations of $h(u)=b\theta\frac{u}{1+Au^2}$ are attracted to its unique fixed point at $A^{-1/2}\sqrt{b\theta-1}$, following the iterations \eqref{E:6.33ue} we then get that, after a finite number of steps, we have
\begin{equation}
(A+2B\epsilon)^{-1/2}\epsilon\le\lambda_k(1)\le (A-2B\epsilon)^{-1/2}\epsilon
\end{equation}
For~$\epsilon$ such that $2B\epsilon<\delta$, this gives \twoeqref{E:6.38}{E:6.38a}. With the help \eqref{E:6.31} and~$\lambda_k\le A_2\epsilon$ we then get \eqref{E:6.26ua} as well.
\end{proofsect}

\subsection{Contractive region}
We now proceed to define a subdomain of~$\Sigma'$ on which we later prove uniform contractivity of the map~$F$. This subdomain will depend on~$\beta>0$, which we assume is such that~$b\theta>1$, and numbers $\delta,\delta'\in(0,1)$ and $A>0$ as
\begin{equation}
\label{E:6.44a}
\begin{aligned}
\Sigma_0:=\Biggl\{\blambda\in\,&\Sigma'\colon \lambda(1)\ge\Bigl(\frac{(b-1)^2}2\frac{(b+1)^3}{b^3-1}+\delta\Bigr)^{-1/2}\sqrt{b\theta-1}
\\&\wedge
\lambda(1)\le\Bigl(\frac{(b-1)^2}2\frac{(b+1)^3}{b^3-1}-\delta\Bigr)^{-1/2}\sqrt{b\theta-1}
\\&\wedge
\,\sup_{q\ge0}\frac{\lambda(q+1)}{\lambda(q)}\le A\sqrt{b\theta-1}\,\,\wedge\,\,\Bigl|\frac{\lambda(2)}{b\theta-1}- \frac1{(b-1)(b+1)^3}\Bigr|\le\delta'\Biggr\}.
\end{aligned}
\end{equation}
Assuming that Assumption~\ref{ass-1} holds, we now summarize the previous observations in:

\begin{lemma}
\label{lemma-6.12}
For all $A>\sqrt b$ and~$\delta,\delta'\in(0,1)$, there exists~$\epsilon>0$ such that, for $\Sigma_0$ defined by~$\delta$,~$\delta'$ and~$A$ as above, the following is true for all~$\beta$ satisfying~$1<b\theta<1+\epsilon$: For all $\blambda_{0}\in\Sigma'$ there exists~$k_3\ge0$ such that $\blambda_k\in\Sigma_0$ holds whenever $\min\{k,n-k\}\ge k_3$
for the iterations  $\{\blambda_k\}_{k=0}^n$  defined from~$\blambda_{0}$ via~\eqref{E:1.5}.
\end{lemma}

\begin{proofsect}{Proof}
This follows from Lemmas~\ref{lemma-6.7a} and~\ref{lemma-6.8} along with the fact that, for $\wt\alpha\in(\sqrt b/A,1)$, the quantity~$A_2$ in \eqref{E:6.19a} tends to $\wt\alpha\sqrt b<A$ as~$b\theta$ decreases to~$1$.
\end{proofsect}

We will henceforth focus on the evolution driven by~$F$, i.e., for~$\theta_k=\theta$ for all~$k$. Here we need to check that~$F$ maps~$\Sigma_0$ into itself.

\begin{lemma}
\label{lemma-6.13}
For each~$A>\sqrt{\frac2{b-1}}$, each $\delta'\in(0,1)$ and each~$\delta\in(0,\frac12b^{-4}\delta')$ there exists $\epsilon>0$ such that
\begin{equation}
F(\Sigma_0)\subseteq\Sigma_0
\end{equation}
holds for all~$\beta>0$ with~$1<b\theta<1+\epsilon$.
\end{lemma}

\begin{proofsect}{Proof}
Fix~$A$,~$\delta'$ and~$\delta$ as above and, abusing notation again, abbreviate $\epsilon:=\sqrt{b\theta-1}$. Note that $A_\star$ from \eqref{E:6.34i} obeys~$A_\star>1$ so the expression on the second line in \eqref{E:6.44a} is meaningful. Pick~$\blambda\in\Sigma_0$ and note that, since~$F(\blambda)\in\Sigma'$ by Corollary~\ref{cor-6.4}, we only need to verify that $F(\blambda)$ obeys the conditions in \eqref{E:6.44a}. 

For the first two conditions in \eqref{E:6.44a}, we repeat the calculations underlying the proof of Lemma~\ref{lemma-6.8} to get
\begin{equation}
F_1(\blambda)= b\theta\frac{\lambda(1)+[\binom{b-1}2+\frac12\frac{(b-1)^2}{b^3-1}+\eta\epsilon^2]\lambda(1)^3}{1+[b(b-1)+ \eta'\epsilon^2]\lambda(1)^2}
\end{equation}
for $\eta$ and~$\eta'$ depending only on $\delta'$ and~$A$. Using the same calculation as that leading to \eqref{E:6.33ue}, this gives 
\begin{equation}
b\theta\frac{\lambda(1)}{ 1+(A_\star+B\epsilon)\lambda(1)^2}\le F_1(\blambda)\le b\theta\frac{\lambda(1)}{ 1+(A_\star -B\epsilon)\lambda(1)^2}
\end{equation}
for a constant~$B\ge0$ derived from~$\eta$ and~$\eta'$. Now check that, for $(A_\star+B\epsilon)\lambda(1)^2\le1$ (which holds for~$\epsilon$ small by~$\blambda\in\Sigma_0$) both sides of this inequality are non-decreasing in~$\lambda(1)$. Noting also that $b\theta=1+\epsilon^2$, it follows that, for~$\epsilon$ so small that~$B\epsilon\le\lambda$, the expression on the left then  preserves the inequality~$\lambda(1)\ge (A_\star +\delta)^{-1}\epsilon$  while that on the right preserves the inequality $\lambda(1)\le (A_\star-\delta)^{-1}\epsilon$. We conclude that~$F_1(\blambda)$ obeys the first two lines in~\eqref{E:6.44a}.

For the third condition in \eqref{E:6.44a} we first note that Lemma~\ref{lemma-3.7} along with the third condition for~$\blambda$ give
\begin{equation}
F_1(\blambda)\le\theta\frac{\lambda(1)}{1+\binom b2\lambda(1)^2}+(b-1)\theta\sup_{n\ge0}\frac{\lambda(n+1)}{\lambda(n)}
\\
\le\biggl(\frac{\theta}{1+\binom b2A^2\epsilon^2}+(b-1)\theta\biggr)A\epsilon,
\end{equation}
where we assumed~$\epsilon$ small enough so that the first term is non-decreasing in~$\lambda(1)\le A\epsilon$.
As a calculation shows, under $(b-1)\theta<1$, the prefactor of~$A\epsilon$ is less than one for~$\epsilon$ sufficiently small if~$A>\sqrt{\frac2{b-1}}$. Since Lemma~\ref{lemma-3.5} gives
\begin{equation}
\label{E:6.56u}
\frac{F_{q+1}(\blambda)}{F_q(\blambda)}\le b\theta^{1+2q}\sup_{\ell\ge0}\frac{\lambda(\ell+1)}{\lambda(\ell)}\le b\theta^3 A\epsilon,\quad q\ge1,
\end{equation}
under~$b\theta^3\le 1$, the third condition in \eqref{E:6.44a} thus applies to~$F(\blambda)$.

Finally, for the last condition in \eqref{E:6.44a} abbreviate $C:=[(b-1)(b+1)^3]^{-1}$. We now proceed as in the derivation of \eqref{E:3.66ii} to get
\begin{equation}
F_2(\blambda)-b\theta^4\Bigl[\lambda(2)-\frac{b-1}2\lambda(1)^2\Bigr]=O(\epsilon^4).
\end{equation}
Noting that~$b\theta^4 = b^{-3}+O(\epsilon^2)$, this implies
\begin{equation}
\bigl|F_2(\blambda)-C\epsilon^2\bigr|\le b^{-3}\bigl|\lambda(2)-C\epsilon^2|+b^{-3}\biggl|\frac{b-1}2\lambda(1)^2-(b^3-1)C\epsilon^2\biggr|+A''\epsilon^4,
\end{equation}
where~$A''$ is a constant that depends only on~$A$. Invoking the conditions from~$\Sigma_0$, a  calculation shows
\begin{equation}
\bigl|\epsilon^{-2}F_2(\blambda)-C\bigr|\le b^{-3}\delta'+\biggl[\frac{(b^3-1)(b-1)}{(b+1)^2}\biggr]^2\delta+A''\epsilon^2.
\end{equation}
Under $\delta\le\frac12 b^{-4}\delta'$, the second term on the right is at most~$\delta'/2$. For~$\epsilon$ sufficiently small, the whole expression is thus at most~$\delta'$ and so~$F(\blambda)$ thus obeys also the last condition in the definition of~$\Sigma_0$. We have shown~$F(\blambda)\in\Sigma_0$ as required.
\end{proofsect}

\subsection{Contractivity}
Having identified $\Sigma_0$ and shown that~$F$ maps it into itself, we now prove that~$F$ is actually contractive on it. Note that, for $\sup_{q\in\Z}\lambda(q+1)/\lambda(q)<1$ which on~$\Sigma_0$ holds whenever~$A\sqrt{b\theta-1}<1$, each component of~$F(\blambda)$ is the ratio of two positive convergent sums with the denominator at least~$1$. It follows that ~$F$ is continuously differentiable in each variable once~$b\theta-1$ is sufficiently small. A natural way to prove contractivity is thus to estimate the derivatives of~$F$ in a suitable norm. However, this will only be useful after we show:

\begin{lemma}
\label{lemma-6.13a}
$\Sigma_0$ is a convex set.
\end{lemma}

\begin{proofsect}{Proof}
The first, second, and the last condition  in the definition of~$\Sigma_0$ are clearly preserved by convex combinations. For the  third  condition we note that
\begin{equation}
t\mapsto\frac{tA+(1-t)A'}{tB+(1-t)B'}
\end{equation}
is, for any~$A,A',B,B'>0$, monotone on~$[0,1]$. If the ratio is less than a constant at~$t=0$ and~$t=1$, it is less than that constant for all $t\in[0,1]$.
\end{proofsect}

We will for brevity write $\partial_\ell F_q$ to denote the partial derivative of the $q$-th component of~$F$ with respect to~$\lambda(\ell)$. We start with estimates on these:
\begin{lemma}
\label{lemma-6.14}
Let~$\beta>\beta_\cc$. For each~$A>0$ there exists~$\delta_0>0$ and $\epsilon_0>0$ such that,  for~$\Sigma_0$ defined by~$A$ and any $\delta,\delta'\in(0,\delta_0]$, 
\begin{equation}
\label{E:6.60}
\bigl|\partial_\ell F_q(\blambda)\bigr|
\le\begin{cases}
b(b-1)\theta^{q^2}\lambda(1)+O(\epsilon^2)
,\qquad&\text{if }|q-\ell|=1,
\\
b\theta^{q^2}+O(\epsilon^2),\qquad&\text{if }q=\ell\ge2,
\\
(b\theta^{q^2}+O(\epsilon^2))(bA\epsilon)^{|q-\ell|},\qquad& \text{if } |q-\ell| \ge 1,
\\
1- \frac12 b(b-1)\lambda(1)^2,\qquad&\text{if }q=\ell=1,
\end{cases}
\end{equation}
holds for all $\blambda\in\Sigma_0$ and all $q,\ell\ge1$ provided that $\epsilon:=\sqrt{b\theta-1}\in(0,\epsilon_0)$. Here the implicit constants in~$O(\epsilon^2)$ terms do not depend on~$q$ and~$\ell$.
\end{lemma}

\begin{proofsect}{Proof}
Fix~$q\ge1$ and~$\ell\ge1$. We start with some general considerations. Denote by~$\wt G_q$ the quantity~$G_q$ with~$b$ replaced by~$b-1$. The symmetry $\lambda(-\ell)=\lambda(\ell)$ then gives
\begin{equation}
\partial_\ell G_q(\blambda) = b\wt G_{q-\ell}(\blambda)+b\wt G_{q+\ell}(\blambda)
\end{equation}
and so, by the quotient rule,
\begin{equation}
\label{E:6.62}
\partial_\ell F_q(\blambda) = \Bigl[\bigl(\wt G_{q-\ell}(\blambda)+\wt G_{q+\ell}(\blambda)\bigr)G_0(\blambda)^{-1}-2\wt G_\ell(\blambda) G_q(\blambda)G_0(\blambda)^{-2}\Bigr] b\theta^{q^2}.
\end{equation}
Note that one term in the square bracket is positive and the other is negative. It thus suffices to estimate each of them separately.

 Moving on to actual estimates, we begin with  the third line in \eqref{E:6.60}. Here we first observe that the argument used in the proof of \cite[Lemma~4.2]{BHu} (which is our Lemma~\ref{lemma-3.5}) combined with the third condition in the definition of~$\Sigma_0$ gives 
\begin{equation}
\frac{G_{q+1}(\blambda)}{G_q(\blambda)}\le b\sup_{\ell\ge1}\frac{\lambda(\ell+1)}{\lambda(\ell)}\le bA\epsilon
\end{equation}
and similarly for the ratio $\wt G_{q+1}(\blambda)/\wt G_q(\blambda)$. This along with $G_0(\blambda)\ge1$ shows that the first term in the square bracket in \eqref{E:6.62} is at most $(bA\epsilon)^{|q-\ell|}+(bA\epsilon)^{q+\ell}$ while the second is most $2(bA\epsilon)^{q+\ell}$.  As $(bA\epsilon)^{q+\ell}=(bA\epsilon)^{|q-\ell|} (bA\epsilon)^{2\min\{\ell,q\}}$ which assuming~$bA\theta\le1$ is $(bA\epsilon)^{|q-\ell|}O(\epsilon^2)$, we get the third line in \eqref{E:6.60} with $O(\epsilon^2)$-term uniform in~$q$ and~$\ell$. 

The first and second lines in \eqref{E:6.60} require explicit treatment of the leading-order term contributing to~$\wt G_{q-\ell}(\blambda)$. This is easy for $q=\ell$ where we only need $\wt G_0(\blambda)\le 1+O(\epsilon^2)$, which is proved from \eqref{E:3.52}. This, along with the aforementioned estimates, bounds the square bracket in \eqref{E:6.62} by the maximum of $1+O(\epsilon^2)$ and~$4(bA\epsilon)^2$. For the first line (i.e., $|q-\ell|=1$) we in turn need
\begin{equation}
\wt G_1(\blambda)\le \bigl(b-1+O(\epsilon^2)\bigr)\lambda(1),
\end{equation}
which is proved similarly as \eqref{E:3.53}.
This dominates the quantity in the square bracket in \eqref{E:6.62} by~$(b-1)\lambda(1)+O(\epsilon^2)$, thus showing that $|\partial_\ell F_q(\blambda)|\le b(b-1)\theta^{q^2}\lambda(1)+O(\epsilon^2)$. The $O(\epsilon^2)$-term is independent of~$q$ and~$\ell$ in both cases. 

Unlike the previous cases, both terms in the square bracket in \eqref{E:6.62} will contribute to the alternative $\ell=q=1$ in \eqref{E:6.60}. This results in potential cancellations that force us to extract terms of order up to~$\epsilon^2$ explicitly. Since the first term in the square bracket in \eqref{E:6.62} is of order unity while the second term is of order~$\epsilon^2$, it suffices to derive an upper bound on~$\partial_1 F_1(\blambda)$ without the absolute value.  Here bounding the first term requires the upper bounds
\begin{equation}
\wt G_0(\blambda)\le 1+ (b-1)(b-2)\lambda(1)^2+O(\epsilon^4)
\end{equation}
and
\begin{equation}
\begin{aligned}
\wt G_2(\blambda) &\le (b-1)\lambda(2)+\binom{b-1}2\lambda(1)^2+O(\epsilon^4)
\\
&\le\biggl(\frac{b-1}2\Bigl[b-2+\frac{b-1}{b^3-1}\Bigr]+\delta''\biggr)\lambda(1)^2+O(\epsilon^4),
\end{aligned}
\end{equation}
where $\delta''$ is a quantity of order of~$\delta'+\delta$. Since we are tracking terms up to order~$\epsilon^2$, bounding~$G_0(\blambda)\ge1$ is not sufficient; instead we need
\begin{equation}
G_0(\blambda)\ge 1+b(b-1)\lambda(1)^2.
\end{equation}
In the second term we in turn need the lower bounds 
\begin{equation}
\begin{aligned}
G_1(\blambda)&\ge b\lambda(1)\\
\wt G_1(\blambda)&\ge (b-1)\lambda(1)
\end{aligned}
\end{equation}
along with the upper bound
\begin{equation}
G_0(\blambda)\le 1+O(\epsilon^2).
\end{equation}
Putting these together bounds $\partial_1F_1(\blambda)$ by~$b\theta$ times
\begin{equation}
\begin{aligned}
&\frac{1+[(b-1)(b-2)+\frac{b-1}2[b-2+\frac{b-1}{b^3-1}]+\delta'']\lambda(1)^2+O(\epsilon^4)}{  1+b(b-1)\lambda(1)^2  }-\frac{b(b-1)\lambda(1)^2}{1+O(\epsilon^2)} 
\\
&\,\,= 1 +\Biggl((b-1)(b-2)+\frac{b-1}2\Bigl[b-2+\frac{b-1}{b^3-1}\Bigr]+\delta''-  2b(b-1) \Biggr)\lambda(1)^2+O(\epsilon^4)
\\
&\,\, = 1-\biggl(\frac{b-1}2\Bigl( b+6  -\frac{b-1}{b^3-1}\Bigr)-\delta''\biggr)\lambda(1)^2+O(\epsilon^4).
\end{aligned}
\end{equation}
Invoking the bound
\begin{equation}
b\theta\le 1+\biggl[\frac{(b-1)^2}2\frac{(b+1)^2}{b^3-1}+\delta\biggr]\lambda(1)^2
\end{equation}
implied by $\blambda\in\Sigma_0$ we thus get
\begin{equation}
\begin{aligned}
\partial_1F_1(\blambda)&\le 1 -\Biggl(\frac{b-1}2\Bigl( b+6  -\frac{b-1}{b^3-1}\Bigr)-\delta''-\frac{(b-1)^2}2\frac{(b+1)^2}{b^3-1}-\delta\Biggr)\lambda(1)^2+O(\epsilon^4)\\
&= 1  -  \frac{b-1}2\biggl(b+6-\frac{b-1}{b^3-1}\bigl((b+1)^2+1\bigr)-\delta''-\delta\biggr)\lambda(1)^2+O(\epsilon^4)
\end{aligned}
\end{equation}
We now check that the quantity in the last large parentheses is at least $b+4-\delta'-\delta$. Writing $O(\epsilon^4)=\lambda(1)^2O(\epsilon^2)$ using the first line in \eqref{E:6.44a}, for~$\delta,\delta'\in(0,1)$ the right-hand side is thus at least~$1-\frac{b-1}2b\lambda(1)^2$ once~$\epsilon$ is sufficiently small, uniformly in~$\blambda\in\Sigma_0$.
\end{proofsect}

We now use these to prove:

\begin{lemma}
\label{lemma-6.16}
Let $\Sigma_0$ be defined using $A>\sqrt{\frac2{b-1}}$ and~$\delta\in(0,\delta_0)$, for~$\delta_0$ as in Lemma~\ref{lemma-6.14}. Given~$t>0$ and $\beta>0$ with~$0<b\theta-1<(t bA^2)^{-1}$, we have
\begin{equation}
\label{E:6.56}
\varrho(\blambda,\blambda'):=\sum_{q\ge1}(t bA)^{q-1}(b\theta-1)^{\frac{q-1}2}\bigl|\lambda(q)-\lambda'(q)\bigr|< \infty
\end{equation}
for all~$\blambda,\blambda'\in\Sigma_0$. Moreover, if~$t$ obeys
\begin{equation}
\label{E:6.78}
\frac1{\sqrt{2}}\,\frac{\sqrt{1-b^{-3}}}{(1-b^{-2})\sqrt{1+b^{-1}}}\, b^2>t A
\end{equation}
and 
\begin{equation}
\label{E:6.79}
\sup_{\ell\ge2}\biggl(b^{1-\ell^2}+\sum_{q=1}^{\ell-1} t^{q-\ell}b^{1-q^2}\biggr)<1,
\end{equation}
then there exist~$\delta_1>0$, $\epsilon_1>0$ and~$\eta>0$ such that
\begin{equation}
\label{E:6.74}
\varrho\bigl(F(\blambda),F(\blambda')\bigr)\le \bigl[1-\eta(b\theta-1)\bigr]\varrho(\blambda,\blambda')
\end{equation}
holds for all~$\blambda,\blambda'\in\Sigma_0$ provided $\delta<\delta_1$ and $b\theta-1\le\epsilon_1$.
\end{lemma}

\begin{proofsect}{Proof}
Let us again abbreviate $\epsilon:=\sqrt{b\theta-1}$ and let~$\blambda,\blambda'\in\Sigma_0$. The conditions defining~$\Sigma_0$ then imply $|\lambda(q)-\lambda'(q)|\le 2(A\epsilon)^{|q|}$ and so the series in \eqref{E:6.56} converges whenever $tbA^2\epsilon^2<1$. Next recall that, by Lemma~\ref{lemma-6.13a}, the convex combination~$\blambda_u:=(1-u)\blambda+u\blambda'$ lies in~$\Sigma_0$ for all~$u\in[0,1]$. Moreover, elementary calculus shows
\begin{equation}
\label{E:6.75}
\begin{aligned}
\bigl|F_q(\blambda')-F_q(\blambda)\bigr|
&=\biggl|\int_0^1 \sum_{\ell\ge1}\partial_\ell F_q(\blambda_u)\bigl(\lambda'(\ell)-\lambda(\ell)\bigr)\textd u\Bigr|
\\
&\le \int_0^1\biggl(\,\sum_{\ell\ge1}\bigl|\partial_\ell F_q(\blambda_u)\bigr|\bigl|\lambda'(\ell)-\lambda(\ell)\bigr|\biggr)\textd u
\end{aligned}
\end{equation}
where the second line follows by the triangle inequality.

Multiplying \eqref{E:6.75} by~$(tbA\epsilon)^{q-1}$ and summing over~$q\ge1$ we find out that, in order to prove \eqref{E:6.74}, it suffices to show that for all~$\ell\ge1$ and all~$\blambda\in\Sigma_0$,
\begin{equation}
\label{E:6.76}
\sum_{q\ge1}(tbA\epsilon)^{q-1}\bigl|\partial_\ell F_q(\blambda)\bigr|\le (1-\eta\epsilon^2)(tbA\epsilon)^{\ell-1}.
\end{equation}
Starting first with the cases~$\ell\ge2$, here we plug in the second and third line in \eqref{E:6.60} with the result 
\begin{equation}
\begin{aligned}
\sum_{q\ge1}&(tbA\epsilon)^{q-1}\bigl|\partial_\ell F_q(\blambda)\bigr|
\\
&\le \biggl(b\theta^{\ell^2}+\sum_{q=1}^{\ell-1}bt^{q-\ell}\theta^{q^2}+ \sum_{q>\ell}   b\bigl(\sqrt t\, bA\epsilon\bigr)^{2(q-\ell)}\theta^{q^2}+O(\epsilon^2)\biggr)(tbA\epsilon)^{\ell-1},
\end{aligned}
\end{equation}
where the $O(\epsilon^2)$ term collects the contribution of $O(\epsilon^2)$-terms in \eqref{E:6.60}. Now observe that, in the limit as~$\epsilon\downarrow0$, the term in the large parenthesis is bounded by the supremum in \eqref{E:6.79}, proving \eqref{E:6.76} for~$\ell\ge2$ once~$\epsilon$ is small enough.

For~$\ell=1$ we in turn invoke the first and last line in \eqref{E:6.60} to the leading order terms and bound the rest using the third line with the result
\begin{equation}
\begin{aligned}
 \sum_{q\ge1}&(tbA\epsilon)^{q-1}  \bigl|\partial_1 F_q(\blambda)\bigr|
\le 1-b(b-1)\lambda(1)^2 
\\
&+\bigl(b(b-1)+O(\epsilon^2)\bigr)\theta^4\lambda(1)(tbA\epsilon)+\sum_{q\ge3} \bigl(\sqrt t\, bA\epsilon\bigr)^{2(q-1)}\bigl(b\theta^{q^2} +  O(\epsilon^2)\bigr)
\end{aligned}
\end{equation}
Invoking the upper and lower bounds on~$\lambda(1)$, the right-hand side is bounded by
\begin{equation}
1-\Biggl(\frac{\frac12 b(b-1)}{\frac12\frac{(b-1)^2(b+1)^3}{b^3-1}+\delta}-tbA\theta^4\frac{b(b-1)}{\sqrt{\frac12\frac{(b-1)^2(b+1)^3}{b^3-1}-\delta}}\Biggr)\epsilon^2+O(\epsilon^3).
\end{equation}
We now check that the term in the large parentheses will be positive for $\theta$ close to~$1/b$ and~$\delta$ sufficiently small if \eqref{E:6.78} holds. For small-enough~$\epsilon$, the  term then dominates the expansion in powers of~$\epsilon$ which validates \eqref{E:6.76} for~$\ell=1$.
\end{proofsect}

\subsection{Convergence proofs}
We are now finally in a position to address the proofs of Theorems~\ref{thm-6.1} and~\ref{thm-3}. One last technical hurdle to get out of the way is the choice of parameters~$t$ and~$A$:

\begin{lemma}
\label{lemma-6.14a}
For all~$b\ge2$, the inequalities \twoeqref{E:6.78}{E:6.79} are true when
\begin{equation}
tA\le 2\sqrt b\quad\text{\rm and}\quad t\ge \frac75.
\end{equation}
\end{lemma}

\begin{proofsect}{Proof}
We start by proving that 
\begin{equation}
\inf_{\begin{subarray}{c}
b\in\N\\b\ge2
\end{subarray}}
\frac1{\sqrt2}\,\frac{\sqrt{1-b^{-3}}}{(1-b^{-2})\sqrt{1+b^{-1}}}\,b^{3/2}> 2.
\end{equation}
Indeed, for~$b\ge3$ we invoke $\sqrt{1-b^{-3}}>\sqrt{1-b^{-2}}$ and $\sqrt{1+b^{-1}}\le\sqrt2$ to dominate the expression by 
$b^{3/2}/2$  from below. Since~$b^{3/2}\ge 3^{3/2}\ge 5$ for~$b\ge3$, we are down to~$b=2$. Here we calculate the expression explicitly to be $\frac{4\sqrt7}{3\sqrt 3}$ which is above $2$, albeit just barely. It follows that \eqref{E:6.78} holds if~$tA\le 2\sqrt b$.

As for the second condition \eqref{E:6.79},  for~$\ell=2$ we need that $b^{-3}+t^{-1}<1$ which is true for all~$b\ge2$ whenever~$t >  8/7$.  For $\ell\ge3$ we bound the expression under the infimum by $b^{-8}+t^{-2}+b^{-3}\sum_{q \ge 1} t^{-q}$. Under the assumption that $t > 1$, the monotonicity in~$b$ shows that it suffices to have 
\begin{equation}
2^{-8} + t^{-2}+\frac18\frac1{t-1}<1
\end{equation}
The left-hand side is decreasing in~$t$ and, as is checked, is less than~$1$ at~$t:=7/5$. Hence \eqref{E:6.79} holds for all~$t\ge 7/5$. 
\end{proofsect}

We are now ready for:

\begin{proofsect}{Proof of Theorem~\ref{thm-6.1}}
We assume throughout that~$b\theta-1$ is positive and small enough so that the statements of above lemmas apply. As to the choice of~$t$ and~$A$, relying on Lemma~\ref{lemma-6.14a}, we set~$t:= 7/5$ and put~$A:=\frac{10}7\sqrt b$. Notice that this enables Lemma~\ref{lemma-6.12} and satisfies the requirements in Lemmas~\ref{lemma-6.12} and ~\ref{lemma-6.13} and where restrictions on~$A$ appeared. Also note that~$tA=2\sqrt b$ so $\varrho$ from \eqref{E:6.56} coincides with the expression in \eqref{E:6.6iu}.

Next observe that~$\varrho$ is a metric on~$\Sigma_0$ and, relying on product topology (and,  for~$b$ odd, also on the completeness of the space of probability measures on~$[0,1)$), that~$(\Sigma_0,\varrho)$ is complete. By Lemma~\ref{lemma-6.16},~$F$ is a strict contraction on~$\Sigma_0$. Using Lemma~\ref{lemma-6.12} along with the Banach contraction principle, iterations of~$F$ on any~$\blambda\in\Sigma'$ thus converge to some~$\blambda_\star\in\Sigma_0$ which is then also a unique fixed point of~$F$ in~$\Sigma'$.

Let us now consider a sequence  $\{\blambda_k\}_{k=0}^n$  obtained by $\blambda_k:=F^{(k)}(\blambda_{k-1})$ starting from some~$\blambda_{0}\in\Sigma'$. In order to control the approach of this sequence to~$\blambda_\star$, we need to compare the action of~$F$ and~$F^{(k)}$. For this we first note that, for all $\blambda\in\Sigma_0$ and all~$q\in\Z$,
\begin{equation}
\begin{aligned}
\bigl|F^{(k)}_q(\blambda)-F_q(\blambda)\bigr|
&=\frac{G_q(\blambda)}{G_0(\blambda)}\bigl|\theta^{q^2\sigma_k^2}-\theta^{q^2}\bigr|
\\
&\le (bA\epsilon)^q\theta^{q^2\min\{1,\sigma_k^2\}} \log(1/\theta) q^2 |\sigma_k^2-1|,
\end{aligned}
\end{equation}
where~$\epsilon:=\sqrt{b\theta-1}$. Assuming~$\sigma_k^2\ge1/2$ and recalling that~$A\le 2\sqrt b$,  we thus get
\begin{equation}
\label{E:6.69}
\begin{aligned}
\varrho\bigl(F^{(k)}(\blambda),F(\blambda)\bigr)
&=\sum_{q\ge1}(2b^{3/2}\epsilon)^{q-1} \bigl|F^{(k)}_q(\blambda)-F_q(\blambda)\bigr|
\\
&\le \biggl(\,\sum_{q\ge1} (2b^{3/2}\epsilon)^{2q-1} q^2 \theta^{q^2/2}  \log(1/\theta)\biggr)  |\sigma_k^2-1|. 
\end{aligned}
\end{equation}
Write~$\wh C$ for the expression in the parenthesis and recall the sequence $\{\frakd_k\}_{k\ge0}$ from Assumption~\ref{ass-1}. Abbreviate $\frakd_k':=\frakd_{\min\{k,n-k\}}$. For all $\blambda\in\Sigma_0$ the triangle inequality along with \eqref{E:6.74} and \eqref{E:6.69}  show
\begin{equation}
\begin{aligned}
\varrho\bigl(F^{(k)}(\blambda),\blambda_\star\bigr)
&\le \varrho\bigl(F^{(k)}(\blambda),F(\blambda)\bigr)+\varrho\bigl(F(\blambda),F(\blambda_\star)\bigr)
\\
&\le \wh C\,\frakd_k'+(1-\eta\epsilon^2\bigr)\varrho(\blambda,\blambda_\star).
\end{aligned}
\end{equation}
Using this for~$\blambda:=\blambda_{k-1}$ yields
\begin{equation}
\label{E:6.89}
\varrho(\blambda_k,\blambda_\star)
\le \wh C\,\frakd_k'+(1-\eta\epsilon^2\bigr)\varrho(\blambda_{k-1},\blambda_\star).
\end{equation}
whenever~$k$ is such that~$\blambda_k\in\Sigma_0$ and~$\sigma_k^2\ge1/2$.

To finish the proof, consider the family of $\{\sigma_k^2\}_{k=0}^n$ conforming to Assumption~\ref{ass-1} with sequence $\{\frakd_k\}_{k\ge0}$.  Denote $\diam(\Sigma_0):=\sup\{\varrho(\blambda,\blambda')\colon \blambda,\blambda'\in\Sigma_0\}$ and, for~$n\ge1$, let
\begin{equation}
k_0:=1+\sup_{n\ge1}\max\bigl\{k\le n/2\colon\{\blambda_k,\blambda_{n-k}\}\not\subseteq\Sigma_0\vee\frakd_k'>1/2\bigr\},
\end{equation}
where the maximum is set to be~$0$ if the set is empty. The above lemmas show that~$k_0<\infty$ for each~$\blambda_{0}\in\Sigma'$. Since~$\blambda_{k_0}\in\Sigma_0$, iterations of \eqref{E:6.89} then show
\begin{equation}
\varrho(\blambda_k,\blambda_\star)
\le\wh C\sum_{j=0}^{k-k_0-1}(1-\eta\epsilon^2)^j\frakd'_{k-j}
+(1-\eta\epsilon^2)^{k-k_0}\diam(\Sigma_0)
\end{equation}
whenever~$\min\{k,n-k\}\ge k_0$. To get the claim set $C:=\max\{\wh C,\diam(\Sigma_0)\}$, write~$\texte^{-\eta}$ instead of $1-\eta\epsilon^2$ and extend the range of the sum to all~$j\le k$.
\end{proofsect}

With this we now quickly finish also:

\begin{proofsect}{Proof of Theorem~\ref{thm-5}}
Let~$\epsilon>0$ be such that Theorem~\ref{thm-6.1} applies. This yields the existence of~$\blambda_\star$ which obeys~\eqref{E:3.28u} by extension of the bounds from Lemma~\ref{lemma-6.7a}. The bound \eqref{Eq:3.31} in turn follows from \eqref{E:6.6iu} by retaining only the term corresponding to~$q$ in the sum on the left and redefining~$C$ correspondingly. 

 Let~$\tilde v_\star$ be as defined in \eqref{E:tilde-v_star}.  For the convergence of~$v_k$ and its derivative, we need the uniform bound
\begin{equation}
\label{E:6.93}
\biggl|\frac{\texte^{-v_k(z)}}{a_k(0)}-\texte^{- \tilde v_\star(z)}\biggr|\le\sum_{q\in\Z}\bigl|\lambda_k(q)-\lambda_\star(q)\bigr|.
\end{equation}
The bounds \eqref{E:3.28u} and \eqref{Eq:3.31} then show that the sum on the right tends to zero as $\min\{k,n-k\}\to\infty$. 
Under the additional assumption that~$\{\frakd_k\}_{k\ge0}$ decays exponentially we can unite the estimates \twoeqref{Eq:3.30w}{Eq:3.31} as
\begin{equation}
\label{E:6.94}
\bigl|\lambda_k(q)-\lambda_\star(q)\bigr|\le C'\texte^{-\eta'\max\{k,|q|\}}.
\end{equation}
This now readily shows that the sum on the right of \eqref{E:6.93} decays exponentially with~$k$, proving
\begin{equation}
\sup_{z\in\R}\bigl|v_k(z)+\log a_k(0)-\tilde v_\star(z)\bigr|\le C'\texte^{-\eta'\min\{k,n-k\}}.
\end{equation}
To derive \eqref{Eq:3.32a} from this, note that a simple telescoping argument gives
\begin{equation}
\begin{aligned}
\Biggl|\frac{a_k(0)}{a_{k-1}(0)^b}\,-&\!\sum_{\begin{subarray}{c}
q_1,\dots,q_b\in\Z\\q_1+\dots+q_b=0
\end{subarray}}
\prod_{i=1}^b \lambda_\star(q_i)\,\Biggr|
\\
&\le\sum_{i=1}^b\Bigl(\,\sum_{q\in\Z}\lambda_{k-1}(q)\Bigr)^{i-1}\Bigl(\,\sum_{q\in\Z} \bigl|\lambda_{k-1}(q)-\lambda_\star(q)\bigr| \Bigr)\Bigl(\,\sum_{q\in\Z}\lambda_\star(q)\Bigr)^{b-i}.
\end{aligned}
\end{equation}
 The right-hand side is now bounded by $C''\texte^{-\eta'\min\{k,n-k\}}$. Since $v_\star(z)-b v_\star(z')$ differs from $\tilde v_\star(z)-b \tilde v_\star(z')$ by the logarithm of the giant sum on the left, \eqref{Eq:3.32a} follows by combining the previous two estimates.

To extend the convergence to the derivatives, we note that $v_\star-\tilde v_\star$ differ by a constant and so~$v_\star'=\tilde v_\star'$. Here we get
\begin{equation}
v_k'(z)=  - \texte^{v_k(z)}  a_k(0)\sum_{q\in\Z}\lambda_k(q)(2\pi\texti q)\texte^{2\pi\texti qz}
\end{equation}
and
\begin{equation}
v_\star'(z) =  - \texte^{\tilde v_\star(z)}  \sum_{q\in\Z}\lambda_\star(q)(2\pi\texti q)\texte^{2\pi\texti qz}
\end{equation}
This implies
\begin{equation}
\begin{aligned}
\bigl|v_k'(z)-v_\star'(z)\bigr|
&\le \texte^{v_k(z)} a_k(0)\texte^{\tilde v_\star(z)}\biggl|\frac{\texte^{-v_k(z)}}{a_k(0)}-\texte^{-\tilde v_\star(z)}\biggr|
\sum_{q\in\Z}\lambda_k(q)2\pi|q|
\\&\qquad+\texte^{\tilde v_\star(z)}\sum_{q\in\Z}\bigl|\lambda_k(q)-\lambda_\star(q)\bigr|(2\pi|q|)
\end{aligned}
\end{equation}
Since \eqref{Eq:3.30w} implies that  $\texte^{-v_k(z)}/a_k(0)\ge 1/2$ and  \eqref{E:3.28u} gives $\texte^{-\tilde v_\star(z)}\ge1/2$ once~$b\theta-1$ is sufficiently small while $\sum_{q\in\Z}\lambda_k(q)|q|$ is bounded uniformly in~$k$, both terms on the right decay to zero as~$\min\{k,n-k\}\to\infty$. The decay is exponentially fast if~$\{\frakd_j\}_{j\ge0}$ decays exponentially. This yields \eqref{Eq:3.32} as desired.
\end{proofsect}

\begin{remark}
\label{rem-large-b}
The above proofs are tailored for the near-critical regime, meaning with~$b$ fixed and~$b\theta-1$ positive but small. Another interesting asymptotic regime which can be analyzed is that of large~$b$. Focussing for simplicity on~$b$ even and~$\theta_k=\theta$, here Lemma~\ref{lemma-1.2} shows that~$\lambda_k(q)\le\theta^{q^2}$ for all~$k\ge0$ and all~$q\in\Z$. As~$\theta$ is close to~$1/b$, this suggests introduction of the scaled variables
\begin{equation}
\lambda_k'(q):=b^{q^2}\lambda_k(q).
\end{equation}
In these reduced variables the functions driving the iterations \eqref{E:1.5} take the form of ratios of Bessel functions; namely, 
\begin{equation}
\label{E:6.99}
\lambda_k'(1) = b\theta \frac{\sum_{\ell\ge0}\frac1{\ell!(\ell+1)!}\lambda_{k-1}'(1)^{2\ell+1}+O(b^{-1})}{ \sum_{\ell\ge0}  \frac1{\ell!\ell!}\lambda_{k-1}'(1)^{2\ell}+O(b^{-1})}
\end{equation}
and, for~$|q|\ge2$, 
\begin{equation}
\label{E:6.100}
\lambda_k'(q) = (b\theta)^{q^2} \frac{\sum_{\ell\ge0}\frac1{\ell!(\ell+|q|)!}\lambda_{k-1}'(1)^{2\ell+|q|}+O(b^{-1})}{ \sum_{\ell\ge0} \frac1{\ell!\ell!}\lambda_{k-1}'(1)^{2\ell}+O(b^{-1})}.
\end{equation}
Using these one can show that~$\lambda_k'(1)$ converges to a positive quantity characterized, up to errors that vanish as~$b\to\infty$, as a fixed point  of the ratio of two modified Bessel functions. Since only $\lambda_{k-1}'(1)$ appears explicitly on the right of \eqref{E:6.100}, the evolution of the other reduced variables is governed by the first variable up to errors that vanish as~$b\to\infty$. 

We have in fact carried our initial proof in this framework except that, in order to overcome the non-linearity of the right-hand side \eqref{E:6.99}, we ultimately also had to assume that~$b\theta$ is close to~$1$. However, we expect that with increasing~$b$ large, one should be able to control larger and larger intervals of~$b\theta$.
\end{remark}


\section*{Acknowledgments}
\nopagebreak\nopagebreak\noindent
This project has been supported in part by the NSF awards DMS-1954343 and DMS-2348113 and the BSF award no.~2018330. The authors are thankful to Christophe Garban and Jiwoon Park for helpful suggestions on the first version of the manuscript. M.B.\ wishes to dedicate this paper to the memory of Pierluigi Falco, a wonderful mathematician and a friend who departed this world way too young more than a decade ago.

\bibliographystyle{abbrv}

\end{document}